\documentclass[a4paper, british]{amsart}

\usepackage{libertine}
\usepackage[libertine]{newtxmath}

\usepackage{babel}
\usepackage{enumitem}
\usepackage{hyperref}
\usepackage[utf8]{inputenc}
\usepackage{newunicodechar}
\usepackage{mathtools}
\usepackage{varioref}
\usepackage[arrow,curve,matrix]{xy}

\usepackage{colortbl}
\usepackage{graphicx}
\usepackage{tikz}

\usepackage{mathrsfs}

\definecolor{linkred}{rgb}{0.7,0.2,0.2}
\definecolor{linkblue}{rgb}{0,0.2,0.6}

\numberwithin{figure}{section}

\usepackage[hyperpageref]{backref}

\setdescription{labelindent=\parindent, leftmargin=2\parindent}
\setitemize[1]{labelindent=\parindent, leftmargin=2\parindent}
\setenumerate[1]{labelindent=0cm, leftmargin=*, widest=iiii}

\newunicodechar{α}{\ensuremath{\alpha}}
\newunicodechar{β}{\ensuremath{\beta}}
\newunicodechar{χ}{\ensuremath{\chi}}
\newunicodechar{δ}{\ensuremath{\delta}}
\newunicodechar{ε}{\ensuremath{\varepsilon}}
\newunicodechar{Δ}{\ensuremath{\Delta}}
\newunicodechar{η}{\ensuremath{\eta}}
\newunicodechar{γ}{\ensuremath{\gamma}}
\newunicodechar{Γ}{\ensuremath{\Gamma}}
\newunicodechar{ι}{\ensuremath{\iota}}
\newunicodechar{κ}{\ensuremath{\kappa}}
\newunicodechar{λ}{\ensuremath{\lambda}}
\newunicodechar{Λ}{\ensuremath{\Lambda}}
\newunicodechar{ν}{\ensuremath{\nu}}
\newunicodechar{μ}{\ensuremath{\mu}}
\newunicodechar{ω}{\ensuremath{\omega}}
\newunicodechar{Ω}{\ensuremath{\Omega}}
\newunicodechar{π}{\ensuremath{\pi}}
\newunicodechar{Π}{\ensuremath{\Pi}}
\newunicodechar{φ}{\ensuremath{\phi}}
\newunicodechar{Φ}{\ensuremath{\Phi}}
\newunicodechar{ψ}{\ensuremath{\psi}}
\newunicodechar{Ψ}{\ensuremath{\Psi}}
\newunicodechar{ρ}{\ensuremath{\rho}}
\newunicodechar{σ}{\ensuremath{\sigma}}
\newunicodechar{Σ}{\ensuremath{\Sigma}}
\newunicodechar{τ}{\ensuremath{\tau}}
\newunicodechar{θ}{\ensuremath{\theta}}
\newunicodechar{Θ}{\ensuremath{\Theta}}
\newunicodechar{ξ}{\ensuremath{\xi}}
\newunicodechar{Ξ}{\ensuremath{\Xi}}
\newunicodechar{ζ}{\ensuremath{\zeta}}

\newunicodechar{ℓ}{\ensuremath{\ell}}
\newunicodechar{ï}{\"{\i}}

\newunicodechar{𝔸}{\ensuremath{\bA}}
\newunicodechar{𝔹}{\ensuremath{\bB}}
\newunicodechar{ℂ}{\ensuremath{\bC}}
\newunicodechar{𝔻}{\ensuremath{\bD}}
\newunicodechar{𝔼}{\ensuremath{\bE}}
\newunicodechar{𝔽}{\ensuremath{\bF}}
\newunicodechar{𝔾}{\ensuremath{\bG}}
\newunicodechar{ℕ}{\ensuremath{\bN}}
\newunicodechar{ℙ}{\ensuremath{\bP}}
\newunicodechar{ℚ}{\ensuremath{\bQ}}
\newunicodechar{ℝ}{\ensuremath{\bR}}
\newunicodechar{𝕏}{\ensuremath{\bX}}
\newunicodechar{ℤ}{\ensuremath{\bZ}}
\newunicodechar{𝒜}{\ensuremath{\sA}}
\newunicodechar{ℬ}{\ensuremath{\sB}}
\newunicodechar{𝒞}{\ensuremath{\sC}}
\newunicodechar{𝒟}{\ensuremath{\sD}}
\newunicodechar{ℰ}{\ensuremath{\sE}}
\newunicodechar{ℱ}{\ensuremath{\sF}}
\newunicodechar{𝒢}{\ensuremath{\sG}}
\newunicodechar{ℋ}{\ensuremath{\sH}}
\newunicodechar{𝒥}{\ensuremath{\sJ}}
\newunicodechar{ℒ}{\ensuremath{\sL}}
\newunicodechar{𝒪}{\ensuremath{\sO}}
\newunicodechar{𝒬}{\ensuremath{\sQ}}
\newunicodechar{𝒮}{\ensuremath{\sS}}
\newunicodechar{𝒯}{\ensuremath{\sT}}
\newunicodechar{𝒲}{\ensuremath{\sW}}

\newunicodechar{∂}{\ensuremath{\partial}}
\newunicodechar{∇}{\ensuremath{\nabla}}

\newunicodechar{↺}{\ensuremath{\circlearrowleft}}
\newunicodechar{∞}{\ensuremath{\infty}}
\newunicodechar{⊕}{\ensuremath{\oplus}}
\newunicodechar{⊗}{\ensuremath{\otimes}}
\newunicodechar{•}{\ensuremath{\bullet}}
\newunicodechar{Λ}{\ensuremath{\wedge}}
\newunicodechar{↪}{\ensuremath{\into}}
\newunicodechar{→}{\ensuremath{\to}}
\newunicodechar{↦}{\ensuremath{\mapsto}}
\newunicodechar{⨯}{\ensuremath{\times}}
\newunicodechar{∪}{\ensuremath{\cup}}
\newunicodechar{∩}{\ensuremath{\cap}}
\newunicodechar{⊋}{\ensuremath{\supsetneq}}
\newunicodechar{⊇}{\ensuremath{\supseteq}}
\newunicodechar{⊃}{\ensuremath{\supset}}
\newunicodechar{⊊}{\ensuremath{\subsetneq}}
\newunicodechar{⊆}{\ensuremath{\subseteq}}
\newunicodechar{⊂}{\ensuremath{\subset}}
\newunicodechar{⊄}{\ensuremath{\not \subset}}
\newunicodechar{≥}{\ensuremath{\geq}}
\newunicodechar{≠}{\ensuremath{\neq}}
\newunicodechar{≫}{\ensuremath{\gg}}
\newunicodechar{≪}{\ensuremath{\ll}}

\newunicodechar{≤}{\ensuremath{\leq}}
\newunicodechar{∈}{\ensuremath{\in}}
\newunicodechar{∉}{\ensuremath{\not \in}}
\newunicodechar{∖}{\ensuremath{\setminus}}
\newunicodechar{◦}{\ensuremath{\circ}}
\newunicodechar{°}{\ensuremath{^\circ}}
\newunicodechar{…}{\ifmmode\mathellipsis\else\textellipsis\fi}
\newunicodechar{·}{\ensuremath{\cdot}}
\newunicodechar{⋯}{\ensuremath{\cdots}}
\newunicodechar{∅}{\ensuremath{\emptyset}}
\newunicodechar{⇒}{\ensuremath{\Rightarrow}}

\newunicodechar{⁰}{\ensuremath{^0}}
\newunicodechar{¹}{\ensuremath{^1}}
\newunicodechar{²}{\ensuremath{^2}}
\newunicodechar{³}{\ensuremath{^3}}
\newunicodechar{⁴}{\ensuremath{^4}}
\newunicodechar{⁵}{\ensuremath{^5}}
\newunicodechar{⁶}{\ensuremath{^6}}
\newunicodechar{⁷}{\ensuremath{^7}}
\newunicodechar{⁸}{\ensuremath{^8}}
\newunicodechar{⁹}{\ensuremath{^9}}
\newunicodechar{ⁱ}{\ensuremath{^i}}

\newunicodechar{⌈}{\ensuremath{\lceil}}
\newunicodechar{⌉}{\ensuremath{\rceil}}
\newunicodechar{⌊}{\ensuremath{\lfloor}}
\newunicodechar{⌋}{\ensuremath{\rfloor}}

\newunicodechar{≅}{\ensuremath{\cong}}
\newunicodechar{⇔}{\ensuremath{\Leftrightarrow}}
\newunicodechar{∃}{\ensuremath{\exists}}
\newunicodechar{±}{\ensuremath{\pm}}

\DeclareFontFamily{OMS}{rsfs}{\skewchar\font'60}
\DeclareFontShape{OMS}{rsfs}{m}{n}{<-5>rsfs5 <5-7>rsfs7 <7->rsfs10 }{}
\DeclareSymbolFont{rsfs}{OMS}{rsfs}{m}{n}
\DeclareSymbolFontAlphabet{\scr}{rsfs}
\DeclareSymbolFontAlphabet{\scr}{rsfs}

\DeclareFontFamily{U}{mathx}{\hyphenchar\font45}
\DeclareFontShape{U}{mathx}{m}{n}{
      <5> <6> <7> <8> <9> <10>
      <10.95> <12> <14.4> <17.28> <20.74> <24.88>
      mathx10
      }{}
\DeclareSymbolFont{mathx}{U}{mathx}{m}{n}
\DeclareFontSubstitution{U}{mathx}{m}{n}
\DeclareMathAccent{\wcheck}{0}{mathx}{"71}

\newcommand{\sA}{\scr{A}}
\newcommand{\sB}{\scr{B}}
\newcommand{\sC}{\scr{C}}
\newcommand{\sD}{\scr{D}}
\newcommand{\sE}{\scr{E}}
\newcommand{\sF}{\scr{F}}
\newcommand{\sG}{\scr{G}}
\newcommand{\sH}{\scr{H}}

\newcommand{\sJ}{\scr{J}}

\newcommand{\sL}{\scr{L}}

\newcommand{\sO}{\scr{O}}

\newcommand{\sQ}{\scr{Q}}

\newcommand{\sS}{\scr{S}}
\newcommand{\sT}{\scr{T}}

\newcommand{\sW}{\scr{W}}

\newcommand{\bA}{\mathbb{A}}
\newcommand{\bB}{\mathbb{B}}
\newcommand{\bC}{\mathbb{C}}
\newcommand{\bD}{\mathbb{D}}
\newcommand{\bE}{\mathbb{E}}
\newcommand{\bF}{\mathbb{F}}
\newcommand{\bG}{\mathbb{G}}

\newcommand{\bN}{\mathbb{N}}

\newcommand{\bP}{\mathbb{P}}
\newcommand{\bQ}{\mathbb{Q}}
\newcommand{\bR}{\mathbb{R}}

\newcommand{\bX}{\mathbb{X}}

\newcommand{\bZ}{\mathbb{Z}}

\theoremstyle{plain}
\newtheorem{Th}{Theorem}[section]

\newtheorem{cor}[Th]{Corollary}
\newtheorem{defn}[Th]{Definition}

\newtheorem{lem}[Th]{Lemma}

\newtheorem{Prop}[Th]{Proposition}

\theoremstyle{remark}

\newtheorem{c-n-d}[Th]{Claim and Definition}

\newtheorem{rem}[Th]{\color{blue}{Remark}}

\newtheorem*{rem-nonumber}{Remark}

\numberwithin{equation}{Th}

\setlist[enumerate]{label=(\thethm.\arabic*), before={\setcounter{enumi}{\value{equation}}}, after={\setcounter{equation}{\value{enumi}}}}

\newcommand{\into}{\hookrightarrow}

\newcommand{\factor}[2]{\left. \raise 2pt\hbox{$#1$} \right/\hskip -2pt\raise -2pt\hbox{$#2$}}
\newcommand{\Publication}[1]{}

\newcommand{\svnid}[1]{}
\newcommand{\approvals}[2][Approval]{}
\renewcommand{\phi}{\varphi}
\graphicspath{{../}}

\usepackage{tikz-cd}
\tikzset{commutative diagrams/arrow style=tikz}

\author{Mohamed Kaddar} %
\address{} %
\email{\href{mailto:mohamed.kaddar@univ-lorraine.fr}{mohamed.kaddar@univ-lorraine.fr}} %


\title[Universally equidimensional morphisms and weakly or strongly rational singularities.]{Universally equidimensional morphisms and weakly or strongly rational singularities.}%
\date{\today}

\makeatletter
\hypersetup{
  pdfauthor={\authors},
  pdftitle={\@title},
  pdfsubject={\@subjclass},
  pdfkeywords={\@keywords},
  pdfstartview={Fit},
  pdfpagelayout={TwoColumnRight},
  pdfpagemode={UseOutlines},
  bookmarks,
  colorlinks,
  linkcolor=linkblue,
  citecolor=linkred,
  urlcolor=linkred
}
\makeatother

\newcommand{\chapref}[1]{\hyperref[#1]{Chapter~\ref*{#1}}}
\newcommand{\lemmaref}[1]{\hyperref[#1]{Lemma~\ref*{#1}}}
\newcommand{\parref}[1]{\hyperref[#1]{Section~\ref*{#1}}}
\newcommand{\theoremref}[1]{\hyperref[#1]{Theorem~\ref*{#1}}}
\newcommand{\definitionref}[1]{\hyperref[#1]{Definition~\ref*{#1}}}
\newcommand{\propositionref}[1]{\hyperref[#1]{Proposition~\ref*{#1}}}
\newcommand{\conjectureref}[1]{\hyperref[#1]{Conjecture~\ref*{#1}}}
\newcommand{\corollaryref}[1]{\hyperref[#1]{Corollary~\ref*{#1}}}
\newcommand{\exampleref}[1]{\hyperref[#1]{Example~\ref*{#1}}}
\newcommand{\exerciseref}[1]{\hyperref[#1]{Exercise~\ref*{#1}}}
\newcommand{\factref}[1]{\hyperref[#1]{Fact~\ref*{#1}}}
\newcommand{\claimref}[1]{\hyperref[#1]{Claim~\ref*{#1}}}
\newcommand{\remarkref}[1]{\hyperref[#1]{Remark~\ref*{#1}}}
\newcommand{\settingref}[1]{\hyperref[#1]{Setting~\ref*{#1}}}
\newcommand{\appendixref}[1]{\hyperref[#1]{Appendix~\ref*{#1}}}

\theoremstyle{plain}

\theoremstyle{remark}

\newcommand{\eq}[1][r]
   {\ar@<-3pt>@{-}[#1]
    \ar@<-1pt>@{}[#1]|<{}="gauche"
    \ar@<+0pt>@{}[#1]|-{}="milieu"
    \ar@<+1pt>@{}[#1]|>{}="droite"
    \ar@/^2pt/@{-}"gauche";"milieu"
    \ar@/_2pt/@{-}"milieu";"droite"}

\begin{document}
\maketitle
\approvals[Approval for Abstract]{Mohamed & yes}
\begin{abstract}

This article aims to study  the behavior of certain types of singularities in a universally equidimensional morphism (i.e., open with constant pure dimensional fibers ). These singularities are those of reduced complex spaces of pure dimension $m$ for which the sheaf ${\mathcal L}^{m}_{Z}$ (whose sections are meromorphic forms that extend analytically over any desingularization of $Z$) has depth at least two and are called {\emph{weakly rational}} and denoted ${\mathfrak{F.R}}$; for spaces whose singular locus is of codimension at least two, they are called {\emph{strongly rational}} and denoted ${\mathfrak{F.R}}^{*}$. Our study focuses on the possibility of transferring this type of singularities from the total space to the base, from the base and the fibers to the total space, and from the latter to the fibers.

\end{abstract}

\tableofcontents
%
%

\noindent
\phantomsection\addcontentsline{toc}{part}{Introduction}

The present work aims to study certain types of singularities in families of subspaces described by universally equidimensional morphisms (i.e., open with constant fiber dimension) using the sheaves $\omega^{\bullet}_{Z}$ and ${\mathcal L}^{\bullet}_{Z}$ studied in \cite{K3} and \cite{K4} (Recall that in the top degree, ${\mathcal L}^{m}_{Z}$ is the sheaf of meromorphic forms whose pullback extends analytically over any desingularization, and $\omega^{m}_{Z}$ is Grothendieck's dualizing sheaf).

It has been known since R. Elkik's article \cite{Elkik78} that any deformation of a rational singularity over a base with rational singularities remains rational; stated in the algebraic framework (schemes of finite type over a field of characteristic zero), it remains valid in the complex analytic setting. Since the problem is local on each fiber and on the total space, the rationality of fiber singularities requires both the  morphism to be flat and the  source to be Cohen-Macaulay. In this construction, flatness plays a crucial role as it allows, on one hand, to lift any non-zero divisor to a non-zero divisor at the level of local rings, and on the other hand, to reveal a non-trivial relationship between the absolute and relative dualizing complexes. This leads us to primarily consider the case of a family of divisors parameterized by a smooth base, and mainly by a smooth curve. The general construction is achieved by induction on the base dimension and by applying Hironaka's desingularization theorem \cite{Hironaka} to reach the case of a singular base.

As in the complex analytic setting, it is rare for a geometric situation without special conditions to exhibit Cohen-Macaulay structures, so it seemed natural to try to free ourselves from this constraint while maintaining a certain type of singularity that is, in some sense, close to rational ones. It is in this spirit that da Silva and Andréatta study in \cite{AS84} the deformation of normal spaces whose sheaves $\omega^{\bullet}_{Z}$ and ${\mathcal L}^{\bullet}_{Z}$ coincide in maximal degree which they call \emph{weakly rational} singularities. More precisely, a normal complex space $X$ of pure dimension $m$ has a weakly rational singularity at a point $x$ if, for a desingularization $\phi:{\tilde X}\rightarrow X$, we have the vanishing $({\rm I}\!{\rm R}^{m-1}{\phi}_{*}{\mathcal O}_{\tilde X})_{x} =0$. They show, following the method developed in \cite{Elkik78}, that any deformation of a weakly rational singularity over a \emph{smooth} base remains weakly rational.

This generalization remains partial as it confines them to the case of a smooth base without hope of addressing the singular case through Hironaka's desingularization theorem \cite{Hironaka} (the arguments developed in \cite{Elkik78} are unsuitable for the weakly rational case).

In the following work and unless explicitly stated otherwise, we only consider reduced complex spaces of pure dimension without any normality condition. If the space is not normal or more generally if its singular locus is not of codimension greater that two, the cohomological vanishing condition above obviously does not imply the equality of the sheaves $\omega^{\bullet}_{Z}$ and ${\mathcal L}^{\bullet}_{Z}$ in maximal degree and consequently, these are two disntinct conditions. However, this is the case if, in addition to this condition, we have an identification of these two sheaves outside codimension at least two. These remarks lead us to denote by $\mathfrak{F.R}$ the class of singularities of {\emph{weakly rational}} type (satisfying only the cohomological condition) and by $\mathfrak{F.R}^{*}$ the class that we will call {\emph{strongly rational}}, that is, the class of reduced spaces of type $\mathfrak{F.R}$ with a singular locus of codimension at least two (note that the latter corresponds to the weakly rational spaces of \cite{AS84} if the space is normal). We provide some properties and characterizations of these classes in the preliminary section (for example, $\mathfrak{F.R}$ corresponds precisely to reduced complex spaces $Z$ of pure dimension $m$ on which the sheaf ${\mathcal L}^{m}_{Z}$ has depth at least two). \vspace{1mm}

\noindent
The main objective of this work is to generalize \cite{AS84} in three possible directions, namely:\par
  $\bullet$ weakening the { {flatness}} condition by considering universally equidimensional morphisms (open morphisms with constant fiber dimension),\par
   $\bullet$ considering singularities of type $\mathfrak{F.R}$ or $\mathfrak{F.R}^{*}$ (cf. \definitionref{def2} and \definitionref{def3}),\par
   $\bullet$ allowing singular bases with these types of singularities.\vspace{1mm}

\noindent Leading to one of the main results stating that:\vspace{1mm}

\noindent {\emph{If the fibers of a universally equidimensional morphism have singularities of type $\mathfrak{F.R}^*$ (resp. $\mathfrak{F.R}$) over a base with singularities of type $\mathfrak{F.R}^*$ (resp. $\mathfrak{F.R}$), then the total space (or source) is also of type $\mathfrak{F.R}^*$ (resp. $\mathfrak{F.R}$)}}.\vspace{2mm}

\noindent 
More precisely, we show:
    \Th{}{}\label{Th1}{\it\color{blue}{Stability under direct image for the classes $\mathfrak{R}$, $\mathfrak{F.R}$ and $\mathfrak{F.R}$}:} 

Let $\pi:Z\rightarrow S$ be a $n$-geometrically flat morphism of reduced complex spaces of pure dimensions $m$ and $r$ respectively. If $Z$ has singularities of rational, weakly rational, or strongly rational type, then $S$ has singularities of the same type.

\Th{}{}\label{Th2} {\it \color{blue}{Universally equidimensional deformation of the class $\mathfrak{F.R}^{*}$}:}

Let $\pi:Z\rightarrow S$ be a universally-$n$-equidimensional morphism of reduced complex spaces of pure dimensions $m$ and $r$ respectively. If $\pi$ is of type $\mathfrak{F.R}^{*}$, then so is $Z$.

\Th{}{}\label{Th3} {\it \color{blue}{Universally equidimensional and reduced deformation of the class $\mathfrak{F.R}$}:}\vspace{1mm}

\noindent Let $\pi:Z\rightarrow S$ be a universally-$n$-equidimensional morphism of reduced complex spaces of pure dimension whose base and fibers are of type $\mathfrak{F.R}$. Then, $Z$ is of type $\mathfrak{F.R}$, if one of the following two conditions holds: \vspace{1mm}

\noindent {\bf(i)} $Z$ is locally of pure dimension $m$ and $\pi$ has reduced fibers, or \vspace{1mm}

\noindent {\bf(ii)} if ${\rm N}{\rm N}(Z)$, the non-normal locus of $Z$, satisfies the incidence condition: 

$$(P)\,\,\,\,\,\widering{\pi^{-1}(s)\cap{\rm N}{\rm N}(Z)}=\emptyset\,\,{\rm in}\,\pi^{-1}(s),\,\,\forall\,s\in S$$
\rm
\rm
\Th{}{}\label{Th4} {\it \color{blue}{A weak converse}:}

Let $\pi:Z\rightarrow S$ be a universally-$n$-equidimensional morphism of reduced complex spaces of pure dimension with $S$ of type (resp. $\mathfrak{F.R}$) admitting a simultaneous resolution. If $Z$ is of type $\mathfrak{F.R}^{*}$ (resp. $\mathfrak{F.R}$) then so is $\pi$ is.
\vfill\eject
\phantomsection\addcontentsline{toc}{part}{Préliminaires.}
%
%

\section{\color{blue} {Notations et définitions.}}

\approvals{Mohamed & yes}
\defn{}{}\label{def1} A complex analytic space $Z$ of dimension $m$ has a \emph{rational singularity} at a point $z$ if for a desingularization $\phi:{\tilde Z}\rightarrow Z$, we have ${\rm I}\!{\rm R}{\phi}_{*}{\mathcal O}_{\tilde Z} \simeq {\mathcal O}_{Z}$ (i.e., $({\rm I}\!{\rm R}^{j}{\phi}_{*}{\mathcal O}_{\tilde Z})_{z}=0$ for all $j>0$).

\defn{}{}\label{def2} A complex analytic space $Z$ of dimension $m$ has a \emph{weakly rational singularity} at $z$ if for a desingularization $\phi:{\tilde Z}\rightarrow Z$, we have $({\rm I}\!{\rm R}^{m-1}{\phi}_{*}{\mathcal O}_{\tilde Z})_{z} =0$.

\defn{}{}\label{def3} A complex analytic space $Z$ of dimension $m$ has a \emph{strongly rational singularity} at $z$ if the canonical injective morphism ${\mathcal L}^{m}_{Z,z}\rightarrow \omega^{m}_{Z,z}$ is surjective.\rm

In the category of reduced complex analytic spaces, we will denote by ${\mathfrak {R}}$ (resp. ${\mathfrak {F.R}}$, resp. ${\mathfrak {F.R}}^{*}$) the class of rational singularities (resp. weakly rational, resp. strongly rational) and we will say that a complex space $Z$ is of type ${\mathfrak {R}}$ (resp. ${\mathfrak {F.R}}$, resp. ${\mathfrak {F.R}}^{*}$) if it has rational singularities (resp. weakly rational, resp. strongly rational).

\defn{}{}{} We will say that a complex space is of type $k-\mathfrak{{F.R}}$ if the sheaves $\mathcal{L}^{j}_{Z}$ have depth at least two for all integers $j\leq k$. It will be of type $k-\mathfrak{{F.R}}^{*}$ if, moreover, its singular locus has codimension at least two (in this case, we have the equality $\mathcal{L}^{j}_{Z}=\omega^{j}_{Z}$ for all $j\leq k$).

\defn{}{}\label{def4} We will say that a morphism $\pi:Z\rightarrow S$ of reduced complex spaces is \emph{reduced}, \emph{normal} or of type ${\mathfrak{F.R}}$ (resp. ${\mathfrak{F.R}}^{*}$) if $S$ and its fibers are reduced, normal or of type ${\mathfrak{F.R}}$ (resp. ${\mathfrak{F.R}}^{*}$). \rm

We say it is generically reduced or that the generic fiber is reduced if it is reduced over a dense subset of $S$.

\defn{}{}\label{def6} Let $n\in {\Bbb N}$. We will say that a morphism $\pi:Z\rightarrow S$ of complex spaces is universally-$n$-equidimensional if it is open with constant fiber dimension $n$. It is geometrically flat if these fibers can be equipped with suitable multiplicities to form an analytic family of $n$-cycles (cf. \cite{BM-1}, \cite{K1}, \cite{K2}).\rm\vspace{1mm}

\noindent Geometric flatness imposes the following property on the morphism:\vspace{1mm}

\noindent
For every point $z$ of $Z$ (naturally in the fiber $\pi^{-1}(\pi(z))$) and every local factorization at $z$ (relative to the fiber $\pi^{-1}(\pi(z))$,  
$$\xymatrix{Z\ar[r]^{f}\ar[rd]_{\pi}&S\times U\ar[d]^{q}\\  
&S}$$  
where $f$ is a finite, surjective, and open morphism, $q$ is the canonical projection, and $U$ is a relatively compact Stein open subset of some numerical space ${\Bbb C}^{n}$, there exists a trace morphism $f_{*}{\mathcal O}_{Z}\rightarrow {\mathcal O}_{S\times U}$ inducing a morphism (also of trace type) $f_{*}f^{*}\Omega^{n}_{S\times U/S}\rightarrow\Omega^{n}_{S\times U/S}$ that extends naturally to a trace morphism  
$$f_{*}\Omega^{n}_{Z/S}\rightarrow\Omega^{n}_{S\times U/S}.$$  
\vspace{1mm}

\noindent
\section{\color{blue} {Some properties of ${\mathfrak{R}}$, ${\mathfrak{F.R}}$ and ${\mathfrak{F.R}}^{*}$ singularities }.}
\approvals{Mohamed & yes}
\vspace {2mm}

\noindent
\subsection{Singularities of ${\mathfrak{R}}$-type.}
\vspace{2mm}

\noindent
One cannot cite all the existing literature on rational singularities, but it is impossible not to recall the now-classic result known as the {\it Kempf criterion} \cite{Kem}, valid in all characteristics:
\Prop{}{}\label{P1}: Let $Z$ be a reduced complex space. Then, the following equivalences hold:
\par\noindent{\bf(i)} $Z$ is of type ${\mathfrak R}$\vspace{1mm}

\noindent
{\bf(ii)} The canonical morphism ${\mathcal O}_{Z}\rightarrow{\rm I}\! {\rm R}{\mathcal H}om( {\mathcal L}^{m}_{Z}, {\mathcal D}^{\bullet}_{Z}[-m])$ is a quasi-isomorphism.\vspace{1mm}

\noindent{\bf(iii)} $Z$ is of type ${\mathfrak {F.R}}^{*}$ and Cohen-Macaulay.\vspace{1mm}

\noindent 
{\bf(iv)} The sheaf ${\mathcal L}^{m}_{Z}$ is maximal Cohen-Macaulay and $Z$ is normal.
\rm\par\noindent
 
\begin{proof} 
The equivalence between {\bf(i)} and {\bf(iii)} is precisely Kempf's criterion.\vspace{1mm}

\indent $\bullet$ {\bf(i)}$\Longleftrightarrow${\bf(ii)}: This is a straightforward application of analytic duality for a proper morphism and the Grauert-Riemenschneider vanishing theorem. Indeed, if $\phi:\tilde{Z}\rightarrow Z$ is a given desingularization, the relative analytic duality for a proper morphism from \cite{RRV} yields the isomorphism 
$${\rm I}\!{\rm R}\phi_{*}{\mathcal O}_{\widetilde{Z}}\simeq{\rm I}\!{\rm R}\phi_{*}{\rm I}\!{\rm R}{\mathcal H}om( \Omega^{m}_{\widetilde{Z}},  \Omega^{m}_{\widetilde{Z}})\simeq {\rm I}\!{\rm R}{\mathcal H}om({\rm I}\!{\rm R}\phi_{*}\Omega^{m}_{\widetilde{Z}}[m], {\mathcal D}^{\bullet}_{Z }) $$
However, the Grauert-Riemenschneider vanishing theorem \cite{GR70} gives the quasi-isomorphism ${\rm I}\!{\rm R}\phi_{*}\Omega^{m}_{\widetilde{Z}}={\mathcal L}^{m}_{Z}$ and, consequently,  
$${\rm I}\!{\rm R}\phi_{*}{\mathcal O}_{\widetilde{Z}}\simeq{\rm I}\!{\rm R}{\mathcal H}om({\mathcal L}^{m}_{Z}, {\mathcal D}^{\bullet}_{Z}[-m])$$
{\bf(ii)}$\Longleftrightarrow${\bf(iii)}: 
\vspace{1mm}

\indent $\bullet$ For the implication, note that by biduality, condition {\bf(ii)} is equivalent to the quasi-isomorphism ${\mathcal L}^{m}_{Z}\simeq{\mathcal D}^{\bullet}_{Z}[-m]$. Consequently, the dualizing complex reduces to a single sheaf in degree $-m$. Now, by construction, $\omega^{m}_{Z}={\mathcal H}^{-m}({\mathcal 
D}^{\bullet}_{{Z}})$ and thus ${\mathcal L}^{m}_{Z}=\omega^{m}_{Z}$. Therefore, $Z$ is Cohen-Macaulay and of type ${\mathfrak {F.R}}^{*}$.   \vspace{1mm} 

\noindent To ensure the normality of $Z$, observe that the canonical morphism
$${\rm I}\!{\rm R}{\mathcal H}om({\mathcal L}^{m}_{Z}, \omega^{m}_{Z})\rightarrow{\rm I}\!{\rm R}{\mathcal H}om({\mathcal L}^{m}_{Z}, {\mathcal D}^{\bullet}_{Z}[-m])$$
induces, at the level of degree $0$ cohomology, the morphism
${\mathcal H}om({\mathcal L}^{m}_{Z}, \omega^{m}_{Z})\rightarrow {\mathcal O}_{Z}$ 
which is necessarily injective since it is bijective between torsion-free sheaves. Since there is already a natural morphism ${\mathcal O}_{Z}\rightarrow {\mathcal H}om({\mathcal L}^{m}_{Z}, \omega^{m}_{Z})$ 
that is also injective for the same reasons, we deduce the isomorphism $${\mathcal O}_{Z}\simeq {\mathcal H}om({\mathcal L}^{m}_{Z}, \omega^{m}_{Z}) $$
From this, considering the natural injections
$${\mathcal O}_{Z}\subset {\mathcal L}^{0}_{Z}\subset{\mathcal H}om(
{\mathcal L}^{m}_{Z}, {\mathcal L}^{m}_{Z})\subset {\mathcal H}om({\mathcal L}^{m}_{Z}, \omega^{m}_{Z})$$ 
the isomorphism ${\mathcal O}_{Z}\simeq{\mathcal L}^{0}_{Z}$ characterizes the normality of $Z$.\vspace{1mm}

\indent $\bullet$ The converse is evident since, $Z$ being Cohen-Macaulay, we have 
$${\mathcal O}_{Z}\simeq{\rm I}\! {\rm R}{\mathcal H}om({\omega}^{m}_{Z}[m], {\omega}^{m}_{Z}[m])\simeq{\rm I}\! {\rm R}{\mathcal H}om({\mathcal L}^{m}_{Z}, {\mathcal D}^{\bullet}_{Z}[-m])$$  
{\bf(iii)}$\Longleftrightarrow${\bf(iv)}:\vspace{1mm}

\indent $\bullet$ The quasi-isomorphism ${\mathcal O}_{Z}\simeq{\rm I}\!{\rm R}{\mathcal H}om({\mathcal L}^{m}_{Z}, {\mathcal D}^{\bullet}_{Z}[-m])$ precisely means that 
$${\mathcal H}^{j}({\rm I}\!{\rm R}{\mathcal H}om({\mathcal L}^{m}_{Z}, {\mathcal D}^{\bullet}_{Z}[-m])=0\,\forall j\not=0\,{\rm and}\, {\mathcal H}^{0}({\rm I}\!{\rm R}{\mathcal H}om({\mathcal L}^{m}_{Z}, {\mathcal D}^{\bullet}_{Z}[-m])\simeq {\mathcal O}_{Z}$$
As seen above, the second condition implies the normality of $Z$. As for the first, it ensures that for every Stein open subset $U$ of $Z$, by the duality of \cite{RR70} or \cite{AK}, we have 
$${\rm Ext}^{j}(U; {\mathcal L}^{m}_{{Z}}, {\mathcal D}^{\bullet}_{Z }[-m])\simeq \big({\rm H}^{m-j}_{c}(U, {\mathcal L}^{m}_{{Z}})\big)^{'} $$
thanks to the duality of \cite{RR70} or \cite{AK}. However, the depth characterization from \cite{Ba1} guarantees that the vanishing conditions $${\rm H}^{k}_{c}(U, {\mathcal L}^{m}_{{Z}})=0,\,\forall\,k\not=m$$ 
are equivalent to ${\mathcal L}^{m}_{{Z}}$ being maximal Cohen-Macaulay.\vspace{1mm}

\indent $\bullet$ Conversely, if ${\mathcal L}^{m}_{{Z}}$ is maximal Cohen-Macaulay, then by the duality mentioned earlier, we have the vanishing conditions ${\mathcal H}^{j}({\rm I}\!{\rm R}{\mathcal H}om({\mathcal L}^{m}_{Z}, {\mathcal D}^{\bullet}_{Z}[-m]))=0\,\forall j\not=0$, and thus $${\rm I}\!{\rm R}{\mathcal H}om({\mathcal L}^{m}_{Z}, {\mathcal D}^{\bullet}_{Z}[-m])\simeq {\mathcal H}^{0}({\rm I}\!{\rm R}{\mathcal H}om({\mathcal L}^{m}_{Z}, {\mathcal D}^{\bullet}_{Z}[-m]))$$
It only remains to prove that ${\mathcal H}^{0}({\rm I}\!{\rm R}{\mathcal H}om({\mathcal L}^{m}_{Z}, {\mathcal D}^{\bullet}_{Z}[-m])\simeq {\mathcal O}_{Z}$. \vspace{1mm}

\noindent 
Since $${\mathcal H}^{j}({\mathcal L}^{m}_{Z})={\mathcal H}^{-j}({\mathcal D}^{\bullet}_{Z}[-m])=0,\,\forall\,j>0$$
there is a canonical morphism
$${\mathcal H}^{0}({\rm I}\!{\rm R}{\mathcal H}om({\mathcal L}^{m}_{Z}, {\mathcal D}^{\bullet}_{Z}[-m]))\rightarrow {\mathcal H}om({\mathcal H}^{0}({\mathcal L}^{m}_{Z}), {\mathcal H}^{0}({\mathcal D}^{\bullet}_{Z}[-m]))={\mathcal H}om({\mathcal L}^{m}_{Z}, \omega^{m}_{Z})$$
noting that there is an equally natural arrow in the other direction, as mentioned in the proof of the equivalence between {\bf(ii)} and {\bf(iii)}. The arguments presented in that proof show that we have an isomorphism
 $$ {\mathcal H}^{0}({\rm I}\!{\rm R}{\mathcal H}om({\mathcal L}^{m}_{Z}, {\mathcal D}^{\bullet}_{Z}[-m]))\simeq {\mathcal H}om({\mathcal L}^{m}_{Z}, \omega^{m}_{Z})$$   
But $Z$ being normal, the natural morphism 
${\mathcal O}_{Z}\rightarrow{\mathcal H}om({\mathcal L}^{m}_{Z}, \omega^{m}_{Z})$ is necessarily bijective since it is so on the regular part and the sheaves have depth at least two. Hence, the isomorphism
${\mathcal H}^{0}({\rm I}\!{\rm R}{\mathcal H}om({\mathcal L}^{m}_{Z}, {\mathcal D}^{\bullet}_{Z}[-m]))\simeq{\mathcal O}_{Z}$ and, consequently, the quasi-isomorphism
$${\mathcal O}_{Z}\simeq{\rm I}\!{\rm R}{\mathcal H}om({\mathcal L}^{m}_{Z}, {\mathcal D}^{\bullet}_{Z}[-m])\,\,\blacksquare$$
\end{proof}
\vspace{1mm}

\noindent
Note in passing that we have implicitly established the
\cor{}{}\label{C1} With the notations and hypotheses of \ref{P1}, the following assertions are equivalent:\vspace{1mm}

\noindent 
{\bf(i)} $Z$ has semi-rational singularities,\par\noindent
{\bf(ii)} the canonical morphism ${\mathcal L}^{0}_{Z}\rightarrow {\rm I}\!{\rm R}{\mathcal H}om({\mathcal L}^{m}_{Z}, {\mathcal D}^{\bullet}_{Z}[-m])$ is a quasi-isomorphism,\par\noindent
{\bf(iii)} the normalization of $Z$ has rational singularities,\par\noindent
{\bf(iv)} ${\mathcal L}^{m}_{Z}$ is a maximal Cohen-Macaulay sheaf.\rm

\begin{proof} {\bf(i)}$\Longleftrightarrow${\bf(ii)}: This is essentially due to the definition. Indeed, if $\phi:\widetilde{Z}\rightarrow Z$ is a given desingularization of $Z$, the semi-rationality of the singularities of $Z$ is equivalent to the vanishing ${\rm I}\!{\rm R}^{j}\phi_{*}{\mathcal O}_{\widetilde{Z}}=0$ for all non-zero integers $j$ and, therefore, to the vanishing of cohomology ${\mathcal H}^{j}({\rm I}\!{\rm R}{\mathcal H}om({\mathcal L}^{m}_{Z}, {\mathcal D}^{\bullet}_{Z}[-m]))=0$ for all non-zero $j$ since ${\rm I}\!{\rm R}\phi_{*}{\mathcal O}_{\widetilde{Z}}\simeq {\rm I}\!{\rm R}{\mathcal H}om({\mathcal L}^{m}_{Z}, {\mathcal D}^{\bullet}_{Z}[-m])$, where the degree $0$ cohomology corresponds precisely to the sheaf ${\mathcal L}^{0}_{Z}$. 
 \vspace{1mm}

 \noindent 
{\bf(i)}$\Longleftrightarrow${\bf(iii)}: It suffices to consider a commutative diagram of the form
$$\xymatrix{&\widetilde{Z}\ar[ld]_{\tilde\phi}\ar[rd]^{\phi}&\\
\overline{Z}\ar[rr]_{\nu}&&Z}$$
where $\nu$ is the normalization of $Z$, and $\phi$ is a desingularization dominating this normalization. Then, by the finiteness of $\nu$, we have 
 $${\rm I}\!{\rm R}^{j}\phi_{*}{\mathcal O}_{\widetilde{Z}}=0, \forall\,j\not=0 \,\Longleftrightarrow {\rm I}\!{\rm R}^{j}{\tilde\phi}_{*}{\mathcal O}_{\widetilde{Z}}=0, \forall\,j\not=0$$
 using the fact that for any finite morphism $f:X\rightarrow Y$ and any coherent sheaf ${\mathcal F}$ on $X$, $f_{*}{\mathcal F}=0\Longrightarrow\,{\mathcal F}=0$ due to the canonical surjection
 $f^*f_{*}{\mathcal F}\rightarrow {\mathcal F}$. \par\noindent
 Since $\overline Z$ is normal, ${\tilde\phi}_{*}{\mathcal O}_{\widetilde Z}={\mathcal O}_{\overline Z}$ and, thus, ${\rm I}\!{\rm R}{\tilde\phi}_{*}{\mathcal O}_{\widetilde{Z}}\simeq {\mathcal O}_{\overline Z}$.\vspace{1mm}
 
 \noindent
 To show that {\bf(iv)} is equivalent to any of the other points, we use classical duality theorems. Since our assertions are local in nature, we may assume $Z$ is Stein. Consider the Leray spectral sequences $${\rm E}^{i,j}_{2}:={\rm H}^{i}(Z, {\rm I}\!{\rm R}^{j}\phi_{*}{\mathcal O}_{\widetilde{Z}})\Longrightarrow {\rm H}^{i+j}(\widetilde{Z}, {\mathcal O}_{\widetilde{Z}}) $$
$${\rm E'}^{i,j}_{2}:={\rm H}^{i}_{c}(Z, {\rm I}\!{\rm R}^{j}\phi_{*}{\Omega}^{m}_{\widetilde{Z}})\Longrightarrow {\rm H}^{i+j}_{c}(\widetilde{Z}, {\Omega}^{m}_{\widetilde{Z}}) $$
which, as is known, degenerate by virtue of Grauert's coherence theorem for higher direct images under a proper morphism (\cite{GR0}) and Cartan's Theorem B (\cite{C}), for the first, and by \cite{GR70} for the second. \vspace{1mm}

\indent 
 Then, the topological duality isomorphisms (\cite{RR70}, \cite{Se})
  $${\rm Ext}^{j}(Z; {\mathcal L}^{m}_{{Z}}, {\mathcal D}^{\bullet}_{Z }[-m])\simeq \big({\rm H}^{m-j}_{c}(Z, {\mathcal L}^{m}_{{Z}})\big)^{'} $$ 
  $${\rm H}^{j}(\widetilde{Z}, {\mathcal O}_{\widetilde{Z}})\simeq \big({\rm H}^{m-j}_{c}(\widetilde{Z}, {\Omega}^{m}_{\widetilde{Z}})\big)^{'}$$
  combined with the topological isomorphisms given by the degeneration of the spectral sequences
$${\rm H}^{0}(Z, {\rm I}\!{\rm R}^{j}\phi_{*}{\mathcal O}_{\widetilde{Z}})\simeq {\rm H}^{j}(\widetilde{Z}, {\mathcal O}_{\widetilde{Z}}),\,\,\,{\rm H}^{i}_{c}(Z, {\mathcal L}^{m}_{Z})\simeq{\rm H}^{i}_{c}(\widetilde{Z}, {\Omega}^{m}_{\widetilde{Z}})$$
  show that, for all non-zero integers $j$,  
$${\rm H}^{m-j}_{c}(Z, {\mathcal L}^{m}_{{Z}})=0\,\Longleftrightarrow {\mathcal H}^{j}({\rm I}\!{\rm R}{\mathcal H}om({\mathcal L}^{m}_{Z}, {\mathcal D}^{\bullet}_{Z}[-m]))=0\, \Longleftrightarrow {\rm H}^{j}(\widetilde{Z}, {\mathcal O}_{\widetilde{Z}})=0\, \Longleftrightarrow {\rm I}\!{\rm R}^{j}\phi_{*}{\mathcal O}_{\widetilde{Z}}=0\,\,\blacksquare$$
\end{proof}
\rm
\vspace{1mm}  

\noindent  
A more or less independent result is given by the
\cor{}{}\label{C2} Let $Z$ be a complex space. Then, we have \vspace{1mm}

\noindent 
{\bf(i)} $Z$ is normal if and only if $\,{\rm Codim}({\rm Sing}(Z))\geq 2$ and 
$\displaystyle{{\mathcal O}_{Z}\simeq{\mathcal H}om(\omega^{m}_{Z}, \omega^{m}_{Z})}$.\vspace{1mm}

\noindent 
{\bf(ii)} $Z$ is Cohen-Macaulay if and only if $\omega^{m}_{Z}$ is maximal Cohen-Macaulay and $\displaystyle{{\mathcal O}_{Z}\simeq{\mathcal H}om(\omega^{m}_{Z}, \omega^{m}_{Z})}$. 
\rm
\begin{proof}
{\bf(i)} The implication is obvious. For the converse, the given isomorphism shows that ${\mathcal O}_{Z}$ has depth at least two and, since by hypothesis ${\rm Sing}(Z)$ has codimension at least two, ${\mathcal H}^{1}_{{\rm Sing}(Z)}({\mathcal O}_{Z})={\mathcal H}^{0}_{{\rm Sing}(Z)}({\mathcal O}_{Z})=$; which characterizes normality.\vspace{1mm}

\noindent 
{\bf(ii)} As with the sheaf ${\mathcal L}^{m}_{Z}$, this maximality property translates into cohomological vanishing and leads, with the same arguments, to the vanishing ${\mathcal H}^{j}({\rm I}\!{\rm R}{\mathcal H}om({\omega}^{m}_{Z}, {\mathcal D}^{\bullet}_{Z}[-m]))=0$ for all $j\not=0$. Hence the equivalence
$$\omega^{m}_{Z}\,{\rm maximal\, Cohen\, Macaulay}\Longleftrightarrow {\mathcal H}om(\omega^{m}_{Z}, \omega^{m}_{Z})\simeq {\rm I}\!{\rm R}{\mathcal H}om({\omega}^{m}_{Z}, {\mathcal D}^{\bullet}_{Z}[-m])$$
If we have ${\mathcal O}_{Z}\simeq{\rm I}\!{\rm R}{\mathcal H}om({\omega}^{m}_{Z}, {\mathcal D}^{\bullet}_{Z}[-m])$, we have necessarily by biduality ${\mathcal D}^{\bullet}_{Z}[-m]\simeq {\omega}^{m}_{Z}$ which means that $Z$ is Cohen Macaulay. \vspace{1mm}

\noindent Conversely, we know that $Z$ is Cohen Macauly if and only if ${\mathcal D}^{\bullet}_{Z}[-m]\simeq \omega^{m}_{Z}$ which implies that ${\rm I}\!{\rm R}{\mathcal H}om({\omega}^{m}_{Z}, {\mathcal D}^{\bullet}_{Z}[-m]))\simeq {\mathcal O}_{Z}$ and then  $\omega^{m}_{Z}$ is maximal Cohen Macaulay sheaf $\,\blacksquare$

\end{proof}\vspace{1mm}

\noindent
\subsubsection{\bf{Quasi-smooth spaces and rational singularities.}}

A simple example of rational singularity is given by quotient singularities (cf. Mumford, Steenbrink-Straten, Burns). We can also mention that a large class of rational singularities is given by quasi-smooth spaces:

\defn{}{}\label{defn5} (\cite{BM1}): A normal and reduced complex analytic space is said to be \emph{quasi-smooth} if, for each point $z$ of $Z$, there exists an open neighborhood $U$ such that there exists a smooth variety $V$ and a finite surjective morphism $f:V\rightarrow U$.

\Prop{}{}\label{P2} Every quasi-smooth complex space $Z$ of dimension $m$ has rational singularities.
\begin{proof} We will use Kempf's criterion recalled in \ref{P1} by showing that $Z$ is necessarily Cohen-Macaulay and that ${\mathcal L}^{m}_{Z}$ is isomorphic to the dualizing sheaf $\omega^{m}_{Z}$.\vspace{1mm}

\noindent To avoid overloading the text, we retain the notations with $Z$ and $V$, keeping in mind that local models must be considered...\vspace{1mm}
 
 \indent
$\bullet$  {\bf{$Z$ is Cohen-Macaulay: }}\vspace{1mm}
 
 \noindent Since the morphism $f:V\rightarrow Z$ is a ramified covering of some degree $k$ (because it is finite over a normal base), it induces a trace morphism  
 ${\mathcal T}{\!r}:f_{*}{\mathcal O}_{V}\rightarrow{\mathcal O}_{Z}$ defined by ${{\mathcal T}r(h)(t):=\sum_{x_{j}\in f^{-1}(t)}h(x_{j})}$ for any holomorphic function $h$ on $V$. As previously stated, the pullback ${\mathcal O}_{Z}\rightarrow f_{*}{\mathcal O}_{V}$ followed by this trace identifies with the morphism $k.{\rm I}d$, and the sheaf ${\mathcal O}_{Z}$ is a direct summand of the sheaf $f_{*}{\mathcal O}_{V}$. We deduce that the vanishing  
 $${\rm H}^{i}_{c}(V, {\mathcal O}_V)={\rm H}^{i}_{c}(Z, f_{*}{\mathcal O}_V)=0,\,\,\forall\,i\not=m$$
 implies the vanishing  
 $${\rm H}^{i}_{c}(Z, {\mathcal O}_Z)=0,\,\,\forall\,i\not=m$$  
 which means that ${\mathcal O}_Z$ is maximal Cohen-Macaulay, i.e., $Z$ is Cohen-Macaulay.\vspace{1mm}

 \indent $\bullet$ {\bf{Isomorphism between ${\mathcal L}^{m}_{Z}$ and $\omega^{m}_{Z}$:}}\vspace{1mm}
 
 \noindent The sheaf ${\mathcal O}_{Z}$ being a direct summand of $f_{*}{\mathcal O}_{V}$, ${\mathcal H}om({\mathcal O}_{Z},\omega^{m}_{Z})$ is also a direct summand of ${\mathcal H}om(f_{*}{\mathcal O}_{V},\omega^{m}_{Z})$, meaning $\omega^{m}_{Z}$ is a direct summand of $f_{*}\Omega^{m}_{V}$. Hence, there is a nontrivial morphism from $\omega^{m}_{Z}$ to $f_{*}\Omega^{m}_{V}$. However, this is possible if and only if ${\mathcal L}^{m}_{Z}=\omega^{m}_{Z}$.\vspace{1mm}
 
 \noindent Indeed, consider the desingularization morphism $\nu: \widetilde{Z}\rightarrow Z$ and the resulting commutative base change diagram:
$$\xymatrix{{\widetilde V}\ar[r]_{\tilde{\nu}}\ar[d]_{\,\,\,\,\,\tilde{f}}&{ V}\ar[d]_{f}\\
\widetilde{Z}\ar[r]_{\nu}&Z}$$
where $\widetilde V$ denotes the fiber product $V\times_{Z}{Z}$.\par\noindent
Then, applying the pullback $\tilde{\nu}^{*}$ to the morphism  
$f^{*}\omega^{m}_{Z}\rightarrow\Omega^{m}_{V}$ and using the commutativity of the diagram, we obtain the morphism $\tilde{f}^{*}{\nu}^{*}\omega^{m}_{Z}\rightarrow\Omega^{m}_{\widetilde V}$, or equivalently, ${\nu}^{*}\omega^{m}_{Z}\rightarrow \tilde{f}_{*}\Omega^{m}_{\widetilde V}$. Since $\tilde{f}$ is also a ramified covering, it induces a trace morphism $\tilde{f}_{*}\Omega^{m}_{\widetilde V}\rightarrow \Omega^{m}_{\widetilde Z}$. Thus, we have the arrow ${\nu}^{*}\omega^{m}_{Z}\rightarrow\Omega^{m}_{\tilde Z}$, which gives rise to a morphism $\omega^{m}_{Z}\rightarrow{\mathcal L}^{m}_{Z}$. This morphism is necessarily injective because these two torsion-free sheaves coincide generically on $Z$. The conclusion follows since the sheaf ${\mathcal L}^{m}_{Z}$ naturally injects into the sheaf $\omega^{m}_{Z}\,\,\blacksquare$
 \end{proof} 
 \vspace{2mm}

 \noindent
\begin{rem} 
 {\bf(i)} To show that $Z$ is Cohen-Macaulay, one can prove that for any local parametrization $g:Z\rightarrow W$ (around a fixed point $z$ in $Z$), the sheaf $g_{*}{\mathcal O}_{Z}$ is locally free. After sufficient localization, we obtain a commutative factorization diagram:
$$\xymatrix{V\ar[rd]_{h}\ar[rr]^{f}&&Z \ar[ld]^{g}\\
&W&}$$
 The morphism $h$ is a ramified covering and thus open. But since $V$ and $W$ are smooth, $h$ is necessarily flat, and consequently, $g_{*}{\mathcal O}_{V}$ is locally free. Since ${\mathcal O}_{Z}$ is a direct summand of $f_{*}{\mathcal O}_{V}$, $g_{*}{\mathcal O}_{Z}$ is then a direct summand of the locally free sheaf $h_{*}{\mathcal O}_{V}$. Hence, $g_{*}{\mathcal O}_{Z}$ is locally free\footnote{This shows that for any local model $f:V\rightarrow Z$, $f$ is always of {\it Tordim}-finite type but not necessarily flat, even if it is open with Cohen-Macaulay fibers and base (it suffices to consider the parametrization of the cone or the Whitney umbrella in ${\Bbb C}^{3}$ given by $(u,v)\rightarrow( u^{2},v^{2},uv)$ and $(u,v)\rightarrow (u^{2},v, uv)$, respectively).}. \vspace{1mm}
 
 \noindent
{\bf(ii)} Another way to see this would be to consider relative duality for a proper morphism to show that the dualizing complex of $Z$ reduces to a single sheaf placed in degree $-m$.\vspace{1mm}

\noindent
Indeed, we have $f^{!}{\mathcal D}^{\bullet}_{Z}={\mathcal D}^{\bullet}_{V}=\Omega^{m}_{V}[m]$, and thus ${\mathcal H}^{j}(f^{!}{\mathcal D}^{\bullet}_{Z}[-m])=0$ for all $j\not= 0$. Since $f$ is finite, we have (by construction) $f_{*}(f^{!}{\mathcal D}^{\bullet}_{Z}[-m])={\rm I}\!\!{\rm R}{\mathcal H}om(f_{*}{\mathcal O}_{V}, {\mathcal D}^{\bullet}_{Z}[-m])$. Since we have a composition ${\mathcal O}_{Z}\rightarrow f_{*}{\mathcal O}_{V}\rightarrow {\mathcal O}_{Z}$ identifying with the identity morphism, the sheaf ${\mathcal O}_{Z}$ is a direct summand of $f_{*}{\mathcal O}_{V}$. It follows that the cohomology of the complex ${\mathcal D}^{\bullet}_{Z}[-m]$ is also a direct summand of the cohomology of the complex ${\rm I}\!\!{\rm R}{\mathcal H}om(f_{*}{\mathcal O}_{V}, {\mathcal D}^{\bullet}_{Z}[-m])$, i.e., of the complex $f_{*}(f^{!}{\mathcal D}^{\bullet}_{Z}[-m])$. However, the vanishing of the cohomology of $f^{!}{\mathcal D}^{\bullet}_{Z}[-m])$ in all nonzero degrees, combined with the finiteness of $f$, leads to the desired conclusion since the dualizing complex of $Z$ has cohomology only in degree $-m$.\vspace{1mm}

\noindent
 \end{rem}
\cor{}{}\label{C3} Let $Z$ be a quasi-smooth complex space. Then:  
\par\noindent
{\bf(i)} For every integer $k\geq 0$, the ${\mathcal O}_{Z}$-coherent sheaves $\omega^{k}_{Z}$ and ${\mathcal L}^{k}_{Z}$ coincide and are maximal Cohen-Macaulay. \par\noindent
{\bf(ii)} Moreover, the canonical morphism  
${\omega^{k}_{Z}\rightarrow {\rm I}\!{\rm R}{\mathcal H}om({\omega}^{m-k}_{Z},  {\omega}^{m}_{Z})}$  
is a quasi-isomorphism.\rm
\begin{proof} \par\noindent
{\bf(i)} To show the equality of sheaves ${\mathcal L}^{k}_{Z}=\omega^{k}_{Z}$ for every integer $k$, it suffices to exhibit a nontrivial morphism from $\omega^{k}_{Z}$ to ${\mathcal L}^{k}_{Z}$. In this context, the existence of such a morphism is equivalent to $\omega^{k}_{Z}$ being a direct summand of the sheaf $f_{*}\Omega^{k}_{V}$.\par\noindent
To see this, first observe that the sheaf ${\mathcal L}^{k}_{Z}$ is a direct summand of $f_{*}\Omega^{k}_{V}$ since it is equipped with a pullback ${\mathcal L}^{k}_{Z}\rightarrow f_{*}\Omega^{k}_{V}$ and a trace morphism $f_{*}\Omega^{k}_{V}\rightarrow{\mathcal L}^{k}_{Z}$ whose composition yields the identity on ${\mathcal L}^{k}_{Z}$. It follows that for every integer $k$, the sheaves ${\mathcal H}om( {\mathcal L}^{m-k}_{Z}, \omega^{m}_{Z})$ are direct summands of ${\mathcal H}om(f_{*}\Omega^{m-k}_{V}, \omega^{m}_{Z})$. But thanks to the isomorphism $\omega^{k}_{Z}\simeq{\mathcal H}om({\mathcal L}^{m-k}_{Z} , \omega^{m}_{Z})$ due to the normality of $Z$, we see that $\omega^{k}_{Z}$ is also a direct summand of $f_{*}\Omega^{k}_{V}$, and we conclude as in \ref{P2} to deduce that ${\mathcal L}^{k}_{Z}\simeq \omega^{k}_{Z}$. \vspace{1mm}
 
 \noindent
{\bf(ii)} Now, let us show that for every integer $k$, the natural morphism $\omega^{k}_{{Z}}\rightarrow {\rm I}\!\!{\rm R}{\mathcal H}om( \omega^{m-k}_{Z}, \omega^{m}_{{Z}})$ is in fact a quasi-isomorphism (which constrains the complex ${\rm I}\!\!{\rm R}{\mathcal H}om( \omega^{m-k}_{Z}, \omega^{m}_{{Z}})$ to be acyclic).  
\par\noindent
Since ${\Omega^{k}_{V}\simeq{\mathcal H}om(\Omega^{m-k}_{V}, \Omega^{m}_{V})\simeq  {\rm I}\!\!{\rm R}{\mathcal H}om(\Omega^{m-j}_{V}, \Omega^{m}_{V})}$, relative duality for the proper morphism $f$ gives ${f_{*}\Omega^{k}_{V}\simeq  {\rm I}\!\!{\rm R}{\mathcal H}om(f_{*}\Omega^{m-k}_{V}, \omega^{m}_{Z})}$, and thus the vanishing of the cohomology sheaves ${\mathcal E}xt^{j}(f_{*}\Omega^{m-k}_{V}, \omega^{m}_{Z})$ for all $j\geq 1$. Moreover, by {\bf(ii)}, the sheaf $\omega^{m-k}_{Z}$ (resp. ${\mathcal L}^{m-k}_{Z}$) is a direct summand of $f_{*}\Omega^{m-k}_{V}$, and consequently, the cohomology sheaf ${\mathcal E}xt^{j}(\omega^{m-k}_{Z}, \omega^{m}_{Z})$ (resp. ${\mathcal E}xt^{j}( {\mathcal L}^{m-k}_{Z}, \omega^{m}_{Z})$) is a direct summand of ${\mathcal E}xt^{j}(f_{*}\Omega^{m-k}_{V}, \omega^{m}_{Z})$, which is known to vanish. \vspace{1mm}
 
 \noindent
Since $Z$ is Cohen-Macaulay and Stein, by duality, this is exactly equivalent to writing  
 $${\rm H}^{j}_{c}(Z,\omega^{k}_{Z})={\rm H}^{j}_{c}(Z,{\mathcal L}^{k}_{Z})=0,\,\,\forall\,k,\,\,\forall\,j\not=m$$  
meaning that these sheaves are maximal Cohen-Macaulay $\,\,\blacksquare$
\end{proof}
\begin{rem} 
Of course, {\bf(i)} can be deduced from (\cite{KeSc}, {\it {Theorem (1.2)}}, {\it {Theorem (1.3)}}), which we recall a little later, but this is not necessary.
\end{rem}
\subsubsection{\bf{Singularities of type ${\mathfrak {F.R}}$, ${\mathfrak {F.R}}^{*}$ .}}
\Prop{}{}\label{P3} Let $Z$ be a complex space. Then, the following equivalences hold:\par\indent
{\bf(i)} $Z$ is of type ${\mathfrak {F.R}}$,\par\indent
{\bf(ii)} ${\mathcal L}^{m}_{Z}$ has depth at least two in $Z$,\par\indent
{\bf(iii)} the normalization of $Z$ is of type ${\mathfrak {F.R}}^{*}$,\par\indent
{\bf(iv)} for every integer $j$, the sheaves ${\mathcal L}^{j}_{Z}$ are $\omega^{m}_{Z}$ (resp. ${\mathcal L}^{m}_{Z}$, resp. ${\mathcal L}^{0}_{Z}$)-reflexive.\rm\vspace{1mm}
\begin{proof} Let $\nu:\overline{Z}\rightarrow Z$ denote the normalization of $Z$. \vspace{1mm}
  
  \noindent {\bf(i)}$\Longleftrightarrow${\bf(ii)}:
  Relative analytic duality for a proper morphism \cite{RRV} and the Grauert-Riemenschneider vanishing theorem \cite{GR70} yield the quasi-isomorphisms
$${\rm I}\!{\rm R}{\mathcal H}om( {\rm I}\!{\rm R}\phi_{*}\Omega^{m}_{\tilde Z} , {\mathcal D}^{\bullet}_{Z})={\rm I}\!{\rm R}\phi_{*}{\rm I}\!{\rm R}{\mathcal H}om( \Omega^{m}_{\tilde Z}, \Omega^{m}_{\tilde Z}[m])={\rm I}\!{\rm R}\phi_{*}{\mathcal O}_{\tilde Z}[m]$$
and consequently, the isomorphism of cohomology sheaves \[{\rm E}xt^{j}({\mathcal L}^{m}_{Z}, {\mathcal D}^{\bullet}_{Z})\simeq {\rm I}\!{\rm R}^{m+j}\phi_{*}{{\mathcal O}_{\tilde Z}}\]
Considering strong duality (cf. \cite{AK}) on these cohomology groups endowed with a separated locally convex topological vector space structure since ${\rm H}^{j}_{c}(U; {\mathcal L}^{m}_{Z})={\rm H}^{j}_{c}(\phi^{-1}(U); {\Omega}^{m}_{\tilde Z})$, we obtain, in particular, for every Stein open subset of $Z$, the isomorphism of separated topological vector spaces
$${\rm E}xt^{j}(U;  {\mathcal L}^{m}_{Z}, {\mathcal D}^{\bullet}_{Z})\simeq ({\rm H}^{-j}_{c}(U, {\mathcal L}^{m}_{Z}))^{'}$$
Now, according to \cite{Ba1}, we have ${\rm Prof}({\mathcal L}^{m}_{Z})\geq 2$ if and only if, for every Stein open subset $U$ of $Z$, ${\rm H}^{1}_{c}(U, {\mathcal L}^{m}_{Z})={\rm H}^{0}_{c}(U, {\mathcal L}^{m}_{Z})=0$. But this is equivalent to the vanishing of the sheaves ${\rm E}xt^{j}(U;  {\mathcal L}^{m}_{Z}, {\mathcal D}^{\bullet}_{Z})=0$ for $j=-1, 0$ and thus to the vanishing of the cohomology groups (Fréchet topological vector spaces) 
$\Gamma(U, {\rm I}\!{\rm R}^{i}\phi_{*}{{\mathcal O}_{\tilde Z}})=0$, for $i=m-1$ and $i=m$. Since $U$ is an arbitrary Stein open subset and $Z$ can always be covered by a countable family of Stein open subsets, this proves the desired equivalence.\vspace{1mm}

 \noindent 
{\bf(ii)}$\Longleftrightarrow${\bf(iii)}: 
This is evident since ${\mathcal L}^{m}_{Z}=\nu_{*}{\mathcal L}^{m}_{\overline Z}$ and depth is preserved under finite direct image. Thus, the equality ${\rm Prof}({\mathcal L}^{m}_{\overline Z})={\rm Prof}({\mathcal L}^{m}_{ Z})$ shows that if $Z$ is of type ${\mathfrak{F.R}}$, the same holds for $\overline Z$. Moreover, the sheaves ${\mathcal L}^{m}_{\overline Z}$ and $\omega^{m}_{\overline Z}$, both of depth at least two, coincide on the regular part and hence everywhere since the singular locus has codimension at least two. Therefore, $\overline Z$ is of type ${\mathfrak{F.R}}^{*}$. \vspace{1mm}

  \noindent Conversely, if $\overline{Z}$ is of type ${\mathfrak{F.R}^{*}}$, we have ${\mathcal L}^{m}_{\overline Z}\simeq\omega^{m}_{\overline Z}$ and thus ${\mathcal L}^{m}_{Z}=\nu_{*}{\mathcal L}^{m}_{\overline Z}\simeq\nu_{*}\omega^{m}_{\overline Z}$. Hence, by the invariance of depth under finite direct image, ${\rm Prof}({\mathcal L}^{m}_{Z})\geq 2$, meaning $Z$ is of type ${\mathfrak{F.R}}$. \vspace{1mm}

\noindent 
{\bf(iii)}$\Longleftrightarrow${\bf(iv)}:\par\indent 
  $\bullet$ The implication {\bf(iv)}$\Longrightarrow${\bf(iii)} is clear since for $j=m$, the $\omega^{m}_{Z}$-reflexivity translates to the isomorphism
$${\mathcal L}^{m}_{Z}\simeq{\mathcal H}om({\mathcal H}om({\mathcal L}^{m}_{Z}, {\omega}^{m}_{Z}), \omega^{m}_{Z})$$  
which implies ${\rm Prof}({\mathcal L}^{m}_{Z})\geq 2$ since this holds for the sheaf $\omega^{m}_{Z}$. As above, we also have ${\rm Prof}({\mathcal L}^{m}_{\overline Z})\geq 2$ and thus $\overline Z$ is of type ${\mathfrak{F.R}}^{*}$.\par\indent
$\bullet$ To establish the converse, note that for any finite modification $f:X\rightarrow Y$ of reduced complex spaces, the canonical morphism
$$f_{*}{\mathcal H}om({\mathcal F}, {\mathcal G})\rightarrow {\mathcal H}om(f_{*}{\mathcal F}, f_{*}{\mathcal G})$$
is an isomorphism if ${\mathcal G}$ is torsion-free. This is easily verified using the natural exact sequence 
$$0\rightarrow {\mathcal N}\rightarrow f^{*}f_{*}{\mathcal F}\rightarrow {\mathcal F}\rightarrow 0$$ 
to which we apply the functor ${\mathcal H}om(-, {\mathcal G})$ to obtain the exact sequence
$$0\rightarrow {\mathcal H}om({\mathcal F}, {\mathcal G})\rightarrow{\mathcal H}om(f^{*}f_{*}{\mathcal F}, {\mathcal G})\rightarrow {\mathcal H}om({\mathcal N}, {\mathcal G})$$
Then, knowing that ${\mathcal N}$ is torsion and ${\mathcal G}$ is not, applying the functor $f_*$ and the adjunction formula, we obtain the desired isomorphism. 
\vspace{1mm}

\noindent On the other hand, it is known that for any finite morphism $f:X\rightarrow Y$ of pure $m$-dimensional complex spaces, we have
$$f_{*}{\mathcal H}om({\mathcal F}, \omega^{m}_{X})\simeq {\mathcal H}om(f_{*}{\mathcal F}, \omega^{m}_{Y})$$
Thus, if $\overline Z$ is of type ${\mathfrak{F.R}}^{*}$, the equality of the sheaves ${\mathcal L}^{m}_{\overline Z}$ and $\omega^{m}_{\overline Z}$ implies the equality of the sheaves ${\mathcal L}^{j}_{\overline Z}$ and $\omega^{j}_{\overline Z}$ for all $j$, and since all these sheaves are $\omega^{m}_{\overline Z}$ (and obviously ${\mathcal L}^{m}_{\overline Z}$)-reflexive, the isomorphisms
$${\mathcal L}^{j}_{\overline Z}\simeq{\mathcal H}om({\mathcal H}om({\mathcal L}^{j}_{\overline Z}, {\mathcal L}^{m}_{\overline Z}), {\mathcal L}^{m}_{\overline Z})\simeq {\mathcal H}om({\mathcal H}om({\mathcal L}^{j}_{\overline Z}, {\omega}^{m}_{\overline Z}), {\omega}^{m}_{\overline Z})$$
yield, upon applying the functor $\nu_{*}$, the isomorphisms
$${\mathcal L}^{j}_{ Z}\simeq{\mathcal H}om({\mathcal H}om({\mathcal L}^{j}_{Z}, {\mathcal L}^{m}_{Z}), {\mathcal L}^{m}_{Z})\simeq {\mathcal H}om({\mathcal H}om({\mathcal L}^{j}_{Z}, {\omega}^{m}_{Z}), {\omega}^{m}_{Z})$$
proving {\bf(iv)}$\,\blacksquare$
\end{proof}

In the special case where the normalization is smooth, we have:
\cor{}{}\label{C4} Let $Z$ be a reduced complex space and $\nu:\overline{Z}\rightarrow Z$ its normalization with $\overline{Z}$ smooth. Then, $Z$ is of type ${\mathfrak{F.R}}$, and the sheaves ${\mathcal L}^{j}_{Z}$ are maximal Cohen-Macaulay for every integer $j\in\{0,1,\cdots,m\}$. Moreover, the sheaves ${\mathcal L}^{j}_{Z}$ are $\omega^{m}_{Z}$ (resp. ${\mathcal L}^{m}_{Z}$, resp. ${\mathcal L}^{0}_{Z}$)-reflexive.\rm \vspace{1mm}
\begin{proof}
By the preceding result, the sheaves ${\mathcal L}^{j}_{Z}$ ($0\leq j\leq m$) all have depth at least two in $Z$, which is consequently of type ${\mathfrak{F.R}}$. To see that they are maximal Cohen-Macaulay, we may assume $Z$ is Stein and hence $\widehat Z$ is as well, and write, using Serre duality, 
$${\rm H}^{i}_{c}(Z, {\mathcal L}^{j}_{Z})\simeq{\rm H}^{i}_{c}(\widehat{Z}, {\Omega}^{j}_{\widehat{Z}})\simeq \big({\rm H}^{m-i}(\widehat{Z}, {\Omega}^{m-j}_{\widehat{Z}})\big)^{'}=0,\,\,\,\forall\,i\not=m,\,\,\forall\,j\in\{0,1,\cdots,m\}.$$ 
Their ${\mathcal L}^{m}_{Z}$ and $\omega^{m}_{Z}$-reflexivity is already established in the previous proposition. Their ${\mathcal L}^{0}_{Z}$-reflexivity follows from the ${\mathcal O}_{\overline Z}$-reflexivity of the sheaves $\Omega^{j}_{\overline Z}$, which is easily deduced from the remark in the preceding result:
$${\mathcal L}^{j}_{Z}=\nu_{*}\Omega^{j}_{\widehat Z}\simeq{\mathcal H}om(\nu_{*}{\Omega}^{m-j}_{\widehat Z}, \nu_{*}{\mathcal O}_{\overline Z}) \simeq {\mathcal H}om({\mathcal L}^{m-j}_{Z}, {\mathcal L}^{0}_{Z})$$
and thus,
  $${\mathcal L}^{j}_{Z}\simeq{\mathcal H}om({\mathcal H}om({\mathcal L}^{j}_{Z}, {\mathcal L}^{0}_{Z}), {\mathcal L}^{0}_{Z})\,\,\blacksquare$$
 \end{proof} 
\vspace{1mm}

 \noindent 
 \begin{rem} {\bf(i)} The sheaf ${\mathcal L}^{0}_{Z}$ always has depth at least two on any arbitrary reduced analytic space. Indeed, if $Z$ is normal, the assertion is trivial since ${\mathcal O}_{Z}\simeq{\mathcal L}^{0}_{Z}$. If not, it suffices to recall that this sheaf is the direct image under the normalization morphism of the structural sheaf (which has depth at least two).\vspace{1mm}

 \noindent 
{\bf(ii)} Note that we always have
 $${\mathcal L}^{m}_{Z}\simeq{\mathcal H}om({\mathcal L}^{0}_{Z}, {\mathcal L}^{m}_{Z}),\,\,\,{\mathcal L}^{0}_{Z}\simeq{\mathcal H}om({\mathcal L}^{m}_{Z}, {\mathcal L}^{m}_{Z})\simeq {\mathcal H}om({\mathcal L}^{m}_{Z}, {\omega}^{m}_{Z})$$
 and if $Z$ is of type ${\mathfrak{F.R}}$,
 $${\mathcal L}^{j}_{Z}\simeq \nu_{*}\omega^{m}_{\overline Z}\simeq {\mathcal H}om({\mathcal L}^{m-j}_{Z}, {\mathcal L}^{m}_{Z})\simeq {\mathcal H}om({\mathcal L}^{m-j}_{Z}, {\omega}^{m}_{Z})$$
\end{rem}
\Prop{}{}\label{P4} Let $Z$ be a reduced complex space. Then, the following are equivalent:\vspace{1mm}

\noindent
{\bf(i)} $Z$ is of type ${\mathfrak {F.R}}$ and its singular locus has codimension at least two,\vspace{1mm}

\noindent
{\bf(ii)} $Z$ is of type ${\mathfrak {F.R}^{*}}$.
\begin{proof} The implication ${\bf(i)}\Longrightarrow {\bf(ii)}$ is evident since by \ref{P2}, the sheaf ${\mathcal L}^{m}_{Z}$ has depth at least two and coincides, on the regular part of $Z$, with the sheaf $\omega^{m}_{Z}$. These two sheaves thus coincide everywhere on $Z$.\vspace{1mm}

\noindent 
Conversely, condition {\bf(ii)} immediately implies that the sheaf ${\mathcal L}^{m}_{Z}$ has depth at least two. To see this, consider the short exact sequence 
$$\xymatrix{0\ar[r]&{\mathcal O}_{Z}\ar[r]&{\mathcal L}^{0}_{Z}\ar[r]&{\mathcal K}\ar[r]&0}$$
where ${\mathcal K}$ is the coherent quotient sheaf whose support is contained in the singular locus $\Sigma:={\rm Sing}(Z))$ of $Z$. We may assume $\Sigma$ is reduced (and smooth if desired) of pure codimension $1$ and ${\rm Supp}({\mathcal K})=\Sigma$. Then, from the distinguished triangle
$$\xymatrix{{\rm I}\!{\rm R}{\mathcal H}om({\mathcal K}, {\mathcal D}^{\bullet}_{Z})\ar[rr]&&{\rm I}\!{\rm R}{\mathcal H}om({\mathcal L}^{0}_{Z}, {\mathcal D}^{\bullet}_{Z})\ar[ld]\\
&{\mathcal D}^{\bullet}_{Z}\ar[lu]^{[1]}&}$$
we obtain the exact sequence 
$$\xymatrix{\cdots\ar[r]&{\mathcal H}^{-m}({\rm I}\!{\rm R}{\mathcal H}om({\mathcal K}, {\mathcal D}^{\bullet}_{Z}))\ar[r]&{\mathcal H}^{-m}({\rm I}\!{\rm R}{\mathcal H}om({\mathcal L}^{0}_{Z}, {\mathcal D}^{\bullet}_{Z}))\ar[r]&{\mathcal H}^{-m}({\mathcal D}^{\bullet}_{Z})\ar[d]\\
&&&{\mathcal H}^{-m+1}({\rm I}\!{\rm R}{\mathcal H}om({\mathcal K}, {\mathcal D}^{\bullet}_{Z}))\ar[d]\\
&&&\vdots}$$
Duality for the morphism $\nu$ and the cohomology vanishing theorem of \cite{RRV} and \cite{GR70} gives the quasi-isomorphism ${\mathcal L}^{m}_{Z}={\rm I}\!{\rm R}\nu_{*}\Omega^{m}_{\hat Z}\simeq {\rm I}\!{\rm R}{\mathcal H}om({\mathcal L}^{0}_{Z}, {\mathcal D}^{\bullet}_{Z})[-m]$ 
and, consequently, the vanishing of the cohomology sheaves ${\mathcal H}^{j}({\rm I}\!{\rm R}{\mathcal H}om({\mathcal L}^{0}_{Z}, {\mathcal D}^{\bullet}_{Z}))$ for all $j\not=-m$. This yields the isomorphism $\displaystyle{{\mathcal H}^{-m+j}({\mathcal D}^{\bullet}_{Z})\simeq{\mathcal H}^{-m+1+j}({\rm I}\!{\rm R}{\mathcal H}om({\mathcal K}, {\mathcal D}^{\bullet}_{Z})),\,\,\,\forall\,j\geq 1}$.\vspace{1mm}

\noindent
Moreover, for dimension reasons of the support\footnote{One can also use duality for the embedding of $\Sigma$ in $Z$, ${{\rm I}\!{\rm R}{\mathcal H}om({\mathcal K}, {\mathcal D}^{\bullet}_{Z})\simeq i_{*}{\rm I}\!{\rm R}{\mathcal H}om({\mathcal K}, {\mathcal D}^{\bullet}_{\Sigma})}$ to see that 
$\displaystyle{{\mathcal H}^{-m}({\rm I}\!{\rm R}{\mathcal H}om({\mathcal K}, {\mathcal D}^{\bullet}_{Z}))= {\mathcal E}xt^{-1}({\mathcal K}, \Omega^{m-}}$.}, we have ${\mathcal H}^{-m}({\rm I}\!{\rm R}{\mathcal H}om({\mathcal K}, {\mathcal D}^{\bullet}_{Z}))=0$ since, for any Stein open subset $U$, 
$\displaystyle{{\rm E}xt^{-m}(U; {\mathcal K}, {\mathcal D}^{\bullet}_{Z})\simeq \big({\rm H}^{m}_{c}(U, {\mathcal K})\big)'}$ and the exact sequence becomes
$$\xymatrix{0\ar[r]&{\mathcal H}^{-m}({\rm I}\!{\rm R}{\mathcal H}om({\mathcal L}^{0}_{Z}, {\mathcal D}^{\bullet}_{Z}))\ar[r]&{\mathcal H}^{-m}({\mathcal D}^{\bullet}_{Z})\ar[r]&{\mathcal H}^{-m+1}({\rm I}\!{\rm R}{\mathcal H}om({\mathcal K}, {\mathcal D}^{\bullet}_{Z}))\ar[r]&0}$$
that is
$$\xymatrix{0\ar[r]&{\mathcal L}^{m}_{Z}\ar[r]&\omega^{m}_{Z}\ar[r]&{\mathcal H}^{-m+1}({\rm I}\!{\rm R}{\mathcal H}om({\mathcal K}, {\mathcal D}^{\bullet}_{Z}))\ar[r]&0}$$ 
\vspace{1mm}

\noindent 
By duality for the embedding of $\Sigma$ in $Z$, we have
$\displaystyle{{\rm I}\!{\rm R}{\mathcal H}om({\mathcal K}, {\mathcal D}^{\bullet}_{Z}[-m+1])\simeq {\rm I}\!{\rm R}{\mathcal H}om({\mathcal K}, {\mathcal D}^{\bullet}_{\Sigma}[-m+1])}.$ Now, as we have seen several times, 
${\mathcal H}^{0}({\rm I}\!{\rm R}{\mathcal H}om({\mathcal K}, {\mathcal D}^{\bullet}_{\Sigma})[-m+1])\simeq{\mathcal H}om({\mathcal K}, \omega^{m-1}_{\Sigma})$ which, when $\Sigma$ is smooth, identifies with ${\mathcal H}om({\mathcal K}, \Omega^{m-1}_{\Sigma})$. Finally, the initial short exact sequence becomes 
$$\xymatrix{0\ar[r]&{\mathcal L}^{m}_{Z}\ar[r]&\omega^{m}_{Z}\ar[r]& {\mathcal H}om({\mathcal K}, \omega^{m-1}_{\Sigma})  \ar[r]&0}$$
  Then, under the hypothesis that ${\mathcal L}^{m}_{Z}$ is isomorphic to $\omega^{m}_{Z}$, we necessarily have ${\mathcal H}om({\mathcal K}, \omega^{m-1}_{\Sigma})=0$; which contradicts our hypothesis\footnote{In general, ${\mathcal H}om({\mathcal K}, \omega^{m-1}_{\Sigma})$ only accounts for the torsion-free part of ${\mathcal K}$ since if ${\mathcal T}$ is its torsion subsheaf, we have ${\mathcal H}om({\mathcal K}, \omega^{m-1}_{\Sigma})\simeq {\mathcal H}om({\mathcal K}/{\mathcal T}, \omega^{m-1}_{\Sigma})$. But 
  for ${\mathcal K}$ torsion-free on $\Sigma$, the canonical injective morphism
 ${\mathcal K}\rightarrow{\mathcal H}om({\mathcal H}om({\mathcal K}, \omega^{m-1}_{\Sigma}), \omega^{m-1}_{\Sigma})$ shows that ${\mathcal K}$ must necessarily be the zero sheaf, meaning $Z$ is normal due to the isomorphism ${\mathcal O}_{Z}\simeq {\mathcal L}^{0}_{Z}$; which contradicts our initial assumption. }. Thus, the singular locus has codimension at least two in $Z\,\blacksquare$  
\end{proof}
\rm
\vspace{1mm}
\begin{rem}
 It is no hard to see that $Z$ is of type ${\mathfrak {F.R}^{*}}$ if and only if the sheaf  ${\mathcal L}^{m}_{Z}$ is $\omega$-reflexive and the singular locus of $Z$ has codimension at least two.\vspace{1mm}

 \indent $\bullet$
By \ref{P4}, we know that the singular locus of $Z$ necessarily has codimension at least two.
Since, by definition, $Z$ is of type ${\mathfrak {F.R}^{*}}$ if and only if the canonical injective morphism ${\mathcal L}^{m}_{Z}\rightarrow\omega^{m}_{Z}$ is bijective, the sheaf ${\mathcal L}^{m}_{Z}$ has depth at least two. Then, the $\omega$-reflexivity of ${\mathcal L}^{m}_{Z}$, i.e., the bijectivity of the canonical injective morphism
$${\mathcal L}^{m}_{Z}\rightarrow{\mathcal H}om({\mathcal H}om({\mathcal L}^{m}_{Z}, {\omega}^{m}_{Z}), \omega^{m}_{Z})$$
is immediate since both sheaves have depth at least two and coincide outside the singular locus which has codimension at least two (moreover, $\omega^{m}_{Z}$ is naturally $\omega$-reflexive).\vspace{1mm}

\indent Conversely, the $\omega$-reflexivity of the sheaf ${\mathcal L}^{m}_{Z}$ implies that its depth is at least two. Since the singular locus has codimension at least two and ${\mathcal L}^{m}_{Z}$ coincides with the sheaf $\omega^{m}_{Z}$ (which has depth at least two), they are globally identified on $Z\,\blacksquare$
 \end{rem}
\noindent One of the consequences of (\cite{KeSc}, {\it {Theorem (1.2)}}, {\it {Theorem (1.3)}}) 
\Th{}{}\label{KeSc} (\cite{KeSc}) \vspace{1mm} 

\noindent
Let $Z$ be a complex space. Then the following are equivalent:\par\noindent
{\bf(i)} $Z$ is of type $\mathfrak{{F.R}^{*}}$\par\noindent
{\bf(ii)} $Z$ is of type $k-\mathfrak{{F.R}^{*}}$ for every integer $k\in\{0, 1,\cdots,m\}$.\rm\vspace{1mm}

\noindent
from which immediately follows the
\cor{}{} \label{C6}Let $Z$ be a complex space. Then the following are equivalent:\par\noindent
{\bf(i)} $Z$ is of type $\mathfrak{{F.R}}$\par\noindent
{\bf(ii)} $Z$ is of type $k-\mathfrak{{F.R}}$ for every integer $k\in\{0, 1,\cdots,m\}$.\rm
\begin{proof}
The implication {\bf(ii)}$\Longrightarrow\,${\bf(i)} is obviously trivial. Conversely, if $Z$ is of type $\mathfrak{{F.R}}$, by preservation of depth under finite direct image, its normalization $\nu:\overline Z\rightarrow Z$ is of type $\mathfrak{{F.R}^{*}}$, and \ref{KeSc} gives the equality of sheaves $\mathcal{L}^{j}_{\overline Z}=\omega^{j}_{\overline Z}$ for all $j\leq k$, and then $\mathcal{L}^{j}_{Z}=\nu_{*}\omega^{j}_{\overline Z}\,\blacksquare$
\end{proof}
\section{\color{blue}{Behavior under proper modification and base change of the classes ${\mathfrak {F.R}}^{*}$ and ${\mathfrak {F.R}}$.}} \noindent 
\lem{}{}\label{L1} Let $\psi: X\rightarrow Y$ be a proper and generically finite morphism of reduced complex spaces of pure dimension $m$. Then, ${\rm I}\!{\rm R}^{j}\psi_{*}{\mathcal L}^{m}_{X}=0,\,\,\,\forall\,j\not=0$ and $\psi_{*}{\mathcal L}^{m}_{X}={\mathcal L}^{m}_{Y}$.\rm
\begin{proof} It suffices to consider the commutative diagram
$$\xymatrix{{\widetilde X}\ar[rd]^{\tilde\phi}\ar[d]_{\phi}&\\
X\ar[r]_{\psi}&Y}$$
in which $\phi$ is a resolution of singularities of $X$ which is, therefore, also a desingularization of $Y$, and the Leray spectral sequence
\[{\rm I}\!{\rm R}^{i}\psi_{*}{\rm I}\!{\rm R}^{j}\phi_{*}\Omega^{m}_{\widetilde X}\Longrightarrow{\rm I}\!{\rm R}^{i+j}\tilde\phi_{*}\Omega^{m}_{\tilde X}\]
degenerates because, by virtue of the Grauert-Riemenschneider theorem \cite{GR70} or Takegoshi \cite{Ta}, we have ${\rm I}\!{\rm R}^{j}\phi_{*}\Omega^{m}_{\widetilde X}=0$ and ${\rm I}\!{\rm R}^{j}\tilde\phi_{*}\Omega^{m}_{\widetilde X}=0$ for all $j$ non-zero.  
Thus, ${\rm I}\!{\rm R}^{i}\psi_{*}{\mathcal L}^{m}_{X}=0$ for all $j\not=0$ and $\psi_{*}{\mathcal L}^{m}_{X}=\psi_{*}(\phi_{*}\Omega^{m}_{\widetilde X})=\tilde\phi_{*}\Omega^{m}_{\widetilde X}={\mathcal L}^{m}_{Y}\,\blacksquare$
\end{proof}

\vspace{1mm}
The following results are three variants of the same result, namely the stability of the $\mathfrak{F.R}$ property under proper modification.

\lem{}{}\label{L2} Let $Z$ be a reduced complex analytic space and $T$ a subspace of $Z$. Let $\Theta: \widetilde{Z} \rightarrow Z$ be a proper modification of $Z$ with center $T$ (with preimage $\widetilde{T}$). Then, if $Z$ is of type $\mathfrak{F.R}$, $\widetilde{Z}$ is also of type $\mathfrak{F.R}$.\rm
\begin{proof} Note that the problem is local on $Z$, and the difficulty lies in studying the depth near points of the subset $\widetilde{T} \cap \text{Sing}(\widetilde{Z}) \cap \Theta^{-1}(\text{Sing}(Z))$. We may assume $\Theta$ is projective, since otherwise it can always be dominated by a blow-up. More precisely, for any Stein open subset $U$ of $Z$, there exists a blow-up with center $T \cap U$, $\phi: V \rightarrow U$, and a blow-up of $\widetilde{U} := \Theta^{-1}(U)$ with center $\widetilde{T} \cap \widetilde{U}$, $\eta: V \rightarrow \widetilde{U}$, such that the following diagram commutes:
$$\xymatrix{V\ar[r]^{\!\!\!\eta}\ar[rd]_{\phi}&\ar[d]^{\Theta}\widetilde{U}\\
&U}$$
If $\Theta$ is a projective morphism, there exists a coherent sheaf $\mathcal{F}$ on $Z$ whose projectivization $\mathbb{P}(\mathcal{F})$ admits the factorization:
$$\xymatrix{\widetilde{Z}\ar[r]^{\!\!\!\sigma}\ar[rd]_{\Theta}&\mathbb{P}(\mathcal{F})\ar[d]^{p}\\
&Z}$$
where $\sigma$ is an embedding and $p$ is the canonical projection. We then reduce to the following situation:
$$\xymatrix{\widetilde{W}\ar[r]^{\!\!\!\sigma}\ar[rd]_{\Theta}&U \times \mathbb{P}_N(\mathbb{C})\ar[d]^{p}\\
&U}$$
Since $p$ is a proper projection, we have:
$$\text{Prof}(\mathcal{L}^m_Z) = \text{Prof}(p_*(\sigma_* \mathcal{L}^m_{\widetilde{Z}})) \leq \text{Prof}(\sigma_* \mathcal{L}^m_{\widetilde{Z}}) = \text{Prof}(\mathcal{L}^m_{\widetilde{Z}})$$
More generally, if the direct image of a coherent sheaf under a proper and flat morphism has depth at least two, the sheaf itself necessarily has depth at least two. \vspace{1mm}

\noindent Indeed, let $\pi: X \rightarrow Y$ be a proper, flat, and surjective morphism of reduced complex spaces, and $\mathcal{F}$ a coherent sheaf on $X$ whose direct image $\pi_* \mathcal{F}$ has depth at least two. Since $\pi$ is, in particular, open, it maps any open neighborhood of a given point $x \in X$ to an open neighborhood of $y := \pi(x)$. The coherent sheaf $\pi_* \mathcal{F}$, being of depth at least two by hypothesis, admits at least two holomorphic function germs $f_y$ and $g_y$ at $y$ forming a regular sequence for the stalk $(\pi_* \mathcal{F})_y$. Taking representatives of $f$ and $g$ on a sufficiently small open neighborhood of $y$, the openness of $\pi$ allows lifting $f$ and $g$ to a regular sequence for $\mathcal{F}$ near any point of the fiber $\pi^{-1}(y)$. Since $x$ lies in this fiber, the conclusion follows (note that the openness of the morphism suffices to reach the conclusion; cf. \cite{Fi}, p.154). $\blacksquare$
\end{proof}

\lem{}{}\label{L3} Let $\phi: X \rightarrow Y$ be a proper modification of reduced complex spaces, with $Y$ of pure dimension $m$. Then, if $Y$ is of type $\mathfrak{F.R}$, $X$ is also of type $\mathfrak{F.R}$.

\begin{proof}\vspace{1mm}

\noindent Consider the commutative diagram:
$$\xymatrix{\widetilde{X}\ar[d]_{\widetilde{\psi}}\ar[rd]_{\rho}\ar[r]^{\widetilde{\phi}}&\widetilde{Y}\ar[d]^{\psi}\\
X\ar[r]_{\phi}&Y}$$
where $\widetilde{Y}$ is a desingularization of $Y$, and $\widetilde{X}$ is a desingularization of the strict transform of $X$ (i.e., the irreducible components of the fiber product $X \times_Y \widetilde{Y}$ that surject onto $\widetilde{Y}$). We also consider the two Leray spectral sequences:
$$\text{E}^{i,j}_2 := \text{IR}^i \psi_* \text{IR}^j \widetilde{\phi}_* \mathcal{O}_{\widetilde{X}}, \quad \text{'E}^{i,j}_2 := \text{IR}^i \phi_* \text{IR}^j \widetilde{\psi}_* \mathcal{O}_{\widetilde{X}}$$
converging to $\text{IR}^{i+j} \rho_* \mathcal{O}_{\widetilde{X}}$.\vspace{1mm}

\noindent Since $\text{E}^{i,j}_2 = 0$ for all $j \neq 0$, we have:
$$\text{IR}^i \psi_* \mathcal{O}_{\widetilde{Y}} \simeq \text{IR}^i \rho_* \mathcal{O}_{\widetilde{X}}$$
because $\widetilde{\phi}_* \mathcal{O}_{\widetilde{X}} \simeq \mathcal{O}_{\widetilde{Y}}$ by normality. \vspace{1mm}

\noindent Grauert's direct image theorem shows that $\Sigma_j := \text{Supp}(\text{IR}^j \widetilde{\psi}_* \mathcal{O}_{\widetilde{X}})$ are analytic subsets of dimension at most $m - j - 1$, and thus:
$$\text{'E}^{i,j}_2 = 0 \quad \forall i + j \geq m.$$
We then deduce:
$$\text{'E}^{0,m-1}_2 \simeq \text{'E}^{0,m-1}_3 \simeq \cdots = \text{'E}^{0,m-1}_{\infty}.$$
By definition, in view of the surjective morphism\footnote{These are general considerations for spectral sequences $\text{E}^{i,j}_r$ in the first quadrant, stabilized with a finite filtration $0 \subseteq F^n \text{H}^n \subseteq \cdots \subseteq F^1 \text{H}^n \subseteq \text{H}^n$ and converging to $\text{H}$, written as $\text{E}^{i,j}_r \Longrightarrow \text{H}^{i+j}$. Then $\text{E}^{i,j}_{\infty} := F^i(\text{H}^{i+j}) / F^{i+1}(\text{H}^{i+j}) = \text{Gr}^i(\text{H}^{i+j})$.}:
$$\text{IR}^{m-1} \rho_* \mathcal{O}_{\widetilde{X}} \twoheadrightarrow \text{'E}^{0,m-1}_{\infty},$$
we note that in general, we only have injections:
$$\text{'E}^{0,m-1}_{\infty} \hookrightarrow \cdots \hookrightarrow \text{'E}^{0,m-1}_3 \hookrightarrow \text{'E}^{0,m-1}_2,$$
and thus:
$$\text{H}^{m-1} \twoheadrightarrow \text{'E}^{0,m-1}_{\infty} \hookrightarrow \text{'E}^{0,m-1}_2$$
for the abutment $\text{H}^{m-1}$. \vspace{1mm}

\noindent In our case\footnote{Note that in general, we only have injections:
$$\text{'E}^{0,m-1}_{\infty} \hookrightarrow \cdots \hookrightarrow \text{'E}^{0,m-1}_3 \hookrightarrow \text{'E}^{0,m-1}_2,$$
and thus:
$$\text{H}^{m-1} \twoheadrightarrow \text{'E}^{0,m-1}_{\infty} \hookrightarrow \text{'E}^{0,m-1}_2$$
for the abutment $\text{H}^{m-1}$.}, we obtain a surjective lateral morphism:
$$\text{IR}^{m-1} \rho_* \mathcal{O}_{\widetilde{X}} \twoheadrightarrow \text{'E}^{0,m-1}_2 = \phi_* \text{IR}^{m-1} \widetilde{\psi}_* \mathcal{O}_{\widetilde{X}}.$$
We now use our hypotheses to see that:
$$\text{IR}^{m-1} \psi_* \mathcal{O}_{\widetilde{Y}} \simeq \text{IR}^{m-1} \rho_* \mathcal{O}_{\widetilde{X}} = 0,$$
and thus $\phi_* \text{IR}^{m-1} \widetilde{\psi}_* \mathcal{O}_{\widetilde{X}} = 0$.\vspace{1mm}

\noindent Since the coherent sheaf $\text{IR}^{m-1} \widetilde{\psi}_* \mathcal{O}_{\widetilde{X}}$ has support of dimension at most $0$, we may assume, after sufficient localization, that it is reduced to a point $x_0$ in the singular locus $\text{Sing}(X)$. Then, by definition, we have the equalities of finite-dimensional vector spaces:
$$(\text{IR}^{m-1} \widetilde{\psi}_* \mathcal{O}_{\widetilde{X}})_{x_0} = \Gamma(X, \text{IR}^{m-1} \widetilde{\psi}_* \mathcal{O}_{\widetilde{X}}) = \Gamma(Y, \phi_* \text{IR}^{m-1} \widetilde{\psi}_* \mathcal{O}_{\widetilde{X}}),$$
which, in view of the vanishing of $\phi_* \text{IR}^{m-1} \widetilde{\psi}_* \mathcal{O}_{\widetilde{X}}$, implies $\text{IR}^{m-1} \widetilde{\psi}_* \mathcal{O}_{\widetilde{X}} = 0$. We conclude that $X$ is of type $\mathfrak{F.R}$. $\blacksquare$
\end{proof}\rm
\vspace{1mm}

\noindent
The following result can be deduced from \lemmaref{L3}, but we provide an independent proof, recalling (cf. [\cite{Fi}, 3.6], [\cite{Si}, 1]) that for any morphism $\pi: X \rightarrow Y$ of complex analytic spaces of dimensions $m$ and $r$ respectively, the degeneracy loci:
$$\text{Deg}^k(\pi) := \{x \in X : \dim X_{\pi(x)} > k\}$$
are analytic. If $X$ and $Y$ are reduced and the morphism is generically open such that every irreducible component $X_j$ of $X$ maps to a single irreducible component $Y_j$ of $Y$, the subsets:
$$\text{F}_j := \{x \in X_j : \dim X_{\pi(x)} > \dim X_j - \dim Y_j\}, \quad \text{F} := \bigcup_j \text{F}_j,$$
and $\pi(\text{F})$ are thin in $Y$ (cf. [\cite{Si}, Lemma 1.1]). If $n$ is the generic dimension of $\pi$, we have a stratification:
$$\cdots \subset \text{Deg}^{k+n+1}(\pi) \subset \text{Deg}^{k+n}(\pi) \subset \cdots \subset \text{Deg}^n(\pi) = \text{F},$$
and a family of $(n+k)$-equidimensional morphisms:
$$\pi^k: \text{Deg}^{k+n-1}(\pi) \setminus \text{Deg}^{k+n}(\pi) \rightarrow Y.$$
Moreover, every point of $\text{Deg}_k(\pi) := \text{Deg}^{k+n-1}(\pi) \setminus \text{Deg}^{k+n}(\pi)$ admits an open neighborhood $U$ in $\text{Deg}_k(\pi)$ such that $\pi^k(U)$ is locally analytic of dimension $\dim(\text{Deg}_k(\pi)) - (n + k) < r$.

\lem{}{}\label{L4} Let $\phi: X \rightarrow Y$ be a proper modification of reduced complex spaces, with $Y$ of pure dimension $m$. If for every subspace $\Sigma$ of codimension at least two in $Y$,
$$\mathcal{H}^i_\Sigma(\mathcal{L}^m_Y) = 0 \quad \text{for } i = 0, 1,$$
then the sheaf $\mathcal{L}^m_X$ satisfies the same property on $X$.\rm

\begin{proof} We wish to show that for every subspace $\widetilde{\Sigma}$ of codimension at least two in $X$,
$$\mathcal{H}^i_{\widetilde{\Sigma}}(\mathcal{L}^m_X) = 0 \quad \text{for } i = 0, 1.$$
For $i = 0$, this is clear since $\mathcal{L}^m_X$ is torsion-free. Moreover, if $\Sigma$ lies outside the exceptional locus $\text{E} := \phi^{-1}(T)$ (where $T$ is the center of the modification, a subspace of codimension at least $1$ in $Y$), the result follows immediately from the hypothesis on $Y$; this reduces the problem to considering only subspaces $\Sigma \subsetneq \text{E}$.\vspace{1mm}

\noindent Applying the functor $\phi_*$ to the injection:
$$\mathcal{H}^1_{\widetilde{\Sigma}}(\mathcal{L}^m_X) \hookrightarrow \mathcal{H}^1_{\phi^{-1}(\phi(\widetilde{\Sigma}))}(\mathcal{L}^m_X),$$
and knowing that:
$$\phi_* \mathcal{H}^1_{\phi^{-1}(\phi(\widetilde{\Sigma}))}(\mathcal{L}^m_X) \simeq \mathcal{H}^1_{\phi(\widetilde{\Sigma})}(\mathcal{L}^m_Y),$$
we obtain the injective morphism:
$$\phi_* \mathcal{H}^1_{\widetilde{\Sigma}}(\mathcal{L}^m_X) \hookrightarrow \mathcal{H}^1_{\phi(\widetilde{\Sigma})}(\mathcal{L}^m_Y).$$
Since $\phi$ is a dominant morphism (i.e., maps irreducible components of $X$ to irreducible components of $Y$), we may assume $X$, $Y$, and $T$ are irreducible. Then:
\par\indent $\bullet$ If $T$ has codimension at least two in $Y$, then $\phi(\Sigma)$ will a fortiori have codimension at least two in $Y$, and thus $\phi_* \mathcal{H}^1_{\widetilde{\Sigma}}(\mathcal{L}^m_X) = 0$.
\par\indent $\bullet$ If $T$ has pure codimension $1$, the restriction of $\phi$ to $\text{E}$ is generically finite, meaning there exists a closed subset $T'$ of $T$ with empty interior such that the restriction $\phi: \text{E} \setminus \phi^{-1}(T') \rightarrow T \setminus T'$ is open, finite, and surjective onto a normal base (i.e., a ramified covering of some degree $k$). In this case, $\phi(\widetilde{\Sigma} \cap \phi^{-1}(T'))$ and $\phi(\widetilde{\Sigma} \cap (\text{E} \setminus \phi^{-1}(T')))$ will have empty interior in $T$, the former because it is contained in $T'$, the latter by purity of the dimension of $\text{E} \setminus \phi^{-1}(T')$ (dimension formula for an open and surjective morphism; cf. \cite{Fi}, pp.143–145). \vspace{1mm}

\noindent In conclusion, we may assume $\phi(\widetilde{\Sigma})$ has codimension at least two in $Y$. Thus, since by hypothesis $\mathcal{H}^1_{\phi(\widetilde{\Sigma})}(\mathcal{L}^m_Y) = 0$, we have $\phi_* \mathcal{H}^1_{\widetilde{\Sigma}}(\mathcal{L}^m_X) = 0$. It remains to verify that this implies the vanishing $\mathcal{H}^1_{\widetilde{\Sigma}}(\mathcal{L}^m_X) = 0$.\vspace{1mm}

\noindent We may assume all spaces $X$, $Y$, $T$, and $\text{E}$ are irreducible and reduced. Then:\vspace{1mm}

\indent $\bullet$ If the restriction $\phi'$ of $\phi$ to $\text{E}$ is equidimensional over $T$, we have:
$$\phi_* \mathcal{H}^1_{\widetilde{\Sigma}}(\mathcal{L}^m_X) = 0 \implies \mathcal{H}^1_{\widetilde{\Sigma}}(\mathcal{L}^m_X) = 0.$$
Indeed, let $\mathcal{F}$ denote the sheaf $\mathcal{H}^1_{\widetilde{\Sigma}}(\mathcal{L}^m_X)$, which is an $\mathcal{O}_X$-module with support in $\text{E}$. Assume $\phi'$ is equidimensional with fibers of pure dimension $n$\footnote{It can locally be factored into a finite morphism followed by a projection around each $x \in \widetilde{\Sigma}$ (cf. \cite{Fi}, \S 3.7, for example):
$$\xymatrix{W \ar@/^1pc/[rr]^{\phi'}\ar[r]_{f}&V \times U\ar[r]_{q}&V},$$
where $W$ is an open neighborhood of $x$, $V$ a neighborhood of $\phi'(x)$, and $U$ a (relatively compact polydisc) open subset in some numerical space $\mathbb{C}^n$, reducing the assertion to the case of a finite, surjective, and open morphism.}. By (\cite{Fi}, \emph{corollary} p.138), $\phi'$ can be assumed open. In this case, the hypothesis $\phi_* \mathcal{F} = 0$ implies that for every open subset $U$ of $T$, $\Gamma(U, \phi_* \mathcal{F}) = 0$, and thus in particular for opens of the form $\phi(V)$. But this implies $\Gamma(V, \mathcal{F}) = 0$; varying $V$ over a base of opens in $\text{E}$, we deduce $\mathcal{F} = 0$.\vspace{1mm}

\noindent If $\phi'$ does not have constant fiber dimension, we stratify it using the degeneracy loci of $\phi'$. More precisely, with the notation of \S 2.1, we consider:
$$\text{Deg}^k(\phi) := \{x \in \text{E} : \dim \text{E}_{\phi(x)} > k\}.$$
Clearly, if $\Sigma \cap \text{E} = \emptyset$, we are back in the equidimensional case. Otherwise, letting $\widetilde{\Sigma}_k := \Sigma \cap \text{Deg}_k(\phi')$, we have:
$$\bigcup_k \mathcal{H}^1_{\widetilde{\Sigma}_k}(\mathcal{L}^m_X) = \mathcal{H}^1_{\widetilde{\Sigma} \cap \text{E}}(\mathcal{L}^m_X).$$
Since for every integer $k$:
$$\phi'^k_* (\mathcal{H}^1_{\widetilde{\Sigma}_k}(\mathcal{L}^m_X)) = 0,$$
and $\phi'^k$ is equidimensional, the result follows. $\blacksquare$
\end{proof} \vspace{1mm}

\noindent
\lem{}{}\label{L5} Let $\phi:X\rightarrow Y$ be a proper modification centered on the singular locus $\Sigma$ of a reduced complex analytic space $Y$ of pure dimension $m$. Let ${\mathcal F}$ be a coherent sheaf on $X$, a subsheaf of the coherent sheaf $\omega^{m}_{X}$, such that the direct image ${\phi_{*}}{\mathcal F}$ is a subsheaf of ${\mathcal L}^{m}_{Y}$. Then, ${\mathcal F}$ is a subsheaf of the sheaf ${\mathcal L}^{m}_{X}$.\rm \vspace{1mm}  \noindent
\begin{proof}
Since $\phi$ is a proper modification of a reduced space of pure dimension, $X$ is also reduced of pure dimension. Since ${\mathcal F}$ is a subsheaf of $\omega^{m}_{X}$, it is torsion-free and ${\rm Supp}({\mathcal F})=X$.  
The assertion will be proven if we show that every section of ${\mathcal F}$ extends to a holomorphic form on any desingularization of $X$. Since $\phi$ is a modification centered on the singular locus of $Y$, it suffices to consider the desingularization of $Y$, as it factors through $\phi$ in the commutative diagram  
$$\xymatrix{&{\tilde Y}\ar[ld]_{\psi'}\ar[rd]^{\psi}&\\  
X\ar[rr]_{\phi}&&Y}$$  
It is then enough to use the definitions and the fact that $\phi_{*}{\mathcal F}\subset \phi_{*}{\mathcal L}^{m}_{X}\,\blacksquare$  
\end{proof}  

\Prop{}{}\label{P7} Let $\pi: Z\rightarrow S$ be a universally $n$-equidimensional morphism with $S$ having singularities of type ${\mathfrak {F.R}}^{*}$. Let $\nu:\widetilde{S}\rightarrow S$ be the normalization morphism and the induced base change diagram 
$$\xymatrix{{\widetilde Z}\ar[r]^{\Theta}\ar[d]_{\tilde\pi}&Z\ar[d]^{\pi}\\  
{\widetilde S}\ar[r]_{\nu}&S}$$  
Then, $${\widetilde Z}\,\,{\rm is\, of\, type}\,\, {\mathfrak {F.R}}^{*}\,({\rm resp.} {\mathfrak {F.R}})\Longleftrightarrow \,Z \,\,{\rm is\, of\, type}\, \,{\mathfrak {F.R}}^{*}\, ({\rm resp.} {\mathfrak {F.R}}).$$  
\begin{proof} Since $\pi$ is open, the total preimage $\widetilde Z:={\widetilde S}\times_{S} Z$ coincides with the strict transform.\vspace{1mm}  

\noindent  
$\bullet$ $\displaystyle{{\widetilde Z}\,\,{\rm is\, of\, type}\,\, {\mathfrak {F.R}}\,\Longleftrightarrow \,Z \,\,{\rm is\, of\, type}\, \,{\mathfrak {F.R}}:}$\vspace{1mm}  

\noindent  
This equivalence is obvious since $\nu$ being a finite modification, the same holds for $\Theta$. Consequently, ${\mathcal L}^{m}_{Z}=\Theta_{*}{\mathcal L}^{m}_{\widetilde Z}$ with depth preservation due to the finiteness of $\Theta$. \vspace{1mm}  

\noindent  
$\bullet$ $\displaystyle{{\widetilde Z}\,\,{\rm is\, of\, type}\,\, {\mathfrak {F.R}}^{*}\,\Longleftrightarrow \,Z \,\,{\rm is\, of\, type}\, \,{\mathfrak {F.R}}^{*}:}$\vspace{1mm}  

\indent $\Rightarrow$: Since $\nu$ is a finite modification centered on the singular locus of $S$, $\Theta$ is a modification centered on $\Sigma:=\pi^{-1}({\rm Sing}(S))$, which is of codimension at least two in $Z$ because $S$ is of type ${\mathfrak {F.R}}^{*}$ and $\pi$ is open with fibers of constant dimension.\vspace{1mm}  

\noindent  
But since ${\widetilde Z}$ is of type ${\mathfrak {F.R}}^{*}$, the sheaf ${\mathcal L}^{m}_{\widetilde Z}$ has depth at least two, and consequently, ${\mathcal L}^{m}_{Z}=\Theta_{*}{\mathcal L}^{m}_{\widetilde Z}$ does as well. Moreover, by hypothesis, $\omega^{m}_{\widetilde Z}$ and ${\mathcal L}^{m}_{\widetilde Z}$ coincide globally on ${\widetilde Z}$ and thus on ${\widetilde Z}\setminus \Theta^{-1}(\Sigma)$. We deduce, therefore, that $\omega^{m}_{ Z}$ and ${\mathcal L}^{m}_{Z}$ coincide on $Z\setminus \Sigma$. But since both have depth at least two, they must necessarily coincide globally on $Z$. Hence the conclusion.\vspace{1mm}  

\indent $\Leftarrow$: It is clear that $\widetilde Z$ is of type ${\mathfrak {F.R}}$ due to the equality ${\mathcal L}^{m}_{Z}=\Theta{*}{\mathcal L}^{m}_{\widetilde Z}$, the finiteness of $\Theta$, and the assumption on $Z$, which requires ${\mathcal L}^{m}_{Z}$ to have depth at least two. On the other hand, since $\Theta$ is a modification, we have an injective direct image (defined in the sense of torsion-free currents) $\Theta_{*}\omega^{m}_{\widetilde Z}\rightarrow\omega^{m}_{ Z} $ and an injective morphism ${\mathcal L}^{m}_{Z}=\Theta{*}{\mathcal L}^{m}_{\tilde Z}\rightarrow {\Theta{*}}\omega^{m}_{\tilde Z}$. Since $Z$ is assumed to be of type ${\mathfrak {F.R}}^{*}$, we deduce the isomorphism
${\mathcal L}^{m}_{Z}\simeq {\Theta_{*}}\omega^{m}_{\tilde Z}\simeq \omega^{m}_{Z}$. But the short exact sequence
$$\xymatrix{0\ar[r]&{\mathcal L}^{m}_{\widetilde Z}\ar[r]&\omega^{m}_{\widetilde Z}\ar[r]&{\mathcal K}\ar[r]&0}$$
where ${\mathcal K}$ is the quotient sheaf, which is torsion and supported on the singular locus of $\widetilde Z$, gives us the short exact sequence
$$\xymatrix{0\ar[r]&\Theta_{*}{\mathcal L}^{m}_{\widetilde Z}\ar[r]&{\Theta_{*}}\omega^{m}_{\widetilde Z}\ar[r]&{\Theta_{*}}{\mathcal K}\ar[r]&0}$$
in which the first arrow is bijective according to what precedes. Hence, $\Theta_{*}{\mathcal K}=0$, which, by virtue of the finiteness of $\Theta$ ensuring the surjectivity of the canonical morphism
$$\Theta^{*}\Theta_{*}{\mathcal K}\rightarrow {\mathcal K}$$
imposes ${\mathcal K}=0$ and, consequently, the desired conclusion $\blacksquare$
\end{proof}\vspace{1mm}
\Prop{}{}\label{P8}
Let $\pi: Z\rightarrow S$ be an $n$-geometrically flat morphism with $S$ having singularities of type ${\mathfrak {F.R}}^{*}$. Let $\phi:{\widetilde S}\rightarrow S$ be a proper modification and  the commutative diagram induced by this base change
$$\xymatrix{{\widetilde Z}\ar[r]^{\tilde\phi}\ar[d]_{\tilde\pi}&Z\ar[d]^{\pi}\\
{\tilde S}\ar[r]_{\phi}&S}$$
 Then, $$(i)\,\,\,{\widetilde Z}\,\,{\rm of \,type}\,\, {\mathfrak {F.R}}^{*}\,\Longrightarrow \,Z \,\,{\rm of \,type}\, \,{\mathfrak {F.R}}^{*}$$ $$(ii)\,\,\,{\widetilde Z}\,\,{\rm of \,type}\,\, {\mathfrak {F.R}}\,\Longleftrightarrow \,Z \,\,{\rm of \,type}\, \,{\mathfrak {F.R}}.$$ 
\rm 
\begin{proof}
{\bf(i) $\displaystyle{{\widetilde Z}\,\,{\rm of~type}\,\, {\mathfrak {F.R}}^{*}\,\Longrightarrow \,Z \,\,{\rm of~type}\, \,{\mathfrak {F.R}}^{*}}$:}
\vspace{1mm}

\noindent 
We deliberately restrict ourselves to proper modifications $\nu$ along the singular locus of $S$ (the general case is strictly identical).\vspace{1mm}

\noindent 
As already mentioned, the openness of $\pi$ ensures the identification of the total and strict preimages. Since $\phi$ is a modification centered at the singular locus ${\rm Sing}(S)$, $\tilde\phi$ is a proper modification centered at $\Sigma:=\pi^{-1}({\rm Sing}(S))$.  
Without any loss of generality, we may assume all spaces  of pure dimension and $S$  normal thanks to the \propositionref{P7}, and that $\phi$ is a desingularization morphism.\vspace{1mm}  

\noindent  
Note that the morphism $\pi$ is then geometrically flat, meaning that its fibers can be endowed with cycle structures by assigning suitable multiplicities (cf. \cite{BM-1}, \cite{K1}, \cite{K2}); this is a notion of local nature on the source and stable under arbitrary base change. This holds, for example, for any equidimensional morphism over a normal base. Regarding our purpose, the subtleties related to this notion are not used, as we only need to preserve the openness and the constancy of the fiber dimension. The morphism $\tilde{\pi}$, obtained by base change, is also geometrically flat (hence open) with fibers of pure dimension $n$. \vspace{1mm}

\noindent
The local nature of the problem allows us to reduce to a diagram of local parameterizations of the form (with some abuse of notation), 
$$\xymatrix{{\widetilde Z}\ar[r]^{\tilde\phi}\ar[d]_{\tilde f}&Z\ar[d]^{f}\\
{\widetilde S}\times U\ar[r]_{{\phi'}:=\phi\otimes {Id}}&S\times U}$$ 
\noindent Our assertion will be proven if we show that every section $\xi$ of the sheaf $\omega^{m}_{Z}$ is a section of the sheaf ${\mathcal L}^{m}_{Z}$. \vspace{1mm}

\noindent
To do this, we use the trace morphism induced by $f$, which, for $\xi$, yields the section ${\rm Tr}_{f}(\xi)$ of the sheaf $\omega^{m}_{S\times U}$. But since $S$ is of type ${\mathfrak {F.R}}^{*}$, this section is, in fact, a section of the sheaf ${\mathcal L}^{m}_{S\times U}$. Since the trace commutes with any base change, we have 
 $${\rm Tr}_{\tilde f}(\tilde\phi^{*}(\xi))={{\phi'}}^{*}{\rm Tr}_{f}(\xi)$$
with ${\phi'}^{*}({\rm Tr}_{f}(\xi))$ naturally defining a section of the sheaf $\Omega^{m}_{{\widetilde S}\times U}$. \vspace{1mm}

\noindent Thus, the meromorphic form ${\tilde\phi}^{*}(\xi)$, a section of the sheaf $\tilde\phi^{*}(\omega^{m}_{Z})$, admits holomorphic traces relative to this type of local parameterizations (pullback of local parameterizations on $Z$). To invoke the universal property of the sheaf $\omega^{m}_{\widetilde Z}$, we would need to ensure that $\tilde\phi^{*}(\xi)$ indeed defines a section of $\omega^{m}_{\widetilde Z}$ and, therefore, verify that its actual poles (excluding virtual poles) lie within ${\rm Sing}(\widetilde Z)$. If this is the case, the analyticity of these traces for arbitrary meromorphic forms of $\omega^{m}_{Z}$ provides us with a pullback morphism
$$\tilde\phi^{*}(\omega^{m}_{Z})\rightarrow \omega^{m}_{\widetilde Z}$$
Since $\widetilde Z$ is of type ${\mathfrak {F.R}}^{*}$, this arrow can also be written as
$\displaystyle{\tilde\phi^{*}(\omega^{m}_{Z})\rightarrow {\mathcal L}^{m}_{\widetilde Z}}$
which, upon applying the direct image functor $\tilde\phi_{*}$ and composing with the canonical and injective morphism
$\omega^{m}_{Z}\rightarrow \tilde\phi_{*}\tilde\phi^{*}(\omega^{m}_{Z})$, yields the morphism
$$\omega^{m}_{Z}\rightarrow {\mathcal L}^{m}_{Z}$$
necessarily injective since generically bijective between torsion-free sheaves. Knowing that ${\mathcal L}^{m}_{Z}$ is a subsheaf of $\omega^{m}_{Z}$, we deduce the desired bijectivity. \vspace{1mm}

\indent
Let us now show that the form $\tilde{\xi}:=\tilde\phi^{*}(\xi)$ is indeed a section of the sheaf $\omega^{m}_{\widetilde Z}$. \vspace{1mm}

\noindent Since it is obviously polar in $\tilde\phi^{-1}({\rm Sing}(Z))$, everything reduces to studying the position of the subspaces $\tilde\phi^{-1}({\rm Sing}(Z)$ and ${\rm Sing}(\widetilde{Z})$, while convincing oneself that the true poles are contained in their intersection. We then distinguish the two fundamental cases: \vspace{1mm}

\indent $\bullet$ $\tilde\phi^{-1}({\rm Sing}(Z))\cap{\rm Sing}(\widetilde Z)=\emptyset$:\vspace{1mm}

\noindent This amounts to saying that $\tilde\phi^{-1}({\rm Sing}(Z))\subset{\rm Reg}(\widetilde Z)$ and, consequently, we find ourselves in the situation of a meromorphic form on a smooth variety, hence with poles in a certain divisor (since in codimension two, holomorphic extension is automatic), with this trace holomorphy condition translating into an ${\bf L}^{2}$ growth condition that guarantees global holomorphic extension\footnote{If $f:V\rightarrow U$ is a finite, surjective morphism and $\xi$ a meromorphic form on $V$ whose direct image in the sense of currents $f_{*}\xi$ is a $\overline\partial$-closed current, then the current defined by $\xi$ is necessarily $\overline\partial$-closed.} thanks to (\cite{Gri}, {\bf Lemma 2.1}, p.349-352). It follows that the poles located in ${\rm Reg}(\tilde Z)$ are purely virtual since the form extends holomorphically there; this brings us back to the essential case.\vspace{1mm}

\indent
$\bullet$ $\tilde\phi^{-1}({\rm Sing}(Z))\cap{\rm Sing}(\widetilde Z)\not=\emptyset$:\vspace{1mm}

\noindent
In this case, it is, a fortiori, polar in ${\rm Sing}(\widetilde Z)$, which is of codimension at least two. Since $\omega^{m}_{\widetilde Z}$ is of depth at least two (and thus satisfies the Hartogs extension principle), $\tilde\phi^{*}(\xi)$ extends to a global section of this sheaf. \vspace{1mm}

\noindent We can therefore consider that the meromorphic form $\tilde\phi^{*}(\xi)$ has poles entirely contained in ${\rm Sing}(\widetilde Z)$. Hence, thanks to the trace characterization, 
$${\rm Tr}_{\tilde f}(\tilde\phi^{*}(\xi))\,{\rm holomorphic }\,\,\Longrightarrow\, \tilde\phi^{*}(\xi)\,{\rm a\,section\,of }\,\omega^{m}_{\widetilde Z}$$
Since $\widetilde{Z}$ is of type ${\mathfrak {F.R}}^{*}$, this means that 
 $\tilde\phi^{*}(\xi)$ is a section of  ${\mathcal L}^{m}_{\widetilde Z}$ and, thus, $\xi$ a section of ${\mathcal L}^{m}_{ Z}$; which shows that $Z$ is indeed of type ${\mathfrak {F.R}}^{*}$. \vspace{1mm}

\noindent Note that if $\phi$ is an arbitrary proper modification in which case $\tilde{S}$ may possibly be singular, we will consider local parameterizations on $\widetilde S$ to obtain local parameterizations on $\widetilde Z$. Since trace formation is compatible with the composition of local parameterization morphisms, we are reduced to the previous case where $\widetilde S$ is smooth.\vspace{2mm}

\noindent
{\bf(ii) $\displaystyle{{\widetilde Z}\,\,{\rm of \,type}\,\, {\mathfrak {F.R}}\,\Longleftrightarrow \,Z \,\,{\rm of \,type}\, \,{\mathfrak {F.R}}}$:}
\vspace{1mm}

\indent $\bullet$  $\Leftarrow$: This implication is given by \lemmaref{L3}.\vspace{1mm}

\indent $\bullet$ $\Rightarrow$:
As noted earlier, we may assume all bases to be normal (cf \propositionref{P7}). Let us revisit the base change diagram and complete it with the normalization $Z_{1}$ (resp. $Z_{2}$) of ${Z}$ (resp. ${\widetilde Z}$):
$$\xymatrix{Z_{2} \ar@/_3pc/[ddd]_{\psi_2}\ar[r]^{\tilde\phi"}\ar[d]_{\alpha}&Z_{1}\ar@{=}[d]\ar@/^3pc/[ddd]^{\psi_1}\\
{\overline{Z}}:={\tilde S}\times_{S} Z_{1}\ar[r]^{\tilde\phi'}\ar[d]_{\phi'_{1}}&Z_{1}\ar[d]^{\phi_1}\\
{\widetilde{Z}}:={\widetilde S}\times_{S} Z\ar[r]^{\tilde\phi}\ar[d]_{\tilde\pi}&Z\ar[d]^{\pi}\\
{\widetilde S}\ar[r]_{\phi}&S}$$ 
in which $\phi'_{1}$ is a finite modification, $\pi\circ \phi_{1}$, $\tilde\pi\circ \phi'_{1}$, $\psi_{1}$, and $\psi_{2}$ are universally $n$-equidimensional by stability under base change of the openness properties and constancy of fiber dimension, and because equidimensionality over a normal base implies the openness of the morphism (in fact, they are even geometrically flat since equidimensional over normal bases). \vspace{1mm}

\noindent
Note that if ${\overline{Z}}:={\widetilde S}\times_{S} Z_{1} =\widetilde{Z}\times_{Z} Z_{1}$ is normal, there is no need to normalize $\widetilde Z$ and the result follows from {\bf(i)} since
$${\widetilde{Z}}\,{\rm of\, type}\, {\mathfrak {F.R}}\Longleftrightarrow {\overline{Z}}\,{\rm of\, type}\, {\mathfrak {F.R}}^{*}\Longrightarrow {{Z_1}}\,{\rm of\, type}\, {\mathfrak {F.R}}^{*}\,\Longleftrightarrow {{Z}}\,{\rm of\, type}\, {\mathfrak {F.R}}$$
In the general case, we have the above commutative diagram, justified by the universal property of normalization, which ensures that the normalization morphism $Z_{2}\rightarrow \widetilde{Z}$ factors through $\overline{Z}$ since $\phi'_{1}$ is a finite modification.\vspace{1mm}

\noindent We want to show that ${\mathcal L}^{m}_{Z_1}=\omega^{m}_{Z_1}$. To do this, we revisit the local factorization diagrams:
$$\xymatrix{Z_{2}\ar[rdd]_{\tilde{f}}\ar[rd]^{\alpha}\ar@/^1pc/[rrd]^{\psi}&&\\
&{\overline Z}\ar[r]^{\tilde\phi'}\ar[d]_{\bar f}&Z_{1}\ar[d]^{f}\\
&\widetilde{S}\times U\ar[r]_{\hat\phi}&S\times U}$$
As in {\bf(i)}, we show that $\phi'^{*}(\xi)$ defines a section of the sheaf $\omega^{m}_{\overline{Z}}$ for any section $\xi$ of $\omega^{m}_{Z_{1}}$. \vspace{1mm}

\noindent 
Thanks to the projection formula for a proper morphism, we obtain the isomorphisms 
 $$\tilde{f}_{*}(\psi^{*}(\omega^{m}_{Z_1}))=\tilde{f}_{*}(\alpha^{*}(\phi'^{*}(\omega^{m}_{Z_1})))\simeq\bar{f}_{*}(\phi'^{*}(\omega^{m}_{Z_1})\otimes\alpha_{*}{\mathcal O}_{Z_2})\simeq \bar{f}_{*}(\phi'^{*}(\omega^{m}_{Z_1})\otimes{\mathcal L}^{0}_{\overline Z}) $$ 
and, consequently, 
 $$\hat\phi_{*}\tilde{f}_{*}(\psi^{*}(\omega^{m}_{Z_1}))\simeq
\hat\phi_{*}\bar{f}_{*}(\phi'^{*}(\omega^{m}_{Z})\otimes {\mathcal L}^{0}_{\overline Z})\simeq f_{*}(\phi'_{*}(\phi'^{*}(\omega^{m}_{Z})\otimes {\mathcal L}^{0}_{\overline Z}))\simeq f_{*}(\omega^{m}_{Z_1})$$
since $\phi'_{*}{\mathcal L}^{0}_{\overline Z}={\mathcal O}_{Z_{1}}$ by the normality of $Z_1$. By composing with the trace morphism induced by $f$, we thus deduce a nontrivial morphism
$$\xymatrix{\hat\phi_{*}\tilde{f}_{*}(\psi^{*}(\omega^{m}_{Z_1}))\ar[rd]\ar[r]&f_{*}(\omega^{m}_{Z_1})\ar[d]^{{\rm Tr_{f}}}\\
&\omega^{m}_{S\times U}}$$
Since $S$ (like $\widetilde S$) is of type ${\mathfrak {F.R}}^{*}$, this gives the morphism 
$$\hat\phi_{*}\tilde{f}_{*}(\psi^{*}(\omega^{m}_{Z_1}))\rightarrow {\mathcal L}^{m}_{S\times U}$$ 
and, since $\hat\phi$ is a modification, this imposes the existence of a morphism 
$$\tilde{f}_{*}(\psi^{*}(\omega^{m}_{Z_1}))\rightarrow {\mathcal L}^{m}_{\widetilde{S}\times U}=\omega^{m}_{\widetilde{S}\times U}$$ 
By construction, this morphism factors through the trace induced by $\tilde{f}$ and, consequently, by the universal property of $\omega^{m}_{Z_2}$, we obtain a morphism 
$\psi^{*}(\omega^{m}_{Z_1})\rightarrow \omega^{m}_{Z_2}$
since it is known that, in general, any finite morphism $h:X\rightarrow Y$ of reduced complex spaces of pure dimension $m$ induces, via relative duality, an isomorphism of coherent sheaves (\cite{K5})
$$f_{*}{\mathcal H}om_{{\mathcal O}_{X}}({\mathcal F}, \omega^{m}_{X})\simeq {\mathcal H}om_{{\mathcal O}_{Y}}(f_{*}{\mathcal F}, \omega^{m}_{Y})$$
However, since $Z_2$ is of type ${\mathfrak {F.R}}^{*}$, this morphism can also be written as
$$\psi^{*}(\omega^{m}_{Z_1})\rightarrow {\mathcal L}^{m}_{Z_2}$$
from which, as before, we deduce an injective morphism $\omega^{m}_{Z_1}\rightarrow {\mathcal L}^{m}_{Z_1}$ and, consequently, the identification $${\mathcal L}^{m}_{Z_1}\simeq\omega^{m}_{Z_1}$$
which of course means that $Z_1$ is of type ${\mathfrak {F.R}}^{*}$. But since a reduced complex space is of type ${\mathfrak {F.R}}$ if and only if its normalization is of type ${\mathfrak {F.R}}^{*}$ according to \corollaryref{C4}, it follows that $Z$ is of type ${\mathfrak {F.R}}\,\blacksquare$
\end{proof}
\vspace{2mm}

\noindent
\rm It is now easy to generalize \cite{AS84} to the case of a singular base.
\cor{}{}\label{C11} Let $\pi:Z\rightarrow S$ be a flat morphism of reduced complex spaces where the base and the fibers are of type ${\mathfrak {F.R}}^{*}$. Then, $Z$ is of type ${\mathfrak {F.R}}^{*}$.\rm 
\begin{proof}
Consider the base change diagram induced by a resolution of singularities $\nu:{\widetilde S}\rightarrow S$, where the base change diagram
$$\xymatrix{{\widetilde Z}\ar[r]^{\theta}\ar[d]_{\tilde\pi}&Z\ar[d]^{\pi}\\
{\widetilde S}\ar[r]_{\nu}&S}$$ 
is such that the strict preimage ${\widetilde Z}$ coincides with the total preimage because $\pi$ is open. Flatness being preserved under base change, and the fibers over ${\tilde S}$ having the same singularity type as those over $S$ (being analytically isomorphic), we reduce to the case treated in \cite{AS84} and conclude by applying \propositionref{P8}$\,\blacksquare$
\end{proof}

\phantomsection\addcontentsline{toc}{part}{Proof of the main results.}

\section{ \color{blue}{Transfer of properties ${\mathfrak{R}}$, ${\mathfrak{F.R}}$, and ${\mathfrak{F.R}}^{*}$ from the source to the base: \theoremref{Th1}.}}

Let $\pi:Z\rightarrow S$ be an $n$-geometrically flat morphism of reduced complex spaces of pure dimensions $m$ and $r$, respectively. We will show that if $Z$ has rational (resp. ${\mathfrak{F.R}}^{*}$, resp. ${\mathfrak{F.R}}$) singularities, then $S$ has rational singularities (${\mathfrak{F.R}}^{*}$, resp. ${\mathfrak{F.R}}$). \vspace{2mm}

\noindent {\bf(i) $Z$ with rational singularities $\Longrightarrow S$ with rational singularities.}\vspace{1mm}

\noindent To do this, we first observe that:\vspace{1mm}

\indent
{\bf a) $Z$ normal $\,\Longrightarrow S$ normal.}\vspace{1mm}

\noindent
Since the problem is local on $S$ and $Z$, and $\pi$ is $n$-geometrically flat, we may assume, given a local factorization around a point $z$ of $Z$:
$$\xymatrix{Z\ar[r]^{f}\ar[rd]_{\pi}&S\times U\ar[d]^{q} \\
&S}$$
where $f$ is finite, open, and surjective, $q$ is the natural projection, and $U$ is a relatively compact Stein open subset of a numerical space ${\Bbb C}^{n}$. Thus, it suffices to verify the assertion for $f$.\vspace{1mm}

\noindent
Consider a meromorphic function $g$ that is locally bounded on $S\times U$. Clearly, its pullback also defines a locally bounded meromorphic function, which must then be holomorphic due to the normality of $Z$. However, since $f$ is $0$-geometrically flat, it induces a trace morphism $f_{*}{\mathcal O}_{Z}\rightarrow {\mathcal O}_{S\times U}$. Thus, $f_{*}(f^{*}(g))$ extends holomorphically over $S$. But since $f$ is a ramified covering of some degree $k$, we have $f_{*}(f^{*}(g))=k.g$. Therefore, ${\frac{1}{k}}f_{*}(f^{*}(g))$ is a holomorphic extension of $g$; hence, $S$ is normal.\vspace{1mm}

\indent  
{\bf(b) $Z$ Cohen-Macaulay $\,\Longrightarrow S$ Cohen-Macaulay.}\vspace{1mm}

\noindent
Let ${\mathcal D}^{\bullet}_{Z}$ (resp. ${\mathcal D}^{\bullet}_{S\times U}$) be the dualizing complex of $Z$ (resp. $S\times U$). By hypothesis, ${\mathcal H}^{j}( {\mathcal D}^{\bullet}_{Z})=0$ for all integers $j \neq -m$, and by relative analytic duality for a finite morphism, we have 
$${\rm I}\!{\rm R}f_{*} {\mathcal D}^{\bullet}_{Z}={\rm I}\!{\rm R}{\mathcal H}om({\rm I}\!{\rm R}f_{*}{\mathcal O}_{Z}, {\mathcal D}^{\bullet}_{S\times U})$$
that is, $f_{*} {\mathcal D}^{\bullet}_{Z}={\rm I}\!{\rm R}{\mathcal H}om(f_{*}{\mathcal O}_{Z}, {\mathcal D}^{\bullet}_{S\times U})$. It follows that the homology of the complex ${\mathcal D}^{\bullet}_{S\times U}$ is a direct summand of the homology of $f_{*} {\mathcal D}^{\bullet}_{Z}$. Indeed, since ${\mathcal O}_{S\times U}$ is a direct summand of $f_{*}{\mathcal O}_{Z}$, ${\mathcal D}^{\bullet}_{S\times U}$ is also a direct summand of ${\rm I}\!{\rm R}{\mathcal H}om(f_{*}{\mathcal O}_{Z}, {\mathcal D}^{\bullet}_{S\times U})$.  
By the degeneration of the Leray spectral sequence for $f$, we have the isomorphism 
$f_{*}{\mathcal H}^{j}( {\mathcal D}^{\bullet}_{Z})\simeq {\mathcal H}^{j}( f_{*}{\mathcal D}^{\bullet}_{Z})$, and thus the vanishing 
$${\mathcal H}^{j}( {\mathcal D}^{\bullet}_{S\times U})=0,\,\forall\,j \neq -m.$$
Hence, $S\times U$ is Cohen-Macaulay, and consequently, $S$ is Cohen-Macaulay (since $U$ is smooth).\vspace{1mm} 

\noindent
An alternative approach would be to show that the structure sheaf ${\mathcal O}_{S}$ is maximally Cohen-Macaulay, characterized by the vanishing of the cohomology 
$${\rm H}^{i}_{c}(S, {\mathcal O}_{S})=0,\,\forall\,i \neq r.$$
To verify this, it suffices to use the hypothesis 
$${\rm H}^{i}_{c}(Z, {\mathcal O}_{Z})=0,\,\forall\,i \neq m$$
and the fact that ${\mathcal O}_{S\times U}$ is a direct summand of $f_{*}{\mathcal O}_{Z}$. Indeed, in this case, ${\rm H}^{i}_{c}(Z, {\mathcal O}_{Z})\simeq {\rm H}^{i}_{c}(S\times U, f_{*}{\mathcal O}_{Z})$, and we have 
$${\rm H}^{i}_{c}(S\times U, {\mathcal O}_{S\times U})=0,\,\forall\,i \neq m:=n+r.$$
Then, a Künneth formula yields the decomposition 
$${\rm H}^{i}_{c}(S\times U, {\mathcal O}_{S\times U})=\bigoplus_{i'+j'=i}{\rm H}^{i'}_{c}(S, {\mathcal O}_{S})\widehat{\otimes_{{\Bbb C}}}{\rm H}^{j'}_{c}(U, {\mathcal O}_{U}).$$
Since ${\rm H}^{j'}_{c}(U, {\mathcal O}_{U})=0$ for all $j' \neq n$, we deduce that 
$${\rm H}^{i'}_{c}(S, {\mathcal O}_{S})=0,\,\forall\,i' \neq r,$$
which implies that $S$ is Cohen-Macaulay.\vspace{2mm}

\noindent
Since $S$ is normal and Cohen-Macaulay, the rationality of its singularities relies on the equality ${\mathcal L}^{r}_{S}=\omega^{r}_{{S}}$ by Kempf's criterion (cf \propositionref{P1}).\vspace{1mm}

\noindent 
By hypothesis, we have the equality ${\mathcal L}^{m}_{Z}=\omega^{m}_{{Z}}$, which, by \cite{KeSc}, {\it theorem (1.2)}, implies the equalities ${\mathcal L}^{k}_{Z}=\omega^{k}_{Z}$ for all $k \leq m$; this is a profound result whose proof requires powerful tools.\vspace{1mm}

\noindent 
Again, by restricting to a suitable neighborhood adapted to a cycle $|Z_{s_{0}}|$, we may assume the morphism is finite. Let $\xi$ be a section of the sheaf $\omega^{r}_{S}$, and let $\xi':=\xi\vert_{{\rm Reg}(S)}$ be its restriction to the regular part. Then, $\pi^{*}(\xi')$ defines an $r$-holomorphic form on  the dense subset  $\pi^{-1}({\rm Reg}(S))=Z\setminus\pi^{-1}({\rm Sing}(S))$, whose complement $\pi^{-1}({\rm Sing}(S))$ is a subspace of codimension at least two in $Z$. Since $\omega^{r}_{Z}$ has depth at least two (satisfying Hartog's extension principle), $\pi^{*}(\xi')$ extends naturally to a global section $\pi^{*}(\xi)$ of $\omega^{r}_{Z}$. However, by \cite{KeSc}, {\it theorem (1.2)}, we have ${\mathcal L}^{r}_{Z}\simeq\omega^{r}_{{Z}}$. Thus, we construct a pullback 
$$\pi^{*}\omega^{r}_{Z}\rightarrow {\mathcal L}^{r}_{Z}.$$
This necessarily implies the desired isomorphism. Indeed, we have the diagram 
$$\xymatrix{\tilde{Z}\ar[r]^{\theta}\ar[d]_{\tilde{f}}&Z\ar[d]_{f}\\
\tilde{S}\ar[r]_{\nu}&S}$$
where $\nu$ is a desingularization and $\tilde{Z}$ is the desingularization of the fiber product $\tilde{S}\times_{S} Z$.\vspace{1mm}

\noindent
We obtain the morphisms 
$$\theta^{*}f^{*}\omega^{r}_{Z}\rightarrow \theta^{*}{\mathcal L}^{r}_{Z}\rightarrow \Omega^{r}_{\tilde{Z}}$$
or 
$$\nu^{*}\omega^{r}_{Z}\rightarrow \tilde{f}_{*}\Omega^{r}_{\tilde{Z}}.$$
But the direct image in the sense of currents yields the morphism $\tilde{f}_{*}\Omega^{r}_{\tilde{Z}}\rightarrow \Omega^{r}_{\tilde{S}}$, and consequently, the morphism 
$$\nu^{*}\omega^{r}_{Z}\rightarrow\Omega^{r}_{\tilde{S}},$$
from which we derive 
$$\omega^{r}_{S}\rightarrow{\mathcal L}^{r}_{S},$$
necessarily injective since generically bijective between torsion-free sheaves. Since ${\mathcal L}^{r}_{S}$ is a subsheaf of $\omega^{r}_{S}$, the conclusion follows, and hence $S$ has rational singularities.\vspace{1mm}  

\noindent
{\bf(ii) $Z$ of type ${\mathfrak{F.R}}^{*}$ $\Longrightarrow S$ of type ${\mathfrak{F.R}}^{*}$.}\par\noindent

Since the problem is local on $Z$ and $S$, we use the local factorization from {\bf(i) a)} and the reasoning from {\bf(i) b)}, reducing the problem essentially to the case of a finite and geometrically flat morphism $\pi:Z\rightarrow S\times U$. Consider a section $\xi$ of the sheaf $\omega^{m}_{S\times U}$. Since it is unclear whether the singular locus of $S$ has codimension at least two, we cannot use the previous extension argument. However, we know that $f^{*}(\xi)$ is holomorphic on $Z\setminus f^{-1}({\rm Sing}(S\times U))$ because $\xi$ is holomorphic on $S\times U\setminus {\rm Sing}(S\times U)$. It is thus necessary to examine the intersection $\Sigma:=f^{-1}({\rm Sing}(S\times U))$ with ${\rm Sing}(Z)$. We observe that if $\Sigma\cap{\rm Sing}(Z)=\emptyset$ or, equivalently, $\Sigma\subset {\rm Reg}(Z)$, we are in the situation of \propositionref{P2}, which guarantees that the dense open subset $f({\rm Reg}(Z))$ has rational singularities and is thus of type ${\mathfrak{F.R}}^{*}$. This leads us to focus primarily on the case where $\Sigma\cap{\rm Sing}(Z)\neq\emptyset$ (particularly when $\Sigma\subseteq{\rm Sing}(Z))$. However, since $Z$ is of type ${\mathfrak{F.R}}^{*}$, ${\rm Sing}(Z)$ has codimension at least two, and hence so does $\Sigma\cap{\rm Sing}(Z)$. Since any meromorphic form holomorphic in codimension two extends naturally to a section of $\omega^{m}_{Z}$, the same holds for $f^{*}(\xi)$. At this stage, we can revisit the arguments from {\bf{(i)b)}} to deduce that the pullback $f^{*}\omega^{m}_{{S\times U}}\rightarrow \omega^{m}_{Z}$ can also be written as $f^{*}\omega^{m}_{{S\times U}}\rightarrow{\mathcal L}^{m}_{{Z}}$, since $Z$ is of type ${\mathfrak{F.R}}^{*}$. As before, this forces the isomorphism ${\mathcal L}^{m}_{{S\times U}}\simeq \omega^{m}_{{S\times U}}$. Denoting the canonical projections by $p_{1}:S\times U\rightarrow S$ and $p_{2}:S\times U\rightarrow U$, we have the natural decompositions 
$${\mathcal L}^{m}_{S\times U}={p_{1}}^{*}({\mathcal L}^{m}_{S})\otimes_{{\mathcal O}_{S\times U}}{p_{1}}^{*}(\Omega^{n}_{U})\simeq {\mathcal L}^{r}_{S}\widehat{\otimes}_{{\Bbb C}}\Omega^{n}_{U}$$
$${\omega}^{m}_{S\times U}={p_{1}}^{*}({\omega}^{r}_{S})\otimes_{{\mathcal O}_{S\times U}}{p_{1}}^{*}(\Omega^{n}_{U})\simeq{\omega}^{r}_{S}\widehat{\otimes}_{{\Bbb C}}\Omega^{n}_{U},$$
from which we immediately derive the isomorphism ${\mathcal L}^{r}_{{S}}\simeq \omega^{r}_{{S}}$.\vspace{1mm} 

\noindent
{\bf(iii) $Z$ of type ${\mathfrak{F.R}}$ $\Longrightarrow S$ of type ${\mathfrak{F.R}}$.}\par\noindent
We can argue globally by considering the commutative diagram 
$$\xymatrix{\overline{Z}\ar[d]_{\eta}\ar[rd]&&\\
\widehat{Z}\ar[r]^{\theta}\ar[d]_{\hat{\pi}}&Z\ar[d]_{\pi}\\
\hat{S}\ar[r]_{\nu}&S}$$
where the square is a base change diagram induced by the normalization $\nu$, and $\overline Z$ is the normalization of $Z$.\vspace{1mm}

\noindent Since $\hat S$ is normal and the composition $\tilde{\pi}\circ \eta$ is equidimensional, it is necessarily geometrically flat. Since $Z$ is of type ${\mathfrak{F.R}}$, $\overline Z$ is of type ${\mathfrak{F.R}}^{*}$. By {\bf(ii)}, this implies that $\hat S$ is also of type ${\mathfrak{F.R}}^{*}$, which means $S$ is of type ${\mathfrak{F.R}}$.\vspace{1mm}

\noindent 
We can also argue locally by reusing the factorization from {\bf(i) a)}, reducing the problem essentially to the case of a finite and geometrically flat morphism $\pi:Z\rightarrow S\times U$. The invariance of depth under finite morphisms and the existence of a trace morphism $\pi_{*}{\mathcal L}^{m}_{Z}\rightarrow {\mathcal L}^{m}_{S\times U}$ show, via the pullback ${\mathcal L}^{m}_{S\times U}\rightarrow \pi_{*}{\mathcal L}^{m}_{Z}$, that the sheaf ${\mathcal L}^{m}_{S\times U}$ is a direct summand of the sheaf $\pi_{*}{\mathcal L}^{m}_{Z}$, which has depth at least two. With the notation of {\bf(ii)}, we have the decomposition 
$${\mathcal L}^{m}_{S\times U}={p_{1}}^{*}({\mathcal L}^{m}_{S})\otimes_{{\mathcal O}_{S\times U}}{p_{1}}^{*}(\Omega^{n}_{U})\simeq {\mathcal L}^{m}_{S}\widehat{\otimes}_{{\Bbb C}}\Omega^{n}_{U},$$
showing that ${\rm Prof}({\mathcal L}^{m}_{S\times U})\geq 2 \Longleftrightarrow {\rm Prof}({\mathcal L}^{r}_{S})\geq 2$
since $\Omega^{n}_{U}$ is locally free. Hence, the conclusion follows. $\blacksquare$ \vspace{1mm}
\vfill\eject 
\section{\color{blue}{Transfer of properties ${\mathfrak{F.R}}^{*}$ from Fibers to the Total Space:  \theoremref{Th2}.}}
We propose in what follows two different proofs of this result. The first is essentially a residue-type method used in \cite{AS84}, but with the notable difference that the lack of flatness requires some additional effort using Grauert's notion of active elements. The second, on the other hand, consists of reducing to a flat morphism by flattening the initial morphism to then apply \cite{AS84} or its generalized form given in \corollaryref{C11}. Depending on whether the morphism is proper or not, we appeal to global or local flattening theorems (\cite{H1}, \cite{HLT}, or \cite{Si}). \vspace{1mm}

\noindent 
Let $\displaystyle{\pi:Z\rightarrow S}$ be a universally $n$-equidimensional morphism (i.e., open with fibers of pure dimension $n$) whose fiber singularities, as well as those of $S$, are of type $\mathfrak{F.R}^{*}$. Our goal is to show that $Z$ has the same type of singularities. \vspace{1mm}

\noindent 
The problem is local in nature on $Z$ (and $S$), as it amounts to showing that for every point $z$ of $Z$ lying on the fiber $Z_{s}:=\pi^{-1}(\pi(z))$, the following implication holds:
$${\mathcal L}^{n}_{Z_{s}, z}\simeq \omega^{n}_{Z_{s}, z}\,\Longrightarrow\, {\mathcal L}^{m}_{Z, z}\simeq \omega^{m}_{Z, z}.$$
Since the properties of being open, having fibers of constant dimension, and being of type ${\mathfrak{F.R}}$ or ${\mathfrak{F.R}}^{*}$ are stable under base change (the new fibers being analytically isomorphic to the original ones), and since, moreover, the ${\mathfrak{F.R}}$ condition is preserved under blowups or proper modifications in general (cf. \lemmaref{L2} and \lemmaref{L3}), we may always assume $S$ is normal or even smooth thanks to \propositionref{P7} or \propositionref{P8}. Let us also recall, as already mentioned several times, that the openness of the morphism $\pi$ (and its transforms under successive base changes) ensures that strict and total transforms coincide. \vspace{1mm}

\noindent
\subsection{{Method 1: Active Elements and Fundamental Exact Sequences}}\vspace{1mm}

\noindent Recall that an element $g$ of an analytic algebra is said to be {\emph{active}} if its residual class in its reduction is a {\emph{non-zero divisor}}. We will primarily use (\cite{Fi}, {\bf{Lemma 2}}, p. 144), which states that the pullback of an active element by an open morphism of analytic spaces (not necessarily reduced) with a pure-dimensional base remains an active element. \vspace{1mm}

\indent
{\bf{(i) Smooth Base of Dimension 1.}} Let $s\in S$, $U$ an open neighborhood of $s$, and $f\in\Gamma(U, {\mathcal O}_{S})$ such that its germ $f_{s}$ generates ${\mathcal M}_{S,s}$ and is thus a non-zero divisor, hence active. Since $\pi$ is an open morphism, we know from (\cite{Fi}, {\bf{Lemma 1, Lemma 2}}, p. 144) that for every $x\in X$ and $s:=\pi(x)$, any active element $f$ of ${\mathcal M}_{S,s}$ lifts to an active element of ${\mathcal M}_{Z,z}$. We will loosely denote this element by $\tilde{f}:=\pi^{*}(f)$. Clearly, for every $x\in \pi^{-1}(U)$, $Z_{s}\vert_{\pi^{-1}(U)}=\{\tilde{f}=0\}$. Then, with some deliberate notational abuse omitting localization, multiplication by $\tilde{f}$ yields the short exact sequence 
$$\xymatrix{0\ar[r]&{\tilde{f}}{\mathcal O}_{Z}\ar[r]^{\bullet {\tilde{f}}}&{\mathcal O}_{Z}\ar[r]&{{ i_{s}}}_{*}{\mathcal O}_{Z_{s}}\ar[r]&0}$$
to which we apply the functors ${\mathcal H}om(-,\omega^{m}_{Z})$ and ${\mathcal H}om(-, {\mathcal L}^{m}_{Z})$, respectively, to obtain the commutative diagram
$$\xymatrix{0\ar[r]&\omega^{m}_{Z}\ar[r]&\omega^{m}_{Z, {\tilde{f}}}\ar[r]&{\mathcal E}xt^{1}({{ i_s}}_{*}{\mathcal O}_{{Z_s}}, {\omega}^{m}_{Z})\ar[r]&0\\
0\ar[r]&{\mathcal L}^{m}_{ Z}\ar@{^{(}->}[u]\ar[r]&{\mathcal L}^{m}_{ Z, {\tilde{f}}}\ar[r]\ar@{^{(}->}[u]&{\mathcal E}xt^{1}({{ i_s}}_{*}{\mathcal O}_{{Z_s}}, {\mathcal L}^{m}_{Z})\ar[u]\ar[r]&0}$$
where $\omega^{m}_{Z, {\tilde{f}}}:={\mathcal H}om({\tilde{f}}.{\mathcal O}_{Z},\omega^{m}_{Z})$ and ${\mathcal L}^{m}_{Z, {\tilde{f}}}:={\mathcal H}om({\tilde{f}}.{\mathcal O}_{Z}, {\mathcal L}^{m}_{Z})$. \vspace{1mm}

\noindent 
We will complete this diagram by showing the existence of two natural and injective morphisms 
$${\mathcal E}xt^{1}({{ i_s}}_{*}{\mathcal O}_{{Z_s}}, {\omega}^{m}_{Z})\rightarrow{{ i_s}}_{*}({\omega}^{m-1}_{Z_s})$$
and 
$${{ i_s}}_{*}({\mathcal L}^{m-1}_{Z_s})\rightarrow{\mathcal E}xt^{1}({{ i_s}}_{*}{\mathcal O}_{{Z_s}}, {\mathcal L}^{m}_{Z})$$
which will allow us to establish our result. \vspace{1mm}

\indent
If $Z$ is Cohen-Macaulay then  this injective morphism  ${\mathcal E}xt^{1}(i_{*}{\mathcal O}_{Z_0}, \omega^{m}_{Z})\rightarrow {i_{0}}_{*}\omega^{m-1}_{{ Z}_{0}}$ is bijective, and we have a commutative diagram  
$$\xymatrix{{i_{0}}_{*}\omega^{m-1}_{{ Z}_{0}}\eq[r]&{\mathcal E}xt^{1}(i_{*}{\mathcal O}_{Z_0}, \omega^{m}_{Z})\\  
{i_{0}}_{*}{\mathcal L}^{m-1}_{Z_0}\ar@{^{(}->}[u]\ar@{^{(}->}[r]&{\mathcal E}xt^{1}({i_{0}}_{*}{\mathcal O}_{Z_0}, {\mathcal L}^{m}_{Z})\ar@{^{(}->}[u]}$$  
legitimately prompting the question of whether the residue morphism of the first line induces a residue morphism on the second line.\vspace{1mm}

\indent 
$\bullet$ For the first, observe that the canonical morphism $\displaystyle{\omega^{m}_{Z}[m]\rightarrow{\mathcal D}^{\bullet}_{Z}}$ induces a natural arrow ${\rm I}\!{\rm R}{\mathcal H}om({{ i_s}}_{*}{\mathcal O}_{{Z_s}}, {\omega}^{m}_{Z}[m])\rightarrow {\rm I}\!{\rm R}{\mathcal H}om({{ i_s}}_{*}{\mathcal O}_{{Z_s}}, {\mathcal D}^{\bullet}_{Z})$, which, at the level of cohomology of degree $-(m-1)$, yields the desired morphism
$${\mathcal  H}^{-m+1}({\rm I}\!{\rm R}{\mathcal H}om({{ i_s}}_{*}{\mathcal O}_{{Z_s}}, {\omega}^{m}_{Z}[m]))\rightarrow {{ i_s}}_{*}{\mathcal  H}^{-m+1}({\mathcal D}^{\bullet}_{Z_s})$$
since, by construction, 
${{ i_s}}_{*} {\mathcal D}^{\bullet}_{Z_s}={\rm I}\!{\rm R}{\mathcal H}om({{ i_s}}_{*}{\mathcal O}_{{Z_s}}, {\mathcal D}^{\bullet}_{Z})$ and ${{ i_s}}_{*}{\mathcal  H}^{-m+1}({\mathcal D}^{\bullet}_{Z_s})={{ i_s}}_{*}{\omega}^{m-1}_{Z_s}$. \vspace{1mm}

\noindent To verify its injectivity, it suffices to examine the low-degree terms of the spectral sequence
$${\rm E}^{i,j}_{2}={\mathcal E}xt^{i}( {{i_s}}_{*}{\mathcal O}_{{Z_s}}, {\mathcal H}^{j}({\mathcal D}^{\bullet}_{Z}[-m]))\Longrightarrow{\mathcal E}xt^{i+j}( {{i_s}}_{*}{\mathcal O}_{{Z_s}}, {\mathcal D}^{\bullet}_{Z}[-m]) $$

$$\xymatrix{{\rm E}^{1,0}_{2}\eq[d]\ar@{^{(}->}[r]&{\mathcal E}xt^{1}( {{i_s}}_{*}{\mathcal O}_{{Z_s}}, {\mathcal D}^{\bullet}_{Z}[-m])\eq[d]\\
{\mathcal E}xt^{1}({{ i_s}}_{*}{\mathcal O}_{{Z_s}}, {\omega}^{m}_{Z})\ar@{^{(}->}[r]&{{ i_s}}_{*}{\omega}^{m-1}_{Z_s} }$$
Note, in passing, that if $Z$ is Cohen-Macaulay, examining the low-degree terms of this spectral sequence shows that the injective morphism ${\mathcal E}xt^{1}(i_{*}{\mathcal O}_{Z_0}, \omega^{m}_{Z})\rightarrow {i_{0}}_{*}\omega^{m-1}_{{ Z}_{0}}$ is bijective, and we have a commutative diagram 
$$\xymatrix{{i_{0}}_{*}\omega^{m-1}_{{ Z}_{0}}\eq[r]&{\mathcal E}xt^{1}(i_{*}{\mathcal O}_{Z_0}, \omega^{m}_{Z})\\
{i_{0}}_{*}{\mathcal L}^{m-1}_{Z_0}\ar@{^{(}->}[u]\ar@{^{(}->}[r]&{\mathcal E}xt^{1}({i_{0}}_{*}{\mathcal O}_{Z_0}, {\mathcal L}^{m}_{Z})\ar[u]}$$
in which the second vertical arrow is not necessarily injective!
\vspace{1mm}

\indent $\bullet$ For the second, consider the commutative diagram
$$\xymatrix{{\widetilde{Z}'_{s}}\ar[rd]_{\phi"}\ar[r]_{\tilde{i'}_{s}}\ar@/^1pc/[rr]^{\sigma_{\!\! s}} &{\widetilde{Z}_{s}}\ar[r]_{\tilde{i}_{s}}\ar[d]_{\phi'}&{\widetilde Z}\ar[d]^{\phi}\\
&Z_{s}\ar[r]_{i_s}&Z}$$
induced by the desingularization of $Z$, where $\widetilde{Z}_s$ denotes the total preimage of ${Z_s}$. We may assume the strict preimage $\widetilde{Z}'_{s}$, which is a union of irreducible components of this total preimage, is smooth. \vspace{1mm} 

\noindent We construct the morphism ${{ i_s}}_{*}({\mathcal L}^{m-1}_{Z_s})\rightarrow{\mathcal E}xt^{1}({{ i_s}}_{*}{\mathcal O}_{{Z_s}}, {\mathcal L}^{m}_{Z})$ using elementary operations and compositions summarized in the diagram $$\xymatrix{
{i_{s}}_{*}{\mathcal L}^{m}_{Z_{s}}\ar[dd]\eq[r]&{i_{s}}_{*}{\rm I}\!{\rm R}\phi"_{\!\! *}\Omega^{m-1}_{\widetilde Z'_{s}}\eq[r]&
{\rm I}\!{\rm R}\phi_{*}{\rm I}\!{\rm R}{\mathcal H}om({{\sigma_{\!\! s}}_{*}}{\mathcal O}_{{\widetilde Z'}_{s}}, \Omega^{m}_{\widetilde Z})[1]\ar[d]\\
&&{\rm I}\!{\rm R}{\mathcal H}om({\rm I}\!{\rm R}{{\phi}_{*}}{{{\sigma_{s}}}_{*}}{\mathcal O}_{{\widetilde Z}_{s}}, {\rm I}\!{\rm R}{{\phi}_{*}}{\Omega}^{m}_{\widetilde Z})[1]\ar[d]\\
{\rm I}\!{\rm R}{\mathcal H}om({{i_{s}}_{*}}{{\mathcal O}_{{Z}_{s}}}, {\mathcal L}^{m}_{Z})[1]&&{\rm I}\!{\rm R}{\mathcal H}om({{i_{s}}_{*}}{\rm I}\!{\rm R}{{\phi"}_{\!*}}{{\mathcal O}_{{\widetilde Z'}_{s}}}, {\mathcal L}^{m}_{Z})[1]\ar[ll]}$$
Recall that since $\widetilde{Z'}_{s}$ is smooth of codimension $1$ in $\widetilde Z$, we have
$${\rm I}\!{\rm R}{\mathcal H}om({{ \sigma}_{\!\! s}}_{*}{\mathcal O}_{{\widetilde Z'}_s}, \Omega^{m}_{\widetilde Z})[1]\simeq {\mathcal E}xt^{1}({{ \sigma_{\!\! s}}}_{*}{\mathcal O}_{{\widetilde{Z'}_s}}, {\Omega}^{m}_{\widetilde Z})\simeq {{ \sigma_{\!\! s}}}_{*}\Omega^{m-1}_{{\widetilde{Z'}_s}}.$$
In particular, taking the degree $0$ cohomology of the first vertical arrow yields the desired morphism
$${{ i_s}}_{*}({\mathcal L}^{m-1}_{Z_s})\rightarrow{\mathcal E}xt^{1}({{ i_s}}_{*}{\mathcal O}_{{Z_s}}, {\mathcal L}^{m}_{Z}).$$
We then complete the previous two-line commutative diagram to obtain 
$$\xymatrix{&&&{i_s}_{*}{\mathcal L}^{m}_{ {Z}_{s}}\ar[d]_{\alpha}\ar@/^4pc/[ddd]^{\sigma}\\
0\ar[r]&{\mathcal L}^{m}_{ Z}\ar@{_{(}->}[d]\ar[r]^{\bullet {f}}&{\mathcal L}^{m}_{ Z, \tilde{f}}\ar@{_{(}->}[d]_{\alpha'}\ar[r]&{\mathcal E}xt^{1}({{ i_s}}_{*}{\mathcal O}_{{Z_s}}, {\mathcal L}^{m}_{Z})\ar[d]_{\beta}\ar[r]&0\\
0\ar[r]&\omega^{m}_{Z}\ar[r]^{\bullet f}&\omega^{m}_{Z, \tilde{f}}\ar[r]&{\mathcal E}xt^{1}({{ i_s}}_{*}{\mathcal O}_{{Z_s}}, {\omega}^{m}_{Z})\ar@{_{(}->}[d]^{\gamma}\ar[r]&0\\
&&&{{ i_s}}_{*}{\omega}^{m}_{Z_s}}$$
and deduce that the bijectivity of the natural injection $\sigma$ forces the surjectivity of $\gamma$ and hence its bijectivity. This in turn implies the surjectivity of $\beta$ and the injectivity of $\alpha$. But this latter surjectivity entails that of $\alpha'$ and thus its bijectivity, which in turn ensures the bijectivity of the first vertical arrow—i.e., the bijectivity of the natural injection of ${\mathcal L}^{m}_{ Z}$ into $\omega^{m}_{Z}$—proving that $Z$ is of type $\mathfrak{F.R}^{*}\,\blacksquare$ \vspace{1mm}

\noindent {\bf{(ii) Smooth base of dimension $r>1$ and general case.}}

 \vspace{1mm}  
 
 \noindent
   We proceed by induction on the dimension $r$ of $S$, assuming the result holds for any base of dimension strictly less than $r$. The idea is to consider a smooth subspace of dimension $r-1$, perform a base change, apply the induction hypothesis to the resulting total space, and then view this subspace itself as a fiber at the origin (for example) of a family of hypersurfaces. These are classical arguments often used in such situations (cf. \cite{AS84} or \cite{Schn}). \par\noindent
 Consider a point $s_{0}$ in $S$ around which a local coordinate system $(s_{1},\cdots, s_{r})$ is given on a sufficiently small Stein open set $V$. Let $S':=\{s\in V: s_{1}(s)=\cdots= s_{r-1}(s)=0\}$. Since $\pi$ is open and $s_{r}$ is a parameter near $s_{0}$ (thus inducing an active element on the local ring ${\mathcal O}_{S,s_{0}}$), $\tilde{s}_{r}:=\pi^{*}(s_{r})$ is again an active element on the local ring of $Z$ near any point with image $s_{0}$, by (\cite{Fi}, {\bf{ Lemma 2}}, p. 144). \vspace{1mm}
 \noindent 
 In the Cartesian diagram of base change 
 $$\xymatrix{Z':=Z\times_{S} S'\ar[r]^{\,\,\,\,\,\,\,\,\,\sigma'}\ar[d]_{\pi'}&Z\ar[d]^{\pi}\\
 S'\ar[r]_{\sigma}&S}$$
 Since $\pi'$ is still flat with reduced fibers (analytically isomorphic to the fibers of $\pi$) and of type $\mathfrak{F.R}$, the induction hypothesis applied to $\pi'$ guarantees that $Z'$ is also of type $\mathfrak{F.R}$. \vspace{1mm}

 \noindent  We then consider the canonical projection $p:{\Bbb C}^{r}\rightarrow {\Bbb C}^{r-1}$ sending $(z_{1},\cdots, z_{r})$ to $(z_{1},\cdots, z_{r-1})$ and a local coordinate change isomorphism $h:V\rightarrow{\Bbb C}^{r}$. Then, if $S":=p(h(V))$, the composite morphism $\psi:Z\rightarrow S"$ is still geometrically flat, being open with fibers of dimension $n+1$ over the smooth $S"$. Moreover, after possibly shrinking the data, we obtain the commutative diagram
$$\xymatrix{Z'\ar@/_2pc/[dd]_{\psi'}\ar[r]^{\sigma'}\ar[d]_{\pi'}&Z\ar[d]^{\pi}\ar@/^2pc/[dd]^{\psi}\\
 S'\ar[d]\ar@{^{(}->}[r]_{\sigma}&S\ar[d]\\
 \{0\}\ar@{^{(}->}[r]&S"}$$
 which shows that the fiber over $\{0\}$ of this $\psi$ is given by the analytic space $Z\times_{S"} \{0\}$. Now, the latter is analytically isomorphic to $Z':=Z\times_{S} S'$, which has singularities of type ${\mathfrak{F.R}}^{*}$. We then conclude by applying the induction hypothesis. From this case, we deduce the general case of singular base of type ${\mathfrak{F.R}}^{*}$ by applying \propositionref{P8}$\,\blacksquare$ \vspace{2mm}

 \noindent 
\centerline{\bf{b) Method 2: By flattening the morphism.} }\vspace{1mm}

\noindent 
Recall that a {\emph{local blow-up}} of an analytic space $Z$ consists of an open subset $U$ of $Z$, a closed analytic subset $Y$ of $U$, and the morphism $\sigma: V \rightarrow Z$ composed of the blow-up of $Y$ in $U$ and the inclusion morphism of $U$ into $Z$. A sequence of local blow-ups of $Z$ is a finite system $(U_{\alpha}, Z_{\alpha}, \sigma_{\alpha})_{1 \leq \alpha \leq r}$ where each $\sigma_{\alpha}$ is a local blow-up and $Z_0 = Z$. \vspace{1mm} \noindent Note that a local blow-up is generally never proper.\vspace{1mm}

\noindent The known flattening theorems are primarily of two types: algebraic and geometric. In the following statements, the first two are of the first category, and the latter is of the second.

\Th{}{}{\color{blue}{Global Flattening}}\label{Aplat1}\cite{H1}:\vspace{1mm} \noindent Let $\pi: X \rightarrow Y$ be a morphism of complex analytic spaces with $Y$ reduced. Let $\mathcal{F}$ be a coherent sheaf on $X$. Then, there exists a proper modification $\sigma: \widetilde{Y} \rightarrow Y$ composed of a finite sequence of blow-ups such that the strict transform $(\widetilde{\mathcal{F}}, \widetilde{\pi})$ of $(\mathcal{F}, \pi)$ is flat everywhere.\rm

\Th{}{}{\color{blue}{Local Flattening}}\label{Aplat2}\cite{HLT}: \vspace{1mm} \noindent Let $\pi: X \rightarrow Y$ be a morphism of complex analytic spaces with $Y$ reduced. Let $y \in Y$ and $L$ be a compact subset of the fiber $\pi^{-1}(y)$. Then, there exists a finite number of morphisms $\sigma_{\alpha}: Y_{\alpha} \rightarrow Y$, each composed of a finite number of local blow-ups whose centers lie in subsets with empty interior, such that:\vspace{1mm}  

\noindent
{\bf(i)} For each $\alpha$, there exists a compact subset $K_{\alpha}$ of $Y_{\alpha}$ such that $\displaystyle{\bigcup_{\alpha}\sigma_{\alpha}(K_{\alpha})}$ is a compact neighborhood of $y$ in $Y$.\vspace{1mm} 

\noindent
{\bf(ii)} For each $\alpha$, the strict transform $\pi_{\alpha}: X_{\alpha} \rightarrow Y_{\alpha}$ of $\pi$ by $\sigma_{\alpha}$ is flat at points corresponding to $L$.\rm\vspace{2mm}

\noindent
\Th{}{}{\color{blue}{Local Geometric Flattening}}\label{Aplat3}\cite{Si}: \vspace{1mm} 
\noindent Let $\pi: X \rightarrow Y$ be a generically open morphism of reduced complex analytic spaces, and let ${\rm K}$ be a compact subset of $X$. Then, there exists a normal complex space $Y'$, a holomorphic map $\sigma: Y' \rightarrow Y$, and a relatively compact open neighborhood $V$ of ${\rm K}$ such that:\vspace{1mm}

\noindent
{\bf(i)} If $X'$ is the subspace of $Y' \times_{Y} V$ consisting of irreducible components over which the transform of $\pi$ (by $\sigma$) remains generically open onto $Y'$, then the strict transform $\pi': X' \rightarrow Y'$ is open and surjective.\vspace{1mm}

\noindent
{\bf(ii)} The induced map $\sigma': X' \rightarrow V$ is semi-proper and surjective.
\vspace{1mm}

\noindent
{\bf(iii)} $\sigma$ is generically open and generically finite.\rm
\rm\vspace{2mm}

\noindent
We have the diagram
$$\xymatrix{X'\ar@/^1pc/[rr]^{\sigma'}\ar[rd]_{\pi'}\ar@{^{(}->}[r]&Y' \times_{Y} V\ar[d]\ar[r]&X\ar[d]^{\pi}\\
&Y'\ar[r]_{\sigma}&Y}$$

Note that $\sigma$ is neither proper nor even semi-proper but preserves density by inverse image. Moreover, if ${\rm E}$ is the degeneracy locus of $\pi$ (the union of subsets consisting of points where the fibers have dimension strictly greater than $\dim X - \dim Y$) and ${\rm N}(Y)$ the non-normal locus of $Y$, then $Y'$ is identified with the fiber product $Y' \times_{Y} V$ outside $\Sigma := \sigma^{-1}(\overline{\pi({\rm E} \cap V)} \cup {\rm N}(Y))$,
$$X' \setminus \pi'^{-1}(\Sigma) \simeq (Y' \setminus \Sigma) \times_{Y} V$$ 
with $\sigma: Y' \setminus \Sigma \rightarrow Y$ (resp. $\pi': X' \setminus \pi'^{-1}(\Sigma) \rightarrow X$) open.\vspace{1mm}

\noindent
This theorem follows from \theoremref{Aplat2}, as noted by Sibert in remark ({\bf(6), p.274}). \vspace{2mm}

\noindent
Let us return to the proof of our theorem, for which we must distinguish the proper case from the non-proper case. \vspace{1mm}

\noindent {\bf{(i)The Proper Case.}}\vspace{1mm} 

\noindent This is an almost immediate corollary of \corollaryref{C11}, as it suffices to flatten the morphism by applying Hironaka's global flattening theorem \theoremref{Aplat1} to obtain the diagram:
 $$\xymatrix{{Z_{1}}\ar[r]^{\Theta}\ar[d]_{\pi_{1}}&Z\ar[d]^{\pi}\\
{S_{1}}\ar[r]_{\eta_{1}}&S}$$ 
noting that the strict preimage ${Z_{1}}$ coincides with the total preimage since $\pi$ is open. Then, since $S$ is of type ${\mathfrak{F.R}}^{*}$ (and hence ${\mathfrak{F.R}}$), ${S}_{1}$ is also of type ${\mathfrak{F.R}}$ by virtue of \lemmaref{L2} or \lemmaref{L3}. \vspace{1mm}

\noindent We may then either normalize or desingularize the base $S_1$ to obtain the diagram:
\[\xymatrix{{Z_{2}}\ar@/^1pc/[rr] \ar[r]_{\Theta_{2}}\ar[d]_{\pi_{2}}&{Z_{1}}\ar[r]_{\Theta_{1}}\ar[d]_{\pi_{1}}&Z\ar[d]^{\pi}\\
{S_{2}}\ar@/_1pc/[rr]\ar[r]^{\eta_{2}}&{S_{1}}\ar[r]^{\eta_{1}}&S}\]
where the morphism $\pi_{2}$ is still flat, with fibers and base of type ${\mathfrak{F.R}}^{*}$. If we normalize, we apply \corollaryref{C11}, and if we desingularize \cite{AS84}, we deduce that ${Z_{2}}$ is of type ${\mathfrak{F.R}}^{*}$. \vspace{1mm}  

\noindent 
{\bf{(ii) The Non-Proper Case.}} Here again, we may assume $S$ is normal without loss of generality by applying \propositionref{P8} to the base change given by the normalization of $S$.\vspace{1mm}

\noindent 
We then use either the local flattening theorem \theoremref{Aplat2} or the geometric flattening theorem \theoremref{Aplat3} from \cite{Si}. We begin by noting the following:\vspace{1mm}
 
 \indent $\bullet$ {\bf{Stability of the ${\mathfrak{F.R}}$ Condition Under Local Blowup.}} Let $\phi:S'\rightarrow S$ be a local blowup, i.e., the composition of an open embedding $j:V\hookrightarrow S$ and a blowup $\sigma: S'\rightarrow V$ whose center lies in a nowhere dense subset of $S$. Clearly, a local blowup is never proper in general. \vspace{1mm}

 \noindent 
The stability under open restriction (which is a flat morphism) and under blowup (which is a proper modification to which we apply \lemmaref{L2} or \lemmaref{L3}) will give us stability under an arbitrary local blowup.\vspace{1mm}

\noindent Note that the inverse image ${\mathcal L}^{m}_{S}\rightarrow\phi_{*}{\mathcal L}^{m}_{S'}$ naturally factors as ${\mathcal L}^{m}_{S}\rightarrow j_{*}{\mathcal L}^{m}_{V} \rightarrow{\phi}_{*}{\mathcal L}^{m}_{S'}$, and that $$\phi_{*}{\mathcal L}^{m}_{S'}=j_{*}{\mathcal L}^{m}_{V}=j_{*}j^{*}{\mathcal L}^{m}_{S}\simeq {\mathcal L}^{m}_{S}$$
since the sheaf ${\mathcal L}^{m}_{S}$ has depth at least two, and the center of local blowups may be chosen within the singular locus of codimension at least two in $S$.\vspace{1mm}

 \noindent
\indent $\bullet$ {\bf{Preimage Under a Local Blowup.}} With the same notation, we have the commutative diagram:
$$\xymatrix{Z"\ar[r]^{\sigma'}\ar[d]_{\pi"}&Z'\ar[d]^{\pi'}\ar[r]^{j'}&Z\ar[d]^{\pi}\\
{S'}\ar[r]_{\sigma}&V\ar[r]_{j}&S}$$ 
where $Z':=Z\times_{S}V$ and $Z":=Z'\times_{V}S'\simeq Z\times_{S} S' $. Since $\pi$ is open with fibers of constant pure dimension $n$, the same holds for the morphisms $\pi'$ and $\pi"$, as these properties are preserved under any base change. Moreover, as already mentioned, the total and strict preimages coincide under these base changes, which is crucial for our purpose since the fibers will always remain analytically isomorphic (in general, on the strict preimage, the new fibers are isomorphic to the old ones only generically over the base).\vspace{1mm}

\indent
$\bullet$ This type of diagram will simply be denoted, for any $\alpha$, as:
$$\xymatrix{{Z_{\alpha}}\ar[r]^{\sigma'_{\alpha}}\ar[d]_{\pi_{\alpha}}&Z\ar[d]^{\pi}\\
{S_{\alpha}}\ar[r]_{\sigma_{\alpha}}&S}$$ 
following the notation of \theoremref{Aplat2}. 
The idea is to proceed by induction on the number of local blowups constituting $\sigma_{\alpha}$ to reduce to the case of a single blowup. To do this, we revisit the arguments developed in \cite{Par} point by point, allowing us to essentially reduce to the case of a standard blowup by localizing sufficiently, taking for the compact $L:=\{z\}$ a point in the fiber $\pi^{-1}(\pi(z))$. \vspace{1mm}  

\noindent
We consider an open neighborhood $U$ of $s:=\pi(z)$ and a compact subneighborhood $K$ of $s$ in $U$. Since the sheaves ${\mathcal L}^{m}_{Z}$ and $\omega^{m}_{Z}$ are compatible with open inclusions, we first replace the data $(\pi,Z,S)$ with $(\pi|_{\pi^{-1}(U)},\pi^{-1}(U), U)$, noting that the conditions of \theoremref{Aplat2} are satisfied for $L:=\{z\}$, $\sigma_{\alpha}:=\sigma_{\alpha}|_{\sigma_{\alpha}^{-1}(U)}$, $K_{\alpha}:=K_{\alpha}\cap \sigma_{\alpha}^{-1}(K)$. The centers of the local blowups $\sigma_{\alpha}|_{\sigma_{\alpha}^{-1}(U)}$ are open subsets of the centers of the $\sigma_{\alpha}$. Moreover, for a blowup $\sigma:S'\rightarrow S$, the strict preimages of the $\sigma_{\alpha}$ are compositions of local blowups with rare centers contained in or equal to certain centers of the $\sigma_{\alpha}$.\vspace{1mm}  
 
 \noindent
If $\nu:S'\rightarrow S$ is a given blowup, we obtain, via base changes, the following diagrams:
 $$\xymatrix{{Z'_{\alpha}}\ar[rd]\ar@[red][dd]_{\color{red}{\pi'_{\alpha}}}\ar@[red][rr]^{\nu'_{\alpha}}&&{Z_{\alpha}}\ar@[red][dd]^{\color{red}{{\pi_{\alpha}}}}\ar[rd]&\\
 &Z'\ar@[blue][rr]^{\!\!\!\!\!\!\!\!\!\!\!\!\!\!\nu'}\ar@[blue][dd]_{\pi'}&&Z\ar@[blue][dd]^{\pi}\\
 {S'_{\alpha}}\ar@[red][rr]^{\nu_{\alpha}}\ar[rd]_{\sigma'_{\alpha}}&&{S_{\alpha}}\ar[rd]^{\sigma_{\alpha}}&\\
 &{S'}\ar@[blue][rr]_{\nu}&&S}$$
where $\nu$, $\nu'$, $\nu_{\alpha}$, and $\nu'_{\alpha}$ are proper morphisms, the pair $(\sigma'_{\alpha}, \nu_{\alpha}^{-1}(K_{\alpha}))$ satisfies the conditions of the local flattening theorem for the morphism $\pi':Z'\rightarrow S'$, and $L':=\nu'^{-1}(L)$.\vspace{1mm} 

\noindent Of course, this remains entirely valid if $\nu$ is a local blowup. We then begin the induction on the number of local blowups, essentially reducing—after sufficient localization—to the case of a single blowup, where we can apply \propositionref{P7}. \vspace{1mm}

\noindent
We may bypass the inductive reasoning and its constraints by adopting the formulation of \theoremref{Aplat3}. For a given $z_0$ in $Z$, we consider a relatively compact open neighborhood $U$ with closure $K:=\overline{U}$. Then, denoting by $V$ the relatively compact neighborhood containing $K$ provided by the theorem, we have the commutative diagram:
$$\xymatrix{S'\times_{S}{Z'}\ar[r]^{\sigma'}\ar[d]_{\pi'}&Z\ar[d]^{\pi}\\
{S'}\ar[r]_{\sigma}&S}$$ 
where $Z':=V\times_{S}S'$, $S'$ is normal, and $\sigma$ and $\sigma'$ are open. \vspace{1mm}

\noindent We address the question directly by revisiting the diagram of local parametrizations (again with simplified notation):
$$\xymatrix{{Z'}\ar[r]^{\sigma'}\ar[d]_{f'}&Z\ar[d]^{f}\\
{ S'}\times U\ar[r]_{{\sigma}}&S\times U}$$ 
where, by slight abuse of notation, we denote the induced morphisms by the same letters $\sigma$ and $\sigma'$. We then easily adapt \propositionref{P8} to deduce that if $Z'$ is of type ${\mathfrak{F.R}}^{*}$, so is $Z$. \vspace{1mm}

 \noindent 
More precisely, if $\xi$ is a section of the sheaf $\omega^{m}_{Z}$, we have:
$${\rm Tr}_{f'}(\sigma'^{*}(\xi))=\sigma^{*}{\rm Tr}_{f}(\xi)$$
Then, since by hypothesis ${\rm Tr}_{f}(\xi)$ defines a section of the sheaf ${\mathcal L}^{m}_{S\times U}$, it follows that ${\rm Tr}_{f'}(\sigma'^{*}(\xi))$ defines a section of the sheaf ${\mathcal L}^{m}_{S'\times U}$. Since local blowups, like open morphisms, preserve density under inverse images, the meromorphic form $\sigma'^{*}(\xi))$ has its poles in a nowhere dense closed subset contained in the nowhere dense closed subset $\sigma'^{-1}({\rm Sing}(Z))$. \vspace{1mm}

 \noindent As in \propositionref{P8}, we show that only the poles contained in $\sigma'^{-1}({\rm Sing}(Z))\cap {\rm Sing}(Z')$ matter. Moreover, if $S'\rightarrow V$ is a local parametrization of $S'$, it yields, by composition, the parametrization:
$$\xymatrix{Z'\ar@/^1pc/[rr]^{h}\ar[r]_{f'}&S'\times U\ar[r]_{g}&V\times U}$$
for which we have:
$${\rm Tr}_{h}(\sigma'^{*}(\xi))={\rm Tr}_{g}({\rm Tr}_{f'}(\sigma'^{*}\xi))$$
which is holomorphic since ${\rm Tr}_{f'}(\sigma'^{*}\xi)$ is a section of the sheaf ${\mathcal L}^{m}_{S'\times U}$, a subsheaf of $\omega^{m}_{S'\times U}$. \vspace{1mm}

 \noindent Thus, the form $\sigma'^{*}(\xi)$ defines a section of the sheaf $\omega^{m}_{Z'}$. But since $Z'$ is of type ${\mathfrak{F.R}}^{*}$, this shows that $\sigma'^{*}(\xi)$ is a section of the sheaf ${\mathcal L}^{m}_{Z'}$. We immediately deduce that $\xi$ must necessarily be a section of the sheaf ${\mathcal L}^{m}_{Z}$ and, consequently,
 $${\mathcal L}^{m}_{Z}\simeq\omega^{m}_{Z}$$
meaning $Z$ is of type ${\mathfrak{F.R}}^{*}\,\blacksquare$\vspace{1mm}
\section{\color{blue}{Transfer of properties ${\mathfrak{F.R}}$ from Fibers to the Total Space:  \theoremref{Th3}.}}
\vspace{2mm}

\noindent
Let $\pi:Z\rightarrow S$ be a universally $n$-equidimensional morphism whose base and fibers have singularities of type $\mathfrak{F.R}$, and we aim to show that $Z$ also possesses this property.  
To do this, we proceed in several steps, starting by noting that since the constancy of fiber dimension and openness are properties invariant under base change, we may assume $S$ is smooth or normal by virtue of \propositionref{P7} and \propositionref{P8}.  

\subsection{{The case of a reduced morphism $\pi$.}}\vspace{2mm}  

\noindent Since the property of being reduced for the fibers is also invariant under base change and we may assume the base is smooth, we reduce to the study of a flat morphism with reduced fibers over a smooth base (since any open morphism over a smooth base with reduced fibers is necessarily flat; cf. \cite{GR64}, \cite{Kan}, \cite{Fi} p.157). In this case, the source is also reduced\footnote{Moreover, since the fibers are reduced, this flat morphism is of type ${\rm S}_{1}$. The set of points of $Z$ that are Cohen-Macaulay in the fibers ${\rm CM}(\pi)$ is a dense open subset of $Z$ whose complement is a closed set with empty interior of codimension at least two fiberwise, as shown in Lemmas {\bf(5.2.2)} and {\bf(5.2.5)} of \cite{A.L}.}.  
\vspace{1mm}  

\noindent  
The following cases are immediate corollaries of \theoremref{Th2}:\vspace{1mm}  

\indent $\star$ If $\pi$ is equinormalizable or admits a simultaneous normalization (cf. \cite{ChLi}, \cite{G}, \cite{LCTr}), then it is induced by the normalization of $Z$, $\nu:\overline{Z}\rightarrow Z$, with the commutative diagram  
$$\xymatrix{\overline{Z}\ar[rr]^{\nu}\ar[rd]_{\bar\pi}&&Z\ar[ld]^{\pi}\\  
&S&}$$  
in which the fiber ${\bar\pi}^{-1}(s)$ is the normalization of the fiber ${\pi}^{-1}(s)$. \vspace{1mm}  

\noindent  
Then, we have  
$${\pi}^{-1}(s)\,{\rm of\, type}\,{\mathfrak{F.R}}\Longleftrightarrow {\bar\pi}^{-1}(s)\,{\rm of\, type}\,{\mathfrak{F.R}}^{*}\Longrightarrow\, \overline{Z}\,{\rm of\, type}\,{\mathfrak{F.R}}^{*}\Longleftrightarrow Z\,{\rm of\, type}\,{\mathfrak{F.R}}$$  

\indent $\star$ If for the normalization $\nu:\overline{Z}\rightarrow Z$, $\overline{Z}_{s}:=\nu^{-1}(\pi^{-1}(s)))$ is reduced with a singular locus of codimension at least two in the fiber for all $s$, then $\bar\pi$ is of type ${\mathfrak{F.R}}^{*}$. This is because the normalization morphism of $Z_s$ necessarily factors through $\overline{Z}_{s}$ (cf. \cite{G}, {\bf{Lemma 2.48}, p. 343}), which is then of type ${\mathfrak{F.R}}$ and thus ${\mathfrak{F.R}}^{*}$ by hypothesis on the codimension of the singular locus. \vspace{1mm}  

\indent We draw the reader's attention to the fact that, in general, the morphism $\bar\pi$ is not reduced but remains open and equidimensional over any normal base (in fact, it is geometrically flat, meaning it describes an analytic family of cycles).  
\vspace{2mm}  

\noindent 
\centerline{{\bf{(1) $\pi:Z\rightarrow S$ reduced of type ${\mathfrak{F.R}}$ with smooth base of dimension $1$:}}}\vspace{2mm}  

\noindent {\bf(a) Case $Z$ normal:}  
After localization if necessary, we essentially reduce to the case where $S$ is the unit disk ${\Bbb D}$, $Z_{s}$ is of type ${\mathfrak{F.R}}^{*}$ for all $s\not=0$, and $Z_{0}$ is of type ${\mathfrak{F.R}}$. \theoremref{Th2} ensures that $Z$ is of type ${\mathfrak{F.R}}^{*}$ on $Z\setminus Z_{0}$, and we may assume all fibers are smooth outside the central fiber. \vspace{1mm}  

\noindent  
Since $Z$ is normal, the intersection ${\rm Sing}(Z)\cap Z_{0}$ has empty interior in $Z_0$ (or equivalently, ${\rm Reg}(Z)\cap Z_{0}$ is a dense open subset of $Z_0$). Given our hypotheses, we only need to study what happens near the points of ${\rm Sing}(Z)\cap {\rm Sing}(Z_{0})$, since there are no problem outside $Z_0$ and on ${\rm Sing}(Z)\cap {\rm Reg}(Z_0)$, we are back in the situation of \theoremref{Th2}. We may therefore assume that $Z$ is irreducible and ${\rm Sing}(Z) = {\rm Sing}(Z_{0})$ (which is of codimension at least two in $Z$).  
\vspace{1mm}  

\noindent 
The property for $Z$ to be of type ${\mathfrak{F.R}}$ is local in nature, as it amounts to showing that the sheaf ${\mathcal L}^{m}_{Z}$ has depth at least two over $Z$, which can be checked stalk by stalk. Since $\pi:Z\rightarrow {\Bbb D}$ is open with reduced fibers and hence flat, it allows lifting, at the level of local rings, any non-zero divisor to a non-zero divisor\footnote{even in a slightly more general setting, it is always possible, according to (\cite{G}, {\bf{Lemma 2.48}}, p.343), to find a germ of holomorphic function on $Z$ that is a non-zero divisor} as per (\cite{Fi}, {\bf Lemma}, p.157). Moreover, its germ $\pi_{z}$, at each point $z$ of $Z$, defines a non-zero divisor in ${\mathcal O}_{Z,z}$.
Thus, every point of $Z$ has an open neighborhood on which we have a representative of this germ, which we denote $f$, inducing the short exact sequence  
$$\xymatrix{0\ar[r]&{\mathcal O}_{ Z}\ar[r]^{f}&{\mathcal O}_{Z}\ar[r]&{i_{0}}_{*}{\mathcal O}_{{Z}_{0}}\ar[r]&0}$$  
to which we apply the functors ${\mathcal H}om(-,\omega^{m}_{Z})$ and ${\mathcal H}om(-, {\mathcal L}^{m}_{Z})$ to obtain the commutative diagram with exact rows:  
$$\xymatrix{0\ar[r]&{\mathcal L}^{m}_{Z}\ar@{^{(}->}[d]\ar[r]^{f}&{\mathcal L}^{m}_{Z}\ar@{^{(}->}[d]\ar[r]&{\mathcal E}xt^{1}({{i_{0}}_{*}}{\mathcal O}_{Z_0}, {\mathcal L}^{m}_{Z})\ar[d]\ar[r]&0\\  
0\ar[r]&\omega^{m}_{Z}\ar[r]^{ f}&\omega^{m}_{Z}\ar[r]&{\mathcal E}xt^{1}({{i_{0}}_{*}}{\mathcal O}_{Z_0}, \omega^{m}_{Z})\ar[r]&0}$$  
where the vertical arrows are generically bijective and naturally injective for the first two. We will show that the same holds for the last one by verifying the absence of ${\mathcal O}_{Z_0}$-torsion for the sheaf ${\mathcal E}xt^{1}(i_{*}{\mathcal O}_{Z_0}, {\mathcal L}^{m}_{Z})$. \vspace{1mm}  

\noindent To see this, consider the diagram  
$$\xymatrix{\widetilde{Z'_{0}}\ar@/^1pc/[rr]^{\tilde{i'}_{0}}\ar@{^{(}->}[r]_{\sigma}\ar[rd]_{\phi'_{0}}&{\widetilde{Z}}_{0}\ar@{^{(}->}[r]_{\tilde{i}_{0}}\ar[d]^{\phi_{0}}&\widetilde{Z}\ar[d]^{\phi}\\  
&Z_{0}\ar@{^{(}->}[r]^{{i}_{0}}&Z}$$  
where $\phi:\widetilde{Z}\rightarrow Z$ is a desingularization such that the total preimage $\widetilde{Z}_{0}$ is a simple normal crossing divisor and the strict preimage $\widetilde{Z'_0}$ is smooth (and a finite union of certain irreducible components of $\widetilde{Z}_{0}$). The horizontal arrows are embeddings, and the others are natural restrictions of $\phi$. \vspace{1mm}  

\noindent  
Abusively denoting $\tilde{f}:=\phi^{*}(f)$, we have the short exact sequence  
$$\xymatrix{0\ar[r]&{\mathcal O}_{\tilde Z}\ar[r]^{\tilde f}&{\mathcal O}_{\tilde Z}\ar[r]&{\tilde{i_{0}}_{*}}{\mathcal O}_{{\widetilde Z}_{0}}\ar[r]&0}$$  
to which we apply the functor ${\mathcal H}om(-,\Omega^{m}_{\widetilde Z})$ to obtain the short exact sequence  
$$\xymatrix{0\ar[r]&\Omega^{m}_{\widetilde Z}\ar[r]^{\tilde f}&\Omega^{m}_{\widetilde Z}\ar[r]&{\tilde{i_{0}}_{*}}{\omega^{m-1}_{{\widetilde Z}_{0}}}\ar[r]&0}$$  
since $\displaystyle{{\tilde{i_{0}}_{*}}\omega^{m-1}_{{\tilde Z}_{0}}\simeq{\mathcal E}xt^{1}({\tilde{i_{0}}}_{*}{\mathcal O}_{{\widetilde Z}_{0}}, \Omega^{m}_{\widetilde Z}) }$. Hence, thanks to \cite{GR70}, the exact sequence  
$$\xymatrix{0\ar[r]&{\mathcal L}^{m}_{Z}\ar[r]^{f}&{\mathcal L}^{m}_{Z}\ar[r]&{i_{0}}_{*}{\phi_{0}}_{*}(\omega^{m-1}_{{\widetilde Z}_{0}})\ar[r]&0}$$  
\vspace{1mm}  

\noindent  
From this, we deduce the isomorphisms  
$${\tilde{i_{0}}_{*}}{\omega}^{m-1}_{{\widetilde Z}_{0}}\simeq {\Omega}^{m}_{{\widetilde Z}}/{\tilde{f}}.{\Omega}^{m}_{\widetilde Z},\,\,\,\,\,\,\,\phi_{*}({\tilde{i_{0}}_{*}}{\omega}^{m-1}_{{\widetilde Z}_{0}})\simeq{{\mathcal L}^{m}_{{Z}}/f.{\mathcal L}^{m}_{Z}}\simeq {\mathcal E}xt^{1}({{i_{0}}_{*}}{\mathcal O}_{Z_0}, {\mathcal L}^{m}_{Z})$$  
which incidentally shows that the sheaves ${{\mathcal L}^{m}_{{Z}}/f.{\mathcal L}^{m}_{Z}}$ and ${\mathcal E}xt^{1}({{i_{0}}_{*}}{\mathcal O}_{Z_0}, {\mathcal L}^{m}_{Z})$ are torsion-free over $Z_0$, being the proper direct image of a torsion-free sheaf. It immediately follows that the last vertical arrow is injective due to its generic bijectivity and the absence of torsion.  
\vspace{1mm}  

\noindent Moreover, we have two natural morphisms
$${{i_{0}}_{*}}{\mathcal L}^{m-1}_{{Z}_{0}}\rightarrow{\mathcal E}xt^{1}({{i_{0}}_{*}}{\mathcal O}_{Z_0}, {\mathcal L}^{m}_{Z})$$  
$${\mathcal E}xt^{1}({{i_{0}}_{*}}{\mathcal O}_{Z_0}, {\mathcal L}^{m}_{Z})\rightarrow {{i_{0}}_{*}}{\mathcal L}^{m-1}_{{Z}_{0}}$$  

\indent $\bullet$ The first one, already mentioned in \S 5.1 {\bf(i)}, is injective and can be deduced, in this new point of view, from the direct image or trace at the level of dualizing complexes\footnote{More simply, one can use the pull back ${\mathcal L}^{m-1}_{{Z}_0}\rightarrow{\phi_{0}}_{*}{\mathcal L}^{m-1}_{\widetilde {Z}_0}$ which is composed with the natural morphism ${\mathcal L}^{m-1}_{\widetilde {Z}_0}\rightarrow \omega^{m-1}_{\widetilde {Z}_0}$ and apply the functor ${{i_{0}}_{*}}$.} $\sigma_{*}{\mathcal D}^{\bullet}_{{\widetilde Z'}_{0}}\rightarrow {\mathcal D}^{\bullet}_{{\widetilde Z}_{0}}$
whose degree $-(m-1)$ cohomology yields the morphism of dualizing sheaves (note that ${\widetilde Z'}_{0}$ and ${\widetilde Z}_{0}$ are Gorenstein)\footnote{Note that since $\widetilde{Z}_{0}$ is a simple normal crossing divisor, its normalization $\nu:\overline{\widetilde{Z}_{0}}\rightarrow \widetilde{Z}_{0}$ is smooth. Consequently, the composition $\overline{\widetilde{Z}_{0}}\rightarrow \widetilde{Z}_{0}\rightarrow Z_0$ is a resolution of $Z_0$, and thus $\phi_{*}{\mathcal L}^{m-1}_{\widetilde {Z}_0}={\mathcal L}^{m-1}_{{Z}_0}$.} ${\sigma_{*}}{\omega}^{m-1}_{{\widetilde Z'}_{0}}\rightarrow{\omega}^{m-1}_{{\widetilde Z}_{0}}$  
to which we apply the functor ${{i_{0}}_{*}}\circ {\phi_{0}}_{*}$ to obtain  
$$\xymatrix{{{i_{0}}_{*}}({\phi'_{0}}_{*}{\omega}^{m-1}_{{\widetilde Z'}_{0}})\eq[d]\ar[r]&{{i_{0}}_{*}}({\phi_{0}}_{*}(\omega^{m-1}_{{\widetilde Z}_{0}}))\eq[d]\\  
{{i_{0}}_{*}}{\mathcal L}^{m-1}_{{Z}_{0}}\ar[r]&{\mathcal E}xt^{1}({\mathcal O}_{Z_0}, {\mathcal L}^{m}_{Z})}$$  
\vspace{1mm}  

\noindent
Since $Z$ is normal and $Z_0$ has codimension $1$, this arrow is a generic isomorphism between torsion-free sheaves on $Z_0$, hence injective.
\vspace{1mm}

\indent $\bullet$ To prove the existence of the second morphism, we first observe that for a finite morphism $f:X\rightarrow Y$ of Gorenstein spaces of the same dimension $m$, there exists a natural morphism $f_{*}{\mathcal O}_{X}\rightarrow {\mathcal O}_{Y}$ induced by the usual trace morphism defined at the level of dualizing complexes $f_{*}{\mathcal D}^{\bullet}_{X}\rightarrow {\mathcal D}^{\bullet}_{Y}$; these dualizing complexes have their only non-zero homology in degree $-m$, which is a locally free dualizing sheaf (of rank $1$). In our particular situation, we can invoke more elementary arguments to define the morphism $\sigma_{*}{\mathcal O}_{{\widetilde{Z}'_{0}}}\rightarrow {\mathcal O}_{{\widetilde{Z}_{0}}}$ by reducing to the case where $\widetilde{Z}'_{0}$ is an irreducible component of $\widetilde{Z}_{0}$, in which case the assertion is trivial. Indeed, if ${\mathcal I}$ is the defining ideal of $\widetilde{Z}'_{0}$ in $\widetilde{Z}_{0}$ such that $\sigma_{*}{\mathcal O}_{{\widetilde{Z}'_{0}}}={\mathcal O}_{\widetilde{Z}_{0}}/{\mathcal I}$, we naturally have the short exact sequence
$$\xymatrix{0\ar[r]&{\mathcal I}\ar[r]&{\mathcal O}_{\widetilde{Z}_{0}}\ar[r]_{\beta}&\ar@/_1pc/[l]^{\alpha}{\mathcal O}_{\widetilde{Z}_{0}}/{\mathcal I}\ar[r]&0}$$
and the morphism $\alpha\circ \beta\in {\mathcal H}om({\mathcal O}_{\widetilde{Z}_{0}}, {\mathcal O}_{\widetilde{Z}_{0}})$ is annihilated by ${\mathcal I}$.\vspace{1mm}

\noindent
Then, primarily using duality for a proper morphism, the adjunction formula, and the previous remark, we easily highlight the sequence of natural compositions:
$$\xymatrix{ {\rm I}\!{\rm R}{\mathcal H}om({i_{0}}_{*}{{\mathcal O}_{\widetilde{Z}_{\!\!0}}}, {\mathcal L}^{m}_{{Z}})\ar[dd]\eq[r]&{\rm I}\!{\rm R}{\phi_{*}}{\rm I}\!{\rm R}{\mathcal H}om({\rm I}\!{\rm L}\phi^{*}({i_{0}}_{*}{{\mathcal O}_{{Z}_{\!\!0}}}), {\mathcal D}^{\bullet}_{\widetilde{Z}}[-m])\eq[r]&{\rm I}\!{\rm R}{\phi_{*}}{\tilde{i_{0}}}_{*}{\rm I}\!{\rm R}{\mathcal H}om({\rm I}\!{\rm L}\phi_{0}^{*}({\mathcal O}_{{Z}_{\!\!0}}), {\mathcal D}^{\bullet}_{\widetilde{Z_0}}[-m])\ar[d]^{u}\\
&&{i_{0}}_{*}{\rm I}\!{\rm R}{{\phi_{0}}_{*}}{\rm I}\!{\rm R}{\mathcal H}om({\rm I}\!{\rm L}{\phi_{0}}^{*}({\mathcal O}_{{Z}_{\!\!0}}), \sigma_{*}{\mathcal D}^{\bullet}_{\widetilde{Z'_0}}[-m])\eq[d]\\
{i_{0}}_{*}{\mathcal L}^{m-1}_{Z_0}[-1]\eq[r]&{\rm I}\!{\rm R}{\mathcal H}om({\mathcal O}_{{Z}_{\!\!0}}), {\rm I}\!{\rm R}{{\phi'_{0}}_{*}}{\mathcal D}^{\bullet}_{\widetilde{Z'_0}}[-m])\eq[r]&{i_{0}}_{*}{\rm I}\!{\rm R}{{\phi_{0}}_{*}}\sigma_{*}{\rm I}\!{\rm R}{\mathcal H}om({\rm I}\!{\rm L}{\phi'_{0}}^{*}({\mathcal O}_{{Z}_{\!\!0}}), {\mathcal D}^{\bullet}_{\widetilde{Z'_0}}[-m])}$$
The vertical arrow $u$ is induced by the morphism $\sigma_{*}{\mathcal O}_{\widetilde{Z}'_0}\rightarrow {\mathcal O}_{\widetilde{Z}_0}$ since it suffices to write
$${\mathcal D}^{\bullet}_{\widetilde{Z_0}}={\rm I}\!{\rm R}{\mathcal H}om({\mathcal O}_{\widetilde {Z_{0}}}, {\mathcal D}^{\bullet}_{\widetilde{Z_0}})\rightarrow {\rm I}\!{\rm R}{\mathcal H}om(\sigma_{*}{\mathcal O}_{{\widetilde Z'}_{\!\!0}}, {\mathcal D}^{\bullet}_{\widetilde{Z_0}})=\sigma_{*}{\mathcal D}^{\bullet}_{\widetilde{Z'_0}}$$\vspace{1mm}

\noindent
We thus define, in particular, a morphism
$${\mathcal E}xt^{1}({{i_{0}}_{*}}{\mathcal O}_{Z_0}, {\mathcal L}^{m}_{Z})\rightarrow{{i_{0}}_{*}}{\mathcal L}^{m-1}_{Z_0}$$
Note that these morphisms are constructed so that the composition
$${{i_{0}}_{*}}{\mathcal L}^{m-1}_{Z_0}\rightarrow {\mathcal E}xt^{1}({{i_{0}}_{*}}{\mathcal O}_{Z_0}, {\mathcal L}^{m}_{Z})\rightarrow{{i_{0}}_{*}}{\mathcal L}^{m-1}_{Z_0}$$
is the natural bijection, thus forcing the first arrow to be injective and the second to be surjective. This suggests that the sheaf ${\mathcal L}^{m}_{Z}$ could be a direct factor of ${\mathcal E}xt^{1}({{i_{0}}_{*}}{\mathcal O}_{Z_0}, {\mathcal L}^{m}_{Z})$.\vspace{1mm}

\noindent
However, since the singular locus of $Z$ has codimension at least two, these are generic isomorphisms, and since the sheaf ${\mathcal E}xt^{1}({{i_{0}}_{*}}{\mathcal O}_{Z_0}, {\mathcal L}^{m}_{Z})$ is torsion-free on $Z_0$, the second morphism is necessarily injective and therefore bijective (since it is already injective in the other direction). The initial short exact sequence can then be written as
$$\xymatrix{0\ar[r]&{\mathcal L}^{m}_{Z}\ar[r]^{f}&{\mathcal L}^{m}_{Z}\ar[r]&{i_{0}}_{*}{\mathcal L}^{m-1}_{Z_0}\ar[r]&0}$$
from which we can deduce, thanks to our assumptions, that the sheaf ${\mathcal L}^{m}_{Z}$ must have depth at least two.
\vspace{1mm}

\noindent For simplicity, we will assume $Z$ and $Z_0$ are Stein since the problem is local and allows us to work on Stein open subsets of $Z$ whose restrictions are also Stein on $Z_0$. Then, using the characterization of depth via cohomology with compact support (cf. \cite{Ba1}), we obtain, since ${\mathcal L}^{m-1}_{Z_0}$ has depth at least two by hypothesis, the isomorphism induced by multiplication by $f$:
$${\rm H}^{1}_{c}(Z, {\mathcal L}^{m}_{Z})\simeq{\rm H}^{1}_{c}(Z, {\mathcal L}^{m}_{Z})$$
But since
$${\rm H}^{1}_{c}(Z, {\mathcal L}^{m}_{Z})\simeq {\rm H}^{1}_{c}(\widetilde{Z}, {\Omega}^{m}_{\widetilde Z})$$
Serre duality gives us
$${\rm H}^{1}_{c}(\widetilde{Z}, {\Omega}^{m}_{\widetilde Z})\simeq \biggl({\rm H}^{m-1}(\widetilde{Z}, {\mathcal O}_{\widetilde Z})\biggr)^{'}$$
and the Leray spectral sequence, in turn, provides the isomorphism
$${\rm H}^{m-1}(\widetilde{Z}, {\mathcal O}_{\widetilde Z})\simeq \Gamma(Z, {\rm I}\!{\rm R}^{m-1}\phi_{*}{\mathcal O}_{\widetilde Z}),$$
all these isomorphisms being topological. We thus reduce to an isomorphism between coherent sheaves (of finite type) induced by multiplication by a certain function:
$$ {\rm I}\!{\rm R}^{m-1}\phi_{*}{\mathcal O}_{\widetilde Z}\simeq  {\rm I}\!{\rm R}^{m-1}\phi_{*}{\mathcal O}_{\widetilde Z}$$
which, by virtue of Nakayama's lemma\footnote{We use the variant stating that if ${\mathcal A}$ is a commutative ring, ${\mathcal I}$ an ideal of ${\mathcal A}$ contained in its Jacobson radical, and ${\mathcal M}$ a finitely generated ${\mathcal A}$-module, then $${\mathcal I}.{\mathcal M}={\mathcal M}\Longrightarrow\,{\mathcal M}=0$$}, implies that
$$ {\rm I}\!{\rm R}^{m-1}\phi_{*}{\mathcal O}_{\widetilde Z}=0$$
or equivalently,
$${\rm H}^{1}_{c}(Z, {\mathcal L}^{m}_{Z})=0$$
that is, ${\mathcal L}^{m}_{Z}$ has depth at least two, and thus $Z$ is of type ${\mathfrak{F.R}}\,\blacksquare$\vspace{2mm}

\noindent \noindent
{\bf(b) Case where $Z$ is pure dimensional and not necessarily normal.}
Consider the normalization $\nu:\overline{Z}\rightarrow Z$ and the commutative diagram 
$$\xymatrix{\overline{Z}\ar[rr]^{\nu}\ar[rd]_{\bar\pi}&&Z\ar[ld]^{\pi}\\
&S&}$$
As usual, we denote by ${\rm N}{\rm N}(Z)$ the non-normal locus of $Z$ and
$${\rm N}{\rm N}(\pi):=\{z\in Z: z\,{\rm is\,non\,normal\, in}\,\pi^{-1}(\pi(z))\}$$
that of $\pi$, which is an analytic subspace of $Z$ (\cite{Fi}, {\bf{Prop.3.22},p.160}). Thanks to the flatness of $\pi$, we have the relations
$${\rm N}{\rm N}(\pi)\cap Z_{s}={\rm N}{\rm N}(Z_{s})$$
and if $S$ is normal (which is our case)
$${\rm N}{\rm N}(Z)\subset {\rm N}{\rm N}(\pi),\,\,{\rm N}{\rm N}(\bar\pi)\cap {\overline{Z}}_{s}\subset \nu^{-1}({\rm N}{\rm N}(Z_{s}))\ $$
We may assume that $\bar\pi$ is a generic simultaneous normalization, meaning that ${\bar\pi}^{-1}(s)$ is the normalization of the fiber ${\pi}^{-1}(s)$ for $s\not=0$. If $\displaystyle{\overline{Z}_{0}:={\bar\pi}^{-1}(0)=\nu^{-1}({Z}_{0})}$ is the central fiber, the restriction $\displaystyle{\nu_{0}:{\overline{Z}}_{0}\rightarrow {{Z}}_{0}}$ of the finite modification $\nu$ is an analytic isomorphism from $\overline{Z}_{0}\setminus\nu^{-1}({\rm N}{\rm N}(Z_{0}))$ to $ {{Z}}_{0}\setminus {\rm N}{\rm N}(Z_{0})$, but it is not a modification if $\overline{Z}_{0}$ is not reduced. If this fiber is reduced, then $\overline{Z}_{0}$ and ${Z}_{0}$ have the same normalization (cf \cite{G}, {\bf{prop2.49}}, p.344), $\bar\pi$ is obviously flat, and we know how to conclude (cf the third remark in
{\bf(1)})\footnote{If $\overline{Z}_{0}$ is generically reduced and if the underlying reduced space is normal, then, by a classical result ([Hi], [EGA IV]{\S 5.12.3}), $\overline{Z}_{0}$ is necessarily normal and $\bar\pi$ flat, which means that $\pi$ is equinormalizable.}.\vspace{1mm}

\noindent 
Again, point {\bf(ii)} of (\cite{G}, {\bf{prop2.49}}, p.344) allows us to see, knowing that $Z$ is of pure dimension and that the dimension formula still holds for $\bar\pi$ which is geometrically flat since equidimensional over a normal base (in particular open), that $\nu^{-1}({\rm N}{\rm N}(Z_{0}))$ has empty interior in $\overline{Z}_{0}$. It follows that $\nu_{0}$ is a finite modification and that, consequently, ${\rm Prof}({\mathcal L}^{m-1}_{\overline{Z}_{0}})={\rm Prof}({\mathcal L}^{m-1}_{{Z}_{0}})\geq 2$, meaning $\overline{Z}_{0}$ is of type ${\mathfrak{F.R}}$. We then deduce from {\bf(i)} that $\overline{Z}$ is of type ${\mathfrak{F.R}}^{*}$ and, consequently, ${Z}$ is of type ${\mathfrak{F.R}}\,\blacksquare$\vspace{2mm}

\centerline{{\bf{(2) $\pi:Z\rightarrow S$ reduced of type ${\mathfrak{F.R}}$ with $S$ smooth of dimension $r>1$ and the general case:}}}
 \vspace{1mm}  
 
 \noindent
 Recall that the morphism $\pi:Z\rightarrow S$ is flat since it is open with reduced fibers over a smooth base. We proceed by induction on the dimension $r$ of $S$ as in the proof of (\theoremref{Th2}, {\bf{(a) ii)}}, p.28), assuming the result holds for any base of dimension strictly less than $r$. Using the notations of this paragraph, we have the Cartesian diagram of base change 
$$\xymatrix{Z'\ar@/_2pc/[dd]_{\psi'}\ar[r]^{\sigma'}\ar[d]_{\pi'}&Z\ar[d]^{\pi}\ar@/^2pc/[dd]^{\psi}\\
 S'\ar[d]\ar@{^{(}->}[r]_{\sigma}&S\ar[d]\\
 \{0\}\ar@{^{(}->}[r]&S"}$$
 in which $S'$ is smooth of dimension $1$, $S"$ appropriately chosen of dimension $r-1$. We see that the fiber over $\{0\}$ of $\psi$ is given by the analytic space $Z\times_{S"} \{0\}$ analytically isomorphic to $Z':=Z\times_{S} S'$ which has singularities of type ${\mathfrak{F.R}}$ since the morphism $\pi'$ is flat with reduced fibers over a smooth base of dimension $1$. We then conclude by applying the induction hypothesis to $\psi$ to deduce that $Z$ is of type ${\mathfrak{F.R}}$.\vspace{1mm}
 
 \noindent The case of singular base is deduced from the smooth base by applying the \propositionref{P7}$\blacksquare$ 
 \subsection{$\pi:Z\rightarrow S$ universally $n$-equidimensional with $Z$ and $S$ reduced and no reduced fibers.}\vspace{1mm}

\noindent 
 An intermediate step toward the case of a morphism with non-reduced fibers is suggested by the case of a flat morphism whose generic fiber is reduced. In this situation and given our reductions (given at the beginning of the paragraph), we may assume our source $Z$  reduced by virtue of (\cite{L}, {\bf{theorem 2.5}}, p.645) and proceed as above without significant modifications.\vspace{1mm}

 \noindent 
In general, one can observe that the incidence condition $(P)$ in \theoremref{Th3}, which requires the intersection ${\rm Sing}(\pi)\cap Z_{s}$ to have empty interior in $Z_s$ for all $s\in S$, is not at all an automatic condition since the analytic closed subset ${\rm Sing}(\pi)$ as well as the non-flatness locus of $\pi$ is generally not of empty interior (meaning that ${\rm Reg}(\pi)$ is not a dense open subset!). However, this is the case if $\pi$ is flat (cf. \cite{G}, p103--108).\vspace{1mm}

\noindent We will begin by noting that this condition $(P)$ then imposes on the non-normal locus of $Z$ to satisfy the same condition:
$$(P')\,\,\,\,\,\widering{\pi^{-1}(s)\cap{\rm N}{\rm N}(Z)}=\emptyset\,\,{\rm in}\,\pi^{-1}(s),\,\,\forall\,s\in S$$
Note that, since $Z$ and $S$ are reduced and of pure dimension, these incidence relations are generically satisfied on $S$, meaning that there exists a dense open subset of $S$ over which all fibers satisfy $(P)$ or $(P')$.
\vspace{1mm}

\noindent Let us verify that $(P)\,\Longrightarrow\,(P')$.\vspace{1mm}

\noindent Since every regular point $z$ for $\pi$ is, a fortiori, a point of flatness and normality on the fiber $\pi^{-1}(\pi(z))$, we have:
$${\rm Reg}(\pi)\subset {\rm N}(\pi)$$
On the other hand, if $z$ is a normal point of the fiber $Z_{\pi(z)}$, then it is necessarily normal in $Z$ because, at this point, the fiber is reduced and the flat morphism implies the normality of $Z$ at $z$. Thus, ${\rm N}(\pi)\subset {\rm N}(Z)$ and, therefore:
$${\rm N}{\rm N}(Z)\cap Z_{s}\subset {\rm N}{\rm N}(\pi)\subset {\rm Sing}(\pi)$$ 
\vspace{1mm}

\noindent 
Indeed, consider the base change diagram:
$$\xymatrix{\theta^{-1}({\rm Reg}(\pi))\ar[dd]_{\bar\pi'}\ar@{_{(}->}[rd]\ar[rrr]^{\theta'}&&&{\rm Reg}(\pi)\ar@{^{(}->}[ld]\ar[dd]^{\pi'}\\
&\overline{Z}\ar[r]^{\theta}\ar[dd]_{\bar\pi}&Z\ar[dd]^{\pi}&\\
\nu^{-1}(S')\ar@{_{(}->}[rd]\ar[rrr]^{\nu'}&&&{S':=\pi({\rm Reg}(\pi))}\ar@{^{(}->}[ld]\\
&\overline{S}\ar[r]_{\nu}&S&}$$
relative to a proper modification with singular locus as center $\nu:\overline{S}\rightarrow S$ and the inclusion ${\rm Reg}(\pi)\into Z$.\vspace{1mm}

\noindent 
Since smooth morphisms are stable under arbitrary base changes, $\bar\pi'$ is smooth because $\pi'$ is smooth, and we have $\theta^{-1}({\rm Reg}(\pi))\subset {\rm Reg}(\bar\pi)$. Then, since $\pi$ is universally $n$-equidimensional, its restriction to the dense open subset ${\rm Reg}(\pi)$ remains a universally $n$-equidimensional morphism and even flat by definition. It follows that $S':=\pi({\rm Reg}(\pi))$ is a dense open subset of $S$ over which the fibers are smooth and hence itself smooth. Consequently, the restriction $\nu'$ (resp. $\theta'$) of the modification $\nu$ (resp. $\theta$) to $\nu^{-1}(S')$ (resp. $\theta^{-1}({\rm Reg}(\pi))$) is an analytic isomorphism.\vspace{1mm}

\noindent
Now, since $\theta$ is proper, the inclusion ${\rm Sing}(\bar\pi)\cap \overline{Z}_{\bar s}\subset \theta^{-1}({\rm Sing}(\pi))\cap \overline{Z}_{\bar s}$ yields the inclusions:
$$\theta({\rm Sing}(\bar\pi)\cap \overline{Z}_{\bar s})\subseteq {\rm Sing}(\pi)\cap \theta(\overline{Z}_{\bar s})={\rm Sing}(\pi)\cap Z_{s}\subsetneq Z_s$$
which shows that $(P)$ is satisfied for $\bar\pi$ since $\overline{Z}_{\bar s}$ is analytically isomorphic to $Z_s$ and thus has the same irreducible components (which prevents ${\rm Sing}(\bar\pi)$ from having a common irreducible component with $\overline{Z}_{\bar s}$).\vspace{1mm}

\noindent More generally, if $F$ is a non-$\pi$-saturated closed subset—meaning that $F\cap \pi^{-1}(\pi(x))$ has no irreducible component in common with $\pi^{-1}(\pi(x))$ for all $x\in F$—then the closed subset $\theta^{-1}(F)$ is also non-$\bar\pi$-saturated. This is easily understood knowing that:
$$\theta^{-1}(F)\cap \overline {Z}_{\bar s}\simeq (F\cap Z_{s})\times_{\overline{S}} S $$
\vspace{2mm}

\noindent 
Thanks to the \propositionref{P8} and the remarks above, we proceed exactly as in \S5.1 p.27 and S6.1 p.34, assuming $S$ and begin by treating: \vspace{1mm}

\indent
{\bf{(a) the case of a 1-dimensional base and $Z$ normal.}} \vspace{1mm}

\noindent\par 
With the same notations and considerations as in \S5.1, we have the short exact sequence
$$\xymatrix{0\ar[r]&{\tilde{f}}.{\mathcal O}_{Z}\ar[r]^{\bullet {\tilde{f}}}&{\mathcal O}_{Z}\ar[r]&{{ i_{0}}}_{*}{\mathcal O}_{Z_{0}}\ar[r]&0}$$
to which we apply the functor ${\mathcal H}om(-, {\mathcal L}^{m}_{Z})$ giving us the exact sequence
$$\xymatrix{0\ar[r]&{\mathcal L}^{m}_{ Z}\ar[r]&{\mathcal L}^{m}_{ Z, {\tilde{f}}}\ar[r]&{\mathcal E}xt^{1}({{ i_0}}_{*}{\mathcal O}_{{Z_0}}, {\mathcal L}^{m}_{Z})\ar[r]&0}$$
with ${\mathcal L}^{m}_{Z, {\tilde{f}}}:={\mathcal H}om({\tilde{f}}.{\mathcal O}_{Z}, {\mathcal L}^{m}_{Z})$.\vspace{1mm}

\noindent Just as in (\S 6.1, {\bf(1).a}), we can easily see that the sheaf ${\mathcal E}xt^{1}({i_{0}}_{*}{\mathcal O}_{Z_0}, {\mathcal L}^{m}_{Z})$ is torsion-free on $Z_0$ since isomorphic to the torsion-free sheaf ${{i_{0}}_{*}}({\phi_{0}}_{*}(\omega^{m-1}_{{\widetilde Z}_{0}}))$. \vspace{1mm}

\noindent While the morphism $${i_{0}}_{*}{\mathcal L}^{m-1}_{Z_0}\rightarrow {\mathcal E}xt^{1}({i_{0}}_{*}{\mathcal O}_{Z_0}, {\mathcal L}^{m}_{Z})$$ 
is always injective, the same is not true for
$${\mathcal E}xt^{1}({i_{0}}_{*}{\mathcal O}_{Z_0}, {\mathcal L}^{m}_{Z})\rightarrow {i_{0}}_{*}{\mathcal L}^{m-1}_{Z_0}$$
whose injectivity requires the existence of a dense open subset of $Z$ whose trace on $Z_0$ is still dense and on which this morphism is bijective.
\vspace{1mm}

\noindent 
Since, by hypothesis, ${\rm Sing}(Z)\cap Z_0$ has empty interior in $Z_0$,  we can then work generically on $Z$ and $Z_0$ and establish the desired injectivity and consequently the isomorphism
$${i_{0}}_{*}{\mathcal L}^{m-1}_{Z_0}\simeq {\mathcal E}xt^{1}({i_{0}}_{*}{\mathcal O}_{Z_0}, {\mathcal L}^{m}_{Z})$$
The sheaf ${\mathcal L}^{m-1}_{Z_0}$ being, by hypothesis, of depth at least two, the above short exact sequence and the vanishing $\displaystyle{{\rm H}^{0}_{c}(Z, {\mathcal L}^{m-1}_{Z_0})\simeq {\rm H}^{1}_{c}(Z, {\mathcal L}^{m-1}_{ Z_0})}$ give us the isomorphism 
$${\rm H}^{1}_{c}(Z, {\mathcal L}^{m}_{Z})\simeq {\rm H}^{1}_{c}(Z, {\mathcal L}^{m}_{ Z, {\tilde{f}}}),$$
To show that this necessarily implies ${\rm H}^{1}_{c}(Z, {\mathcal L}^{m}_{Z})=0$, we proceed as in the case of the reduced morphism where we had obtained
$${\rm H}^{1}_{c}(Z, {\mathcal L}^{m}_{Z})\simeq {\rm H}^{1}_{c}(\widetilde{Z}, {\Omega}^{m}_{\widetilde Z})\simeq \biggl({\rm H}^{m-1}(\widetilde{Z}, {\mathcal O}_{\widetilde Z})\biggr)^{'}\simeq \biggl(\Gamma(Z, {\rm I}\!{\rm R}^{m-1}\phi_{*}{\mathcal O}_{\widetilde Z})\biggr)^{'} $$
Then, denoting by ${\mathcal I}$ (resp. $\tilde{\mathcal I}$) the ideal generated by $f$ (resp. $\phi^{*}(f)$), we have
$${\rm H}^{1}_{c}(Z, {\mathcal L}^{m}_{ Z, {\tilde{f}}})={\rm H}^{1}_{c}(Z, {\mathcal H}om({\mathcal{I}}{\mathcal O}_{Z}, {\mathcal L}^{m}_{Z}))\simeq {\rm H}^{1}_{c}(Z, \phi_{*}{\mathcal H}om({\tilde{\mathcal I}}{\mathcal O}_{\widetilde Z}, {\Omega}^{m}_{\widetilde Z}))$$
Since $${\rm I}\!{\rm R}^{i}\phi_{*}(({\tilde{\mathcal I}}{\mathcal O}_{\widetilde Z})^{*}\otimes {\Omega}^{m}_{\widetilde Z}))=,\,\forall\,i\not=0$$
we obtain
$${\rm H}^{1}_{c}(Z, {\mathcal L}^{m}_{ Z, {\tilde{f}}})\simeq {\rm H}^{1}_{c}(\widetilde{Z}, {\tilde{\mathcal I}}{\mathcal O}_{\widetilde Z})^{*}\otimes {\Omega}^{m}_{\widetilde Z})\simeq \biggl({\rm H}^{m-1}(\widetilde{Z}, \tilde{\mathcal I}{\mathcal O}_{\widetilde Z})\biggr)^{'}\simeq \biggl(\Gamma(Z, {\mathcal I}\otimes{\rm I}\!{\rm R}^{m-1}\phi_{*}{\mathcal O}_{\widetilde Z})\biggr)^{'} $$
 The initial isomorphism is therefore equivalent to the isomorphism between coherent sheaves (of finite type) induced by multiplication by ${\mathcal I}$
 $$ {\rm I}\!{\rm R}^{m-1}\phi_{*}{\mathcal O}_{\widetilde Z}\simeq  {\mathcal I}\otimes{\rm I}\!{\rm R}^{m-1}\phi_{*}{\mathcal O}_{\widetilde Z}$$
We then conclude by applying Nakayama's lemma as used above, concluding with the vanishing of ${\rm I}\!{\rm R}^{m-1}\phi_{*}{\mathcal O}_{\widetilde Z}$ and thus that of ${\rm H}^{1}_{c}(Z, {\mathcal L}^{m}_{ Z})$.\vspace{1mm}

\indent {\bf(b) The case of a base of dimension $1$ with $Z$ reduced and pure-dimensional. } \vspace{1mm}

\noindent 

We can revisit the arguments developed in (\S6.1, {\bf 1.b)}, p.38) by again considering the commutative diagram 
$$\xymatrix{\overline{Z}\ar[rr]^{\nu}\ar[rd]_{\bar\pi}&&Z\ar[ld]^{\pi}\\
&S&}$$
where $\nu$ is the normalization morphism. We note that if $z$ is a normal point of the fiber $Z_{\pi(z)}$, then it is necessarily normal in $Z$ since, at this point, the fiber is reduced and, consequently, the morphism is flat at $z$. In this case, $Z$ is normal at $z$. Thus, ${\rm N}{\rm N}(Z)\cap Z_{s}\subset {\rm N}{\rm N}(Z_s)$ and, therefore, the finite morphism induced by $\nu$ on $\overline{Z_s}$ is a biholomorphism from $\overline{Z_s}\setminus \nu^{-1}({\rm N}{\rm N}(Z)\cap Z_{s})$ to ${Z_s}\setminus {\rm N}{\rm N}(Z)\cap Z_{s}$. Since, by hypothesis, ${\rm N}{\rm N}(Z)\cap Z_s$ has empty interior in $Z_s$, we conclude in the same way as in (\S6.1, {\bf(b)}, p.38) by showing that $\nu_{s}:\overline{Z_s}\rightarrow Z_s $ is a modification. We then deduce, thanks to case {\bf(a)}, that $Z$ is of type ${\mathfrak{F.R}}\,\blacksquare$.\vspace{1mm}

\noindent 
{\bf(c) The case of a general reduced base.} As in (\S6.1, {\bf(2)}, p.38), we extend to the case of a smooth base of dimension greater than $1$ by noting that the condition $(P)$ is stable by base change between smooth bases. Indeed if $S$ is a smooth open subset and $S'$ a smooth subvariety equipped with the embedding $S'\into S$, and possibly after restricting the data, we can assume $S$  equipped with a system of parameters $(s_1,s_2,\cdots,s_r)$ and $S'$ a hypersurface given by the vanishing of a certain coordinate and then see that this condition is obviously satisfied over $S'$.  
hypersurface given by the vanishing of a certain coordinate.\vspace{1mm}

\noindent 
The case of a reduced base of type ${\mathfrak{F.R}}$ is deduced from the smooth base thanks to the \propositionref{P8}$\,\blacksquare$ \vspace{1mm}

\noindent
\begin{rem} 
{\bf(i)} Instead of ${\rm Reg}(\pi)$, one can take the Cohen-Macaulay locus ${\rm CM}(\pi)$, which is a dense open subset of $Z$ by virtue of Frisch's theorem (\cite{Fri}) and the fact that ${\rm Reg}(\pi)\subset {\rm CM}(\pi)$.\vspace{1mm}

\noindent {\bf(ii)} If $\pi$ is geometrically flat, its restriction to the dense open subset ${\rm Reg}(\pi)$ remains geometrically flat, and the dense open subset $\pi({\rm Reg}(\pi))$ has rational singularities by virtue of \theoremref{Th1}.
\end{rem}
\noindent \section{\color{blue}{Transfer of properties $\mathfrak{F.R}^{*}$ and ${\mathfrak{F.R}}$ from the total space to the fibers: \theoremref{Th4}.}}
 Let $Z$ and $S$ be two reduced complex spaces of pure dimensions $m$ and $r$, respectively. Let $\pi:Z\rightarrow S$ be an open morphism with fibers of pure dimension $n$ and reduced. We assume that $\pi$ admits a simultaneous resolution, and that $Z$ and $S$ are of type $\mathfrak{F.R}^{*}$ (resp. $\mathfrak{F.R}$). Let us show that $\pi$ is of type $\mathfrak{F.R}^{*}$ (resp. $\mathfrak{F.R}$).\rm\vspace{1mm}

 \noindent 
First, note that we may assume $S$ is smooth. Indeed, let $\phi:\widetilde{Z}\rightarrow Z$ be a simultaneous resolution of singularities for $\pi$, and let $\nu:\widetilde{S}\rightarrow S$ be a desingularization of $S$. Then, we have the commutative diagram of base change 
$$\xymatrix{\widetilde{{S}}\times_{{S}}\widetilde{Z}\ar[d]_{\tilde{\phi}}\ar[r]_{\theta}\ar@/_2pc/[dd]_{\tilde\Psi}&\widetilde{Z}\ar[d]^{\phi}\ar@/^2pc/[dd]^{\Psi}\\
\widetilde{{S}}\times_{{S}}{Z}\ar[d]_{\tilde{\pi}}\ar[r]_{\tilde\nu}
&Z\ar[d]^{\pi}\\
\widetilde{{S}}\ar[r]_{\nu}&{S}}$$ 
Since smooth morphisms are stable under arbitrary base change, it follows that $\tilde\Psi$ is a smooth morphism and, consequently, $\tilde\phi$ is a simultaneous resolution for $\tilde\pi$. Thus, thanks to \lemmaref{L3} or \lemmaref{L2} and \propositionref{P8}, we see that if the result holds for $\tilde\pi$, it will also hold for $\pi$, whose fibers are analytically isomorphic to those of $\tilde\pi$.\vspace{1mm}

\noindent {\bf(i) Case $\mathfrak{F.R}^{*}$:} This is an immediate consequence of (\cite{AS84}, {\bf{lemma 5.8},p.74}), which gives the existence of a commutative diagram 
$$\xymatrix{{\mathcal L}^{m}_{Z}\ar[r]\ar[d]&{\mathcal L}^{n}_{Z_s}\ar[r]\ar[d]&0\\
{\omega}^{m}_{Z}\ar[r]&{\omega}^{n}_{Z_s}\ar[r]&0}$$
\vspace{1mm}

\noindent {\bf(ii) Case $\mathfrak{F.R}$:} We factorize $\phi$ through the normalization $\nu:\overline{Z}\rightarrow Z$ to obtain the commutative diagram 
$$\xymatrix{&\widetilde{Z}\ar[ld]_{\tilde{\phi}}\ar[rd]^{\phi}\ar@/_5pc/[dd]_{\tilde\Psi}\ar@/^5pc/[dd]^{\Psi}&\\
\overline{Z}\ar[rd]_{\bar{\pi}}\ar[rr]^{\nu}
&&Z\ar[ld]^{\pi}\\
&{S}&}$$ 
The morphism $\pi$ is flat because it is open with reduced fibers over a smooth base, the morphism $\Psi$, by hypothesis, is smooth, and the fibers ${\widetilde Z}_{s}$ are desingularizations of the fibers ${Z}_{s}$. The morphism $\bar\pi$ is also open with fibers of pure dimension $n$ because it is equidimensional over a smooth base and thus flat. Moreover, these fibers are necessarily reduced since their preimages by $\tilde\phi$ are smooth and hence reduced. It follows that $\tilde\phi$ is a simultaneous resolution of $\bar\pi$.\vspace{1mm}

\noindent Then, we have 
$$\xymatrix{{Z\,{\rm of\, type}\, \mathfrak{F.R}}\ar@{<=>}[r]&{\overline{Z}\,{\rm of\, type\,}\mathfrak{F.R}^{*}}\ar@{=>}[r]&\overline{Z}_{s}\,{\rm of\, type}\, \mathfrak{F.R}^{*}\ar@{=>}[r]&{{Z}_{s}{\rm \,of\, type}\, \mathfrak{F.R}}}$$ 
where the first equivalence is given by \propositionref{P3}, the second implication is deduced from case {\bf(i)}, and the last is due to the fact that the morphism induced by $\nu$ on each fiber $\overline{Z}_{s}$ is a finite modification dominated, moreover, by the normalization of ${Z}_{s}$ (which is also that of $\overline{Z}_{s}\,\blacksquare$ 
\section{\color{blue}{General remarks and some examples.}}

{\bf(i)} If $Z_0$ is an hypersurface of the Cohen Macaulay reduced complex spaces $Z$, the hypothesis ${\rm Prof}({\mathcal L}^{m-1}_{Z_{0}}):={\rm Prof}({\mathcal L}^{m-1}_{Z}\vert_{Z_{0}})\geq 2$ implies that ${\rm Prof}({\mathcal L}^{m-1}_{Z})\geq 3$ due to the well-known fact (in algebraic or analytic geometry) that if $Z_0$ is a Cartier divisor of a reduced analytic pure dimensional complex space $Z$  and ${\mathcal F}$ is an $S_1$-coherent sheaf on $Z$, then, at a point $z$,  
$${\rm Prof}({\mathcal F}_{z})\geq d+1\,\Longleftrightarrow\,{\rm Prof}(({\mathcal F}\vert_{Z_{0}})_{z})\geq d$$  
Thus, for $Z$  normal, we have  
$${\mathcal L}^{m-1}_{Z}\simeq \omega^{m-1}_{Z}$$  
which, by virtue of (\cite{KeSc}, {\bf{theorem (1.2), theorem 4}}), implies the isomorphisms  
$${\mathcal L}^{j}_{Z}\simeq \omega^{j}_{Z},\,\forall\,j\in\{0,1,\cdots, m-1\}$$  
enforcing the inclusions \footnote{These equalities are not sufficient to enforce the equality ${\mathcal L}^{m}_{Z}\simeq \omega^{m}_{Z}$, as can be seen by revisiting the two-parameter family from (\cite{SS85}, p. 104) given by ${\mathcal Z}:=\{(x,y,z,s,t)\in{\Bbb C}^{3}\times {\Bbb C}^{2}: \frac{x^{7}}{7} + \frac{y^{3}}{3}+\frac{z^{3}}{3}-syx^{5} -tyzx^{4}=0\}$, for which the arithmetic genus is zero for $s=0$ and $t\not=0$ and the geometric genus is non-zero for every pair $(s,t)$; this translates into the equalities ${\mathcal L}^{j}_{Z}\simeq \omega^{j}_{Z}$ for $j=0,1$ and the strict inclusion ${\mathcal L}^{2}_{Z}\subset \omega^{2}_{Z}$.}  
$${\mathcal L}^{j}_{Z}\subset{\mathcal H}om({\mathcal L}^{m-j}_{Z}, {\mathcal L}^{m}_{Z})\subset{\mathcal H}om({\Omega}^{m-j}_{Z}, {\mathcal L}^{m}_{Z}) \subset {\mathcal H}om({\Omega}^{m-j}_{Z}, {\omega}^{m}_{Z})\simeq \omega^{j}_{Z}$$  
to be bijective for all $j\in\{0,1,\cdots, m-1\}$. \vspace{1mm}  

\noindent
{\bf(ii)} The arrow ${\mathcal E}xt^{1}({\mathcal O}_{Z_0}, {\mathcal L}^{m}_{Z})\rightarrow{i_{0}}_{*}{\mathcal L}^{m-1}_{Z_0}$  
If one has a residue morphism ${\rm Res}_{Z_0}(\omega):{\mathcal E}xt^{1}({\mathcal O}_{Z_0}, {\omega}^{m}_{Z})\rightarrow {i_{0}}_{*}{\omega}^{m-1}_{Z_0}$, the arrow ${\mathcal E}xt^{1}({\mathcal O}_{Z_0}, {\mathcal L}^{m}_{Z})\rightarrow{i_{0}}_{*}{\mathcal L}^{m-1}_{Z_0}$ can only be the restriction of this residue morphism, especially since, for $Z$ and $Z_0$ being Gorenstein, it is easy to construct the commutative diagram  
$$\xymatrix{0\ar[r]&{\mathcal L}^{m}_{Z}\ar@{^{(}->}[d]\ar[r]&{\mathcal L}^{m}_{Z}[Z_0]\ar@{^{(}->}[d]\ar[r]& {\mathcal L}^{m-1}_{Z_0}\ar@{^{(}->}[d]\ar[r]&0\\  
0\ar[r]&{\omega}^{m}_{Z}\ar[r]&{\omega}^{m}_{Z}[Z_0]\ar[r]& {\omega}^{m-1}_{Z_0}\ar[r]&0}$$  
One can make the calculations explicit by considering a local embedding of $Z$ into a Stein manifold $W$ and an embedded desingularization. By fixing a generator $\omega_{0}$ of $\omega^{m}_{Z}$ that allows describing all sections of the canonical sheaf in the form $g.\omega_{0}$ with $g$ an arbitrary holomorphic function on $Z$. Then, it is known that the sections of ${\mathcal L}^{m}_{Z}$ are exactly the forms $g.\omega_{0}$ with $g$ in the adjunction ideal of $Z$, i.e., the holomorphic functions $g$ such that $g.\omega_{0}$ is square integrable. It is clear that for any element of ${\mathcal E}xt^{1}({\mathcal O}_{Z_0}, {\mathcal L}^{m}_{Z})$, the relation  
$${\rm Res}_{Z_{0}}\left(\frac{g\omega_{0}}{f}\right)=g{\rm Res}_{Z_{0}}\left(\frac{\omega_{0}}{f}\right)$$  
shows that ${\rm Res}_{Z_{0}}\left(\frac{g\omega_{0}}{f}\right)$ is a section of the sheaf $\omega^{m-1}_{Z_0}$\footnote{For a local embedding $Z_{0}\into Z\into W$ such that $Z:=\{z:=(z_1,\cdots, z_{m+1})/ g(z)=0\}$ and $Z_{0}:=\{z:=(z_1,\cdots, z_{m+1})\in W/ g(z)=f(z)=0\}$  
we can  use the short exact sequence  
$$\xymatrix{0\ar[r]&{\Omega}^{m+1}_{W}\ar[r]&{\rm Adj}(Z){\Omega}^{m+1}_{W}[Z]\ar[r]& {\mathcal L}^{m-1}_{Z}\ar[r]&0}$$}S and as already mentioned, we can assume the situation that the singular locus of $Z$ is the singular locus of $Z_0$. With this assumption, it's obvious that the restriction to $Z_0$ of the holomorphic function $g$ still in the adjunction ideal of $Z_0$  which precisely means that this residue defines a section of the sheaf ${\mathcal L}^{m-1}_{Z_0}$. Since every element of ${\mathcal E}xt^{1}({\mathcal O}_{Z_0}, {\mathcal L}^{m}_{Z})$ can be described in this way, the conclusion follows.  
\vspace{1mm}

\noindent  {\bf(iii)} Some examples of particular morphisms:
Consider the hypersurfaces in ${\Bbb C}^{4}$ given respectively by
$Z_{1}:=\{(x,y,z,t)\in {\Bbb C}^{4}: x^{2} + zy^{2} +txz=0\}$, $Z_{2}:=\{(x,y,z,t)\in {\Bbb C}^{4}: x^{2} + zy^{3} +tz^{2}=0\}$, $Z_{3}:=\{(x,y,z,t)\in {\Bbb C}^{4}: z^{2} + xy^{2} +tx^{2}z=0\}$, and 
$Z_{5}:=\{(x,y,z,t)\in {\Bbb C}^{4}: x^{4} + y^{3} +tz^{3}=0\}$  
the morphisms $\pi_{j}:Z_{j}\rightarrow {\Bbb C}$ induced  by the canonical projection ${\Bbb C}^{4}\rightarrow {\Bbb C}$ sending $(x,y,z,t)$ to $t$.\vspace{1mm}

\noindent 
Explicit calculations show that these Gorenstein morphisms are, for the first three, generically normal fibers of type ${\mathfrak{F.R}}^*$ (hence with rational singularities) with a non-normal central fiber of type $\frak{FR}$ and a normal total space of type ${\mathfrak{F.R}}^*$ (thus also with rational singularities since normal); for the fourth, generically normal fibers of type ${\mathfrak{F.R}}^*$, a central fiber that is normal but not of type $\frak{FR}$, and a total space of type ${\mathfrak{F.R}}^*$.\vspace{1mm}

\noindent The fifth has generically normal fibers not of type $\frak{FR}$, a non-normal central fiber of type $\frak{FR}$, and a normal total space not of type $\frak{FR}$.

\end{document}